\documentclass[12pt,reqno]{amsart}%{article}
\usepackage{amsfonts,amsmath,amssymb}
\usepackage{tikz}
\usetikzlibrary{decorations.pathmorphing}

\usepackage{graphicx}

\def\bim{\begin{itemize}\item[]}
	
	\def\eim{\end{itemize}}

\allowdisplaybreaks
\advance\textheight by \topskip

\usepackage{mathrsfs}

\usepackage[latin1]{inputenc}
\usepackage{graphicx}

\newtheorem{theorem}{Theorem}
\usepackage{euscript}

\def\[{[\! [}
\def\]{]\! ]}

\title{A walk in my lattice path garden}

\author[H.~Prodinger]{Helmut Prodinger}

\address{Helmut Prodinger,
Mathematics Department, Stellenbosch University,
7602 Stellenbosch, South Africa.}
\email{hproding@sun.ac.za}

\dedicatory{Dedicated to Neil Peart.}

\date{\today}

\begin{document}

\begin{figure}
	\includegraphics[width=10cm]{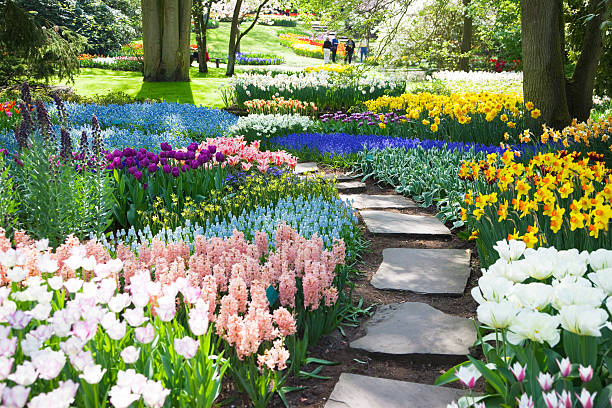}
	%	\caption{A boat.}
	%	\label{fig:boat1}
\end{figure}

\begin{abstract}
Various lattice path models are reviewed. The enumeration is done using generating functions. A few
bijective considerations are woven in as well. The kernel method is often used. Computer algebra was an essential tool.
Some results are new, some have appeared before. 

The lattice path models we treated, are: Hoppy's walks, the combinatorics of sequence A002212 in \cite{OEIS} 
(skew Dyck paths,  Schr\"oder paths, Hex-trees, decorated ordered trees, multi-edge trees, etc.) Weighted unary-binary trees also
occur, and we could improve on our old paper on Horton-Strahler numbers \cite{FlPr86}, by using a different substitution.
Some material on ternary trees appears as well, as on Motzkin numbers and paths (a model due to Retakh), and a new concept
called amplitude that was found in \cite{irene}. Some new results on Deutsch paths in a strip are included as well. During the Covid period,
I spent much time with this beautiful concept that I dare to call Deutsch paths, since Emeric Deutsch stands at the beginning with a
problem that he posted in the American Mathematical Monthly some 20 years ago. Peaks and valleys, studied by Rainer Kemp 40 years under the names \textsc{max}-turns and \textsc{min}-turns,
are revisited with a more modern approach, streamlining the analysis, relying on the `subcritical case' (named so by Philippe Flajolet),
the adding a new slice technique and once again the kernel method.
\end{abstract}

\maketitle

\newpage

\setcounter{tocdepth}{1}
\tableofcontents
%\newpage

\section{Introduction}

Around 20 years ago I published a collection of examples about applications of the kernel method
in the present journal. The success of this enterprise was unexpected and came as a very pleasant 
surprise. My current plan is to present again a collection of subjects, loosely related as they all
have a lattice path flavour (trees are also allowed in my private book when lattice paths are mentioned).
The subjects cover my last 2 or 3 years of research; some results were only posted on the arxiv, some
were submitted but no evaluations were available up to this day, some have been appearing, and some are
completely new. In each instance I try to make clear the current status. 

As in the predecessor paper, the kernel method plays a role again, but also analytic techniques like
singularity analysis and Mellin transform, as well as bijective results.

	\section{Hoppy walks}

Deng and Mansour \cite{deng} introduce a rabbit named Hoppy and let him move according to certain rules.
While the story about Hobby is charming and entertaining, we do not need this here and move straight ahead
to the enumeration issues. Eventually, the enumeration problem is one about $k$-Dyck paths.
The up-steps are $(1,k)$ and the down-steps are $(1,-1)$.
The model that has $(1,1)$ as up-step and the down-step are $(1,-k)$ will also be called $k$-Dyck paths.

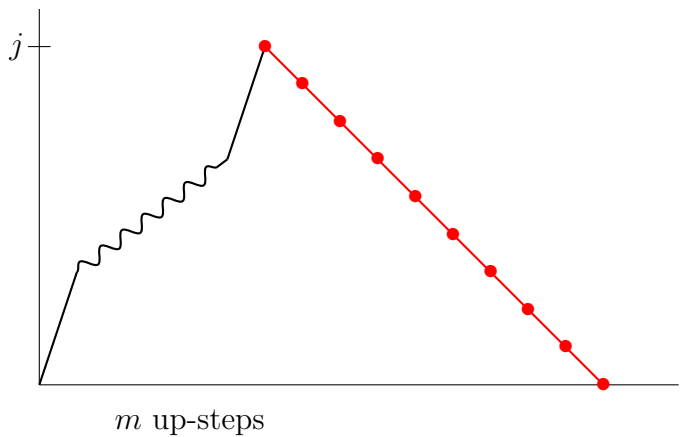
\begin{figure}[h]
\begin{center}
	\begin{tikzpicture}[scale=0.5]
		\draw (0,0) -- (17,0);
		\draw (0,0) -- (0,10);
		\draw[thick](0,0)--(1,3);
		\draw[thick](5,6)--(6,9);
		\draw [decorate,decoration=snake,thick] (1,3) -- (5,6);
		\foreach \i in {0,...,8}
		{\draw[thick,red](6+\i,9-\i)--(7+\i,8-\i);
			\node[thick, red] at (6+\i,9-\i){$\bullet$};
		}
		\node[thick, red] at (6+5+4,9-5-4){$\bullet$};

		\node at (4,-1){$m$ up-steps};
		
		%\draw (-0.3,6) --(0.3,6);
		%\draw (-0.3,9) --(0.3,9);
		\node[thick] at (-1+0.4,9){$j$};

		\draw (-0.3,9) --(0.3,9);
	\end{tikzpicture}
\end{center}
\caption{The number of final down-steps}
\end{figure}	
The question is about the length of the sequence of down-steps printed in red. Or, phrased differently, how many $k$-Dyck paths end on level $j$, after $m$ up-steps, the last step being an up-step. The recent paper \cite{jcmcc} contains similar computations, although without the restriction that the last step must be an up-step. 

Counting the number of up-steps is enough, since in total, there are $m+km=(k+1)m$ steps. The original description of Deng and Mansour is a reflection of this picture, with up-steps of size 1 and down-steps of sice $-k$, but we prefer it as given here, since we are going to use the adding-a-new-slice method, see \cite{FP, Prodinger-handbook}. A slice is here a run of down-steps, followed by an up-step. 
The first up-step is treated separately, and then $m-1$ new slices are added. We keep track of the level after each slice, using a variable $u$. The variable $z$ is used to count the number of up-steps.

Deng and Mansour work out a formula which comprises $O(m)$ terms. Our method leads only to a sum of $O(j)$ terms.

We will use the adding a new slice technique \cite{Prodinger-handbook} together with the kernel method to find useful bivariate generating functions.

The following substitution is essential for adding a new slice (a maximal sequence of down-steps followed by one up-step followed):
\begin{equation*}
	u^j\longrightarrow z\sum_{0\le h \le j} u^{h+k}=\frac{zu^k}{1-u}(1-u^{j+1}).
\end{equation*}
Now let $F_m(z,u)$ be the generating function according to $m$ runs of down-steps. The substitution leads to
\begin{equation*}
	F_{m+1}(z,u)=\frac{zu^k}{1-u}F_m(z,1)-\frac{zu^{k+1}}{1-u}F_m(z,u),\quad F_0(z,u)=zu^k.
\end{equation*}
Let $F=\sum_{m\ge0}F_m$, so that we don't care about the number $m$ anymore; then
\begin{equation*}
	F(z,u)=zu^k+\frac{zu^k}{1-u}F(z,1)-\frac{zu^{k+1}}{1-u}F(z,u),
\end{equation*}
or
\begin{equation*}
	F(z,u)\frac{1-u+zu^{k+1}}{1-u}=zu^k+\frac{zu^k}{1-u}F(z,1).
\end{equation*}	
The equation $1-u+zu^{k+1}=0$ is famous when enumerating $(k+1)$-ary trees (or $k$-Dyck paths). Its relevant combinatorial solution (also the only one being analytic at the origin) is
\begin{equation*}
	\overline{u}=\sum_{\ell\ge0}\frac1{1+\ell(k+1)}\binom{1+\ell(k+1)}{\ell}z^\ell.
\end{equation*}
Since $u-\overline{u}$ is a factor of the LHS, is must also be a factor of the RHS, and we can compute (by dividing out the factor
$(u-\overline{u})$) that
\begin{equation*}
	\frac{zu^k(1-u+F(z,1))}{u-\overline{u}}=-zu^k.
\end{equation*}
Thus
\begin{equation*}
	F(z,u)=zu^k\frac{\overline{u}-u}{1-u+zu^{k+1}}.
\end{equation*}	

%\begin{equation*}
%	F(z,u)=-1+\frac{1-u+zu^k\overline{u}}{1-u+zu^{k+1}}.
%\end{equation*}	

The first factor has even a combinatorial interpretation, as a description of the first step of the path. It is also clear from this that the level reached is $\ge k$ after each slice. We don't care about the factor $zu^k$ anymore, as it produces only a simple shift. The main interest is now how to get to the coefficients of
\begin{equation*}
	\frac{\overline{u}-u}{1-u+zu^{k+1}}
\end{equation*}	
in an efficient way. There is also the formula
\begin{equation*}
	1-u+zu^{k+1}=(\overline{u}-u)\Big(1-z\frac{u^{k+1}-\overline{u}^{k+1}}{u-\overline{u}}\Big),
\end{equation*}
but it does not seem to be useful here.

First we deal with the denominators
\begin{equation*}
	S_j:=[u^j]\frac1{1-u+zu^{k+1}}=\sum_{0\le m\le j/k}(-1)^m\binom{j-km}{m}z^m.
\end{equation*}
One way to see this formula is to prove by induction that the sums $S_j$ satisfy the recursion
\begin{equation*}
	S_j-S_{j-1}+zS_{j-k-1}=0
\end{equation*}
and initial conditions $S_0=\dots =S_k=1$. In \cite{jcmcc} such expressions also appear as determinants.
Summarizing,	
\begin{equation*}
	\frac1{1-u+zu^{k+1}}=\sum_{m\ge0}(-1)^mz^m\sum_{j\ge km }\binom{j-km}{m}u^j.
\end{equation*}
Now we read off coefficients:
\begin{equation*}
	[u^j]\frac{\overline{u}}{1-u+zu^{k+1}}=\sum_{0\le m\le j/k}(-1)^m\binom{j-km}{m}z^m
	\sum_{\ell\ge0}\frac1{1+\ell(k+1)}\binom{1+\ell(k+1)}{\ell}z^\ell
\end{equation*}
and further
\begin{multline*}
	[z^n][u^j]\frac{\overline{u}}{1-u+zu^{k+1}}\\=\sum_{0\le m\le j/k}(-1)^m\binom{j-km}{m}
	\frac1{1+(n-m)(k+1)}\binom{1+(n-m)(k+1)}{n-m}.
\end{multline*}
The final answer to the Deng-Mansour enumeration (without the shift) is
\begin{multline*}
	\sum_{0\le m\le j/k}(-1)^m\binom{j-km}{m}
	\frac1{1+(n-m)(k+1)}\binom{1+(n-m)(k+1)}{n-m}\\*
	-(-1)^n\binom{j-1-kn}{n}.
\end{multline*}
If one wants to take care of the factor $zu^k$ as well, one needs to do the replacements $n\to n+1$ and $j\to j+k$ in the formula just derived. That enumerates then the $k$-Dyck paths ending at level $j$ after $n$ up-steps, where the last step is an up-step.

\subsection*{An application}

The encyclopedia of integer sequences \cite{OEIS} has the sequences A334680, A334682, A334719, (with a reference to \cite{AHS}) which is the total number of down-steps of the last down-run, for $k=2,3,4$. So, if the path ends on level $j$, the contribution to the total is $j$. 

All we have to do here is to differentiate
\begin{equation*}
	F(z,u)=zu^k\frac{\overline{u}-u}{1-u+zu^{k+1}}.
\end{equation*}	
w.r.t.\ $u$, and then replace $u$ by 1. The result is
\begin{equation*}
	\frac{\overline{u}}z-\overline{u}-\frac1z,
\end{equation*}
and the coefficient of $z^m$ therein is
\begin{align*}
	\frac1{1+(m+1)(k+1)}\binom{1+(m+1)(k+1)}{m+1}-\frac1{1+m(k+1)}\binom{1+m(k+1)}{m}.
\end{align*}
The bivariate generating function does this enumeration cleanly and quickly.

\subsection*{Hoppy's early adventures}

Now we investigate what Hoppy does after his first up-step; he might follow with $0,1,\dots,k$ down-steps. Eventually, we want to sum all these steps (red in the picture).
\begin{center}
	\begin{tikzpicture}[scale=0.5]
		\draw (0,0) -- (17,0);
		\draw (0,0) -- (0,9);
		\draw[thick](0,0)--(1,7);
		\foreach \i in {0,...,4}
		{\draw[thick,red](1+\i,7-\i)--(2+\i,6-\i);
			\node[thick, red] at (1+\i,7-\i){$\bullet$};
		}
		\node[thick, red] at (6,2){$\bullet$};
		\draw[thick](6,2)--(7,9);
		\draw [decorate,decoration=snake,thick] (7,9) -- (17,0);
		
		\draw(-0.2,2)--(0.2,2);
		\node[thick,red] at (-2,2){$k-i=h$};
		\draw[ dashed](6,0)--(6,9);
		\node[ ] at (3,-0.8){one up-step};
		\node[ ] at (11,-0.8){$m$ up-steps};
	\end{tikzpicture}
\end{center}

A new slice is now an up-step, followed by a sequence of down-steps. The substitution of interest is:
\begin{equation*}
	u^i\rightarrow  z\sum_{0\le h\le i+k} u^h=\frac{z}{1-u}-\frac{zu^{i+k+1}}{1-u}.
\end{equation*}
Furthermore
\begin{equation*}
	F_{h+1}(z,u)=\frac{z}{1-u}F_h(z,1)-\frac{zu^{k+1}}{1-u}F_h(z,u),
\end{equation*}
and $F_0=u^h$, the starting level.

We have
\begin{equation*}
	H(z,u)=\sum_{h\ge0}F_h(z,u)=u^h+\frac{z}{1-u}H(z,1)-\frac{zu^{k+1}}{1-u}H(z,u)
\end{equation*}
or
\begin{equation*}
	H(z,u)(1-u+zu^{k+1}) =u^h(1-u)+zH(z,1).
\end{equation*}
Plugging in $\overline{u}$ into the RHS gives 0:
\begin{equation*}
	zH(z,1)=-\overline{u}^h(1-\overline{u}),
\end{equation*}
and
\begin{equation*}
	H(z,u) =\frac{u^h(1-u)-\overline{u}^h(1-\overline{u})}{1-u+zu^{k+1}}.
\end{equation*}
But we only need $H(z,0)$, since we return to the $x$-axis at the end:
\begin{equation*}
	H(z,0) = [h=0]+\overline{u}^{h+1}-\overline{u}^h  .
\end{equation*}
The total contribution of red steps is then
\begin{equation*}
	k+\sum _{h=0}^k(k-h)(\overline{u}^{h+1}-\overline{u}^h)=\sum _{h=1}^k\overline{u}^{h};
\end{equation*}
the coefficient of $z^m$ in this is the total contribution. 
Since $\overline{u}=1+z\overline{u}^{k+1}$,
there is the further simplification
\begin{equation*}
	-1+\frac1z+\frac1{1-\overline{u}}=\sum_{m\ge1}\frac{k}{m+1}\binom{(k+1)m}{m}z^m.
\end{equation*}
The proof of this is as follows. Let $m\ge1$, then
\begin{align*}
	[z^m]	\Bigl(-1+\frac1z+\frac1{1-\overline{u}}\Bigr)&=-[z^m]\frac1{z\overline{u}^{k+1}}\\
	&=-[z^{m+1}]\sum_{\ell\ge0}\frac {-(k+1)}{(k+1)\ell -(k+1)}\binom{(k+1)\ell -(k+1)}{\ell}z^\ell\\
	&=[z^{m+1}]\sum_{\ell\ge0}\frac {(k+1)}{(k+1)(\ell-1)}\binom{(k+1)(\ell-1)}{\ell}z^\ell\\
	&= \frac {(k+1)}{(k+1)m}\binom{(k+1)m}{m+1}=\frac{k}{m+1}\binom{(k+1)m}{m}.
\end{align*}

We did not expect such a simple answer  $\frac{k}{m+1}\binom{(k+1)m}{m}$ to this question about Hoppy's early adventures!

This analysis of Hoppy's early adventures covers sequences
A007226, A007228, A124724 of \cite{OEIS}, with references to \cite{AHS}.

\subsection*{Hoppy walks into negative territory}

Hoppy is now adventurous and allows himself to go to level $-1$ as well, but not deeper. The setup with generating functions is the same, but the $u$-variable counts the level relative to the $-1$ level, so this has to be corrected later. 

Hoppy, after some initial frustration discovers that he can now start with an up-step or a down-step!

First, let us start Hoppy with an up-step:
\begin{equation*}
	F(z,u)=zu^{k+1}+\frac{zu^k}{1-u}F(z,1)-\frac{zu^{k+1}}{1-u}F(z,u),
\end{equation*}
or
\begin{equation*}
	F(z,u)\frac{1-u+zu^{k+1}}{1-u}=zu^{k+1}+\frac{zu^k}{1-u}F(z,1).
\end{equation*}

Since $u-\overline{u}$ is a factor of the LHS, is must also be a factor of the RHS, and we can compute   that
\begin{equation*}
	\overline{u}(1-\overline{u})+F(z,1)=0.
\end{equation*}
The RHS is
\begin{equation*}
	\frac{zu^k}{1-u}\Big(u(1-u)   -\overline{u}(1-\overline{u}) \Big)
\end{equation*}
and
\begin{equation*}
	F(z,u)=\frac{zu^k}{1-u+zu^{k+1}}\Big(u(1-u)   -\overline{u}(1-\overline{u}) \Big).
\end{equation*}
But Hoppy can also start with a downstep! So we have to add the result of the previous computation, and get finally
\begin{equation*}
	G(z,u)=\frac{zu^k}{1-u+zu^{k+1}}\Big(u(1-u)   -\overline{u}(1-\overline{u}) \Big)+
	\frac{zu^k}{1-u+zu^{k+1}}\Big(\overline{u}-u \Big)
\end{equation*}
or better
\begin{equation*}
	G(z,u)=\frac{zu^k}{1-u+zu^{k+1}}\Big(\overline{u}^2-u^2 \Big)
\end{equation*}
Now we need
\begin{equation*}
	\frac{\partial}{\partial u}G(z,1)-G(z,1).
\end{equation*}
This substraction is necessary, since the contribution of $u^j$ is not $j$ as before but only $j-1$.
The result is
\begin{equation*}
	\frac{\overline{u}^2}{z}-2\overline{u}^2-\frac1z.
\end{equation*}

Hoppy knows that $\overline{u}^d$  has beautiful coefficients:
\begin{equation*}
\overline{u}^d=\sum_{\ell\ge0}\binom{d-1+(k+1)\ell}{\ell}\frac{d}{k\ell+d}
\end{equation*}
and he  inserts
$k=2$ (A030983):
\begin{equation*}
	3 z+16 z^2+83 z^3+442 z^4+2420 z^5+\cdots
\end{equation*}
$k=3$ (A334608):
\begin{equation*}
	5 z+34 z^2+236 z^3+1714 z^4+12922 z^5+\cdots
\end{equation*}
$k=4$ (A334610):
\begin{equation*}
	7 z+58 z^2+505 z^3+4650 z^4+44677 z^5+\cdots
\end{equation*}
In general,
\begin{equation*}
	\frac{\overline{u}^2}{z}-2\overline{u}^2-\frac1z=
	\sum_{\ell\ge0}\bigg[\binom{1+(k+1)(\ell+1)}{\ell+1}\frac{2}{k(\ell+1)+2}-2\binom{1+(k+1)\ell}{\ell}\frac{2}{k\ell+2}\bigg]z^{\ell}.
\end{equation*}
Happy Hoppy decides to stop this line of computations here.

\section{Combinatorics of sequence A002212 in \cite{OEIS}}

The following sections gives some (mostly new) results about the sequence
\begin{equation*}
	1, 1, 3, 10, 36, 137, 543, 2219, 9285, 39587, 171369, 751236, 3328218, 14878455,\dots,
\end{equation*}
which is A002212 in \cite{OEIS}.

Here is the plan about the structures enumerated by this sequence:

Hex-trees \cite{KimStanley}; they are identified as weighted unary-binary trees, with weight one. Apart from left and right branches, as
in binary trees, there are also unary branches, and they can come in different colours, here in just one colour.
Unary-binary trees played a role in the present authors scientific development, as documented in \cite{FlPr86}, a paper written with
the late and great Philippe Flajolet, about the register function (Horton-Strahler numbers) of unary-binary trees. Here, we can offer an improvement, using
a ``better'' substitution than in \cite{FlPr86}. The results can now be made fully explicit. As a by-product, this provides a definition and analysis
of the Horton-Strahler numbers of Hex-trees. An introductory section (about binary trees) provides all the basics.

Then we move to skew Dyck paths \cite{Deutsch-italy}. They are like Dyck paths, but allow for an extra step $(-1,-1)$, provided that the path does not intersect itself.
An equivalent model, defined and described using a bijection, is from \cite{Deutsch-italy}: Marked ordered trees. They are like ordered trees, with an additional feature,
namely each rightmost edge (except for one that leads to a leaf) can be coloured with two colours. Since we find this class of trees to be interesting, we
analyze two parameters of them: number of leaves and height. While the number of leaves for ordered trees is about $n/2$, it is only $n/10$ in the new model.
For the height, the leading term $\sqrt{\pi n}$ drops to $\frac{2}{\sqrt 5}\sqrt{\pi n}$. Of course, many more parameters of this new class of trees could be investigated,
which we encourage to do.

More about skew Dyck paths appears in a later section.

The next two classes of structures are multi-edge trees. Our interest in them was already triggered in an earlier publication \cite{HPW}. They may be seen as ordered trees, but
with weighted edges. The weights are integers $\ge1$, and a weight $a$ may be interpreted as $a$ parallel edges. The other class are 3-Motzkin paths. They are like
Motzkin paths (Dyck paths plus horizontal steps); however, the horizontal steps come in three different colours. A bijection is described. Since 3-Motzkin paths and multi-edge trees
are very much alike (using a variation of the classical rotation correspondence), all the structures that are discussed in this paper can be linked via bijections.

\section{Binary trees and Horton-Strahler numbers}

This section is classical and serves as the basis of some new developments about weighted unary-binary trees.
A full account can be found here: \cite{EATCS}.

Binary trees may be expressed by the following symbolic equation, which says that they include the empty tree 
and trees recursively built from a root followed by two subtrees, which are binary trees:
\begin{center}\small
	\begin{tikzpicture}
		[inner sep=1.3mm,
		s1/.style={circle=10pt,draw=black!90,thick},
		s2/.style={rectangle,draw=black!50,thick},scale=0.5]
		
		\node at ( -4.8,0) { $\mathscr{B}$};
		
		\node at (-3,0) { $=$};
		\node(c) at (-1.5,0){ $\qed$};
		\node at (0.7,0) {$+$};
		%\node at (1.5,0) {$2$};
		\node(d) at (3,1)[s1]{};
		\node(e) at (2,-1){ $\mathscr{B}$};
		\node(f) at (4,-1){ $\mathscr{B}$};
		\path [draw,-,black!90] (d) -- (e) node{};%
		\path [draw,-,black!90] (d) -- (f) node{};%
		
	\end{tikzpicture}
\end{center}

Binary trees are counted by Catalan numbers and there is an important  parameter \textsf{reg}, which 
in Computer Science is called the register function. It associates to each binary tree (which is used to code an arithmetic expression,
with data in the leaves and operators in the internal nodes) the minimal number of extra registers that is needed to evaluate the tree. The optimal strategy is to evaluate difficult subtrees first, and use one register to keep its value, which does not hurt, if the other subtree requires less registers. If both subtrees are equally difficult, then one more register is used, compared to the requirements of the subtrees. This natural parameter is among Combinatorialists known as the Horton-Strahler numbers, and we will adopt this
name throughout this paper.

There is a recursive description of this function: $\textsf{reg}(\square)=0$, and if tree $t$ has subtrees $t_1$ and $t_2$, then
\begin{equation*}
	\textsf{reg}(t)=
	\begin{cases}
		\max\{\textsf{reg}(t_1),\textsf{reg}(t_2)\}&\text{ if } \textsf{reg}(t_1)\ne\textsf{reg}(t_2),\\
		1+\textsf{reg}(t_1)&\text{ otherwise}.
	\end{cases}
\end{equation*}

The recursive description attaches numbers to the nodes, starting with 0's at the leaves and
then going up; the number appearing at the root is the Horton-Strahler number of the tree.
\begin{center}\tiny
	\begin{tikzpicture}
		[scale=0.4,inner sep=0.7mm,
		s1/.style={circle,draw=black!90,thick},
		s2/.style={rectangle,draw=black!90,thick}]
		\node(a) at ( 0,8) [s1] [text=black]{$\boldsymbol{2}$};
		\node(b) at ( -4,6) [s1] [text=black]{$1$};
		\node(c) at ( 4,6) [s1] [text=black]{$2$};
		\node(d) at ( -6,4) [s2] [text=black]{$0$};
		\node(e) at ( -2,4) [s1] [text=black]{$1$};
		\node(f) at ( 2,4) [s1] [text=black]{$1$};
		\node(g) at ( 6,4) [s1] [text=black]{$1$};
		\node(h) at ( -3,2) [s2] [text=black]{$0$};
		\node(i) at ( -1,2) [s2] [text=black]{$0$};
		\node(j) at ( 1,2) [s2] [text=black]{$0$};
		\node(k) at ( 3,2) [s2] [text=black]{$0$};
		\node(l) at ( 5,2) [s2] [text=black]{$0$};
		\node(m) at ( 7,2) [s2] [text=black]{$0$};
		\path [draw,-,black!90] (a) -- (b) node{};%
		\path [draw,-,black!90] (a) -- (c) node{};%
		\path [draw,-,black!90] (b) -- (d) node{};%
		\path [draw,-,black!90] (b) -- (e) node{};%
		\path [draw,-,black!90] (c) -- (f) node{};%
		\path [draw,-,black!90] (c) -- (g) node{};%
		\path [draw,-,black!90] (e) -- (h) node{};%
		\path [draw,-,black!90] (e) -- (i) node{};%
		\path [draw,-,black!90] (f) -- (j) node{};%
		\path [draw,-,black!90] (f) -- (k) node{};%
		\path [draw,-,black!90] (g) -- (l) node{};%
		\path [draw,-,black!90] (g) -- (m) node{};%
	\end{tikzpicture}
\end{center}
%\end{figure}

Let $\mathscr{R}_{p}$ denote the family of
trees with Horton-Strahler number $=p$, then one gets immediately from the recursive 
definition:
\begin{center}\small
	\begin{tikzpicture}
		[inner sep=1.3mm,
		s1/.style={circle=10pt,draw=black!90,thick},
		s2/.style={rectangle,draw=black!50,thick},scale=0.5]
		
		\node at ( -5,0) { $\mathscr{R}_p$};
		
		\node at (-4,0) { $=$};
		\node(a) at (-2,1)[s1]{};
		\node(b) at (-3,-1){ $\mathscr{R}_{p-1}$};
		\node(c) at (-1,-1){ $\mathscr{R}_{p-1}$};
		\path [draw,-,black!90] (a) -- (b) node{};%
		\path [draw,-,black!90] (a) -- (c) node{};%
		\node at (0.7,0) {$+$};
		%\node at (1.5,0) {$2$};
		\node(d) at (3,1)[s1]{};
		\node(e) at (2,-1){ $\mathscr{R}_{p}$};
		\node(f) at (4,-1.2){ $\sum\limits_{j<p}\mathscr{R}_{j} $};
		\path [draw,-,black!90] (d) -- (e) node{};%
		\path [draw,-,black!90] (d) -- (f) node{};%
		\node at (5+0.7,0) {$+$};
		%5\node at (5+1.5,0) {$2$};
		\node(dd) at (5.5+3,1)[s1]{};
		\node(ee) at (5.5+2,-1.2){ $\sum\limits_{j<p}\mathscr{R}_{j}$};
		\node(ff) at (5.5+4,-1){ $\mathscr{R}_{p}$};
		\path [draw,-,black!90] (dd) -- (ee) node{};%
		\path [draw,-,black!90] (dd) -- (ff) node{};%
	\end{tikzpicture}
\end{center}
In terms of generating functions, these equations read as
\begin{equation*}
	R_p(z)=zR_{p-1}^2(z)+2zR_p(z)\sum_{j<p}R_j(z);
\end{equation*}
the variable $z$ is used to mark the size (i.~e., the number of internal nodes) of the binary tree.

A historic account of these concepts, from the angle of Philippe Flajolet, who was one of the pioneers
is \cite{register-introduction}, compare also \cite{ECA-historic}.

Amazingly, the recursion for the generating functions $R_p(z)$  can be solved explicitly! The substitution
\begin{equation*}
	z=\frac{u}{(1+u)^2}
\end{equation*}
that de Bruijn, Knuth, and Rice~\cite{BrKnRi72} also used, produces the nice expression
\begin{equation*}
	R_p(z)=\frac{1-u^2}{u}\frac{u^{2^p}}{1-u^{2^{p+1}}}.
\end{equation*}
Of course, once this is \emph{known}, it can be proved by induction, using the recursive formula. For the readers benefit, this will be sketched now. 

We start with the auxiliary formula
\begin{equation*}
	\sum_{0\le j<p}\frac{u^{2^j}}{1-u^{2^{j+1}}}=\frac{u}{1-u}-\frac{u^{2^p}}{1-u^{2^{p}}},
\end{equation*}
which we will  prove now by induction: For $p=0$, the formula $0=\frac{u}{1-u}-\frac{u}{1-u}$ is correct, and then
\begin{align*}
	\sum_{0\le j<p+1}
	\frac{u^{2^j}}{1-u^{2^{j+1}}}
	&=\frac{u}{1-u}-\frac{u^{2^p}}{1-u^{2^{p}}}+\frac{u^{2^p}}{1-u^{2^{p+1}}}\\
	&=\frac{u}{1-u}-\frac{u^{2^p}(1+u^{2^p})}{1-u^{2^{p+1}}}+\frac{u^{2^p}}{1-u^{2^{p+1}}}
	=\frac{u}{1-u}-\frac{u^{2^{p+1}}}{1-u^{2^{p+1}}}.
\end{align*}
Now the formula for $R_p(z)$ can also be proved by induction. First, $R_0(z)=\frac{1-u^2}{u}\frac{u}{1-u^{2}}=1$, as it should, and
\begin{align*}
	zR_{p-1}^2(z)&+2zR_p(z)\sum_{j<p}R_j(z)\\
	&=\frac{u}{(1+u)^2}\frac{(1-u^2)^2}{u^2}\frac{u^{2^{p}}}{(1-u^{2^{p}})^2}
	+\frac{2u}{(1+u)^2}R_p(z)\sum_{j<p}\frac{1-u^2}{u}\frac{u^{2^j}}{1-u^{2^{j+1}}}\\
	&= \frac{u^{2^{p}-1}(1-u)^2}{(1-u^{2^{p}})^2}
	+\frac{2(1-u)}{(1+u)}R_p(z)\sum_{j<p}\frac{u^{2^j}}{1-u^{2^{j+1}}}.
\end{align*}
Solving
\begin{align*}
	R_p(z)= \frac{u^{2^{p}-1}(1-u)^2}{(1-u^{2^{p}})^2}
	+\frac{2(1-u)}{(1+u)}R_p(z)\bigg[\frac{u}{1-u}-\frac{u^{2^p}}{1-u^{2^{p}}}\bigg]
\end{align*}
leads to
\begin{align*}
	R_p(z)\frac{1-u}{1+u}\bigg[1+2\frac{u^{2^p}}{1-u^{2^{p}}}\bigg]= \frac{u^{2^{p}-1}(1-u)^2}{(1-u^{2^{p}})^2},
\end{align*}
or, simplified
\begin{align*}
	R_p(z)= \frac{u^{2^{p}-1}(1-u^2)}{(1-u^{2^{p}})(1+u^{2^{p}})}
	=\frac{1-u^2}{u}\frac{u^{2^{p}}}{(1-u^{2^{p+1}})},
\end{align*}
which is the formula that we needed to prove. \qed

\subsection*{Weighted unary-binary trees and Horton-Strahler numbers}

%\section{Unary-binary trees and Hex-trees}

The family of unary-binary trees ${\mathscr{M}}$ might be defined by the symbolic equation
\begin{center}
	\begin{tikzpicture}
		[inner sep=1.3mm,
		s1/.style={circle=10pt,draw=black!90,thick},
		s2/.style={rectangle,draw=black!50,thick},scale=0.5]
		
		\node at ( -12.8,0.1) { ${\mathscr{M}}$};
		
		\node at (-11.2,0) { $=$};
		\node at (-4.5,0) { $+$};
		\node(a) at (-2,1)[s1]{};
		\node(b) at (-3,-1){ ${\mathscr{M}}$};
		\node(c) at (-1,-1){ ${\mathscr{M}}$};
		\path [draw,-,black!90] (a) -- (b) node{};%
		\path [draw,-,black!90] (a) -- (c) node{};%
		\node(a1) at (-6.5,1)[s1]{};
		\node(b1) at (-6.5,-1){ $ {\mathscr{M}}\setminus\{\square\}$};
		\path [draw,-,black!90] (a1) -- (b1) node{};%
		\node at (-9.0,0) { $\square\ \ +$};
	\end{tikzpicture}
\end{center}
The equation for the generating function is
\begin{equation*}
	M=1+z(M-1)+zM^2
\end{equation*}
with the solution
\begin{equation*}
	M=M(z)=\frac{1-z-\sqrt{1-6z+5z^2}}{2z}=1+z+3{z}^{2}+10{z}^{3}+36{z}^{4}+\cdots;
\end{equation*}
the coefficients form again sequence A002212 in \cite{OEIS} and enumerate Schr\"oder paths, among many other things.
We will come to equivalent structures a bit later. 

In the instance of unary-binary trees, we can also work with a substitution: Set $z=\frac{u}{1+3u+u^2}$, then
$M(z)=1+u$. Unary-binary trees and the register function were investigated in \cite{FlPr86}, but the present favourable substitution was not used.
Therefore, in this previous paper, asymptotic results were available but no explicit formulae.

This works also with a weighted version, where we allow unary edges with $a$ different colours. Then
\begin{center}
	\begin{tikzpicture}
		[inner sep=1.3mm,
		s1/.style={circle=10pt,draw=black!90,thick},
		s2/.style={rectangle,draw=black!50,thick},scale=0.5]
		
		\node at ( -12.8,0.1) { ${\mathscr{N}}$};
		
		\node at (-11.5,0) { $=$};
		\node at (-4.5,0) { $+$};
		\node(a) at (-2,1)[s1]{};
		\node(b) at (-3,-1){ ${\mathscr{N}}$};
		\node(c) at (-1,-1){ ${\mathscr{N}}$};
		\path [draw,-,black!90] (a) -- (b) node{};%
		\path [draw,-,black!90] (a) -- (c) node{};%
		\node at (-7.5,0){$a\ \cdot$};
		\node(a1) at (-6.5,1)[s1]{};
		\node(b1) at (-6.5,-1){ $ {\mathscr{N}}\setminus\{\square\}$};
		\path [draw,-,black!90] (a1) -- (b1) node{};%
		\node at (-9.5,0) { $\square\ \ +$};
	\end{tikzpicture}
\end{center}
and with the substitution $z=\frac{u}{1+(a+2)u+u^2}$, the generating function is beautifully expressed as $N(z)=1+u$.
For $a=0$, this covers also binary trees.

We will consider  the Horton-Strahler numbers of unary-binary trees in the sequel. The definition is naturally extended by

\begin{center}
	\vspace{0pt}\begin{tikzpicture}
		[inner sep=0.6mm,
		s1/.style={circle=1pt,draw=black!90,thick}]
		\node[] at ( -0.300,-0.10) {\textsf{reg}\bigg(};
		\node[] at ( 0.6500,-0.10) {\bigg)};
		\node[] at ( 1.5500,-0.10) {=\ \textsf{reg}(t).};
		\path [draw,-,black!90 ] (0.3,0.34) -- (0.3,-0.350) ;
		\node [s1]at ( 0.300,0.4) { };
		\node[] at ( 0.300,-0.60) {$t$ };
	\end{tikzpicture} 
\end{center}

Now we can move again to $R_p(z)$, the generating funciton of (generalized) unary-binary trees with Horton-Strahler number $=p$.
The recursion (for $p\ge1$) is
\begin{center}\small
	\begin{tikzpicture}
		[inner sep=1.3mm,
		s1/.style={circle=10pt,draw=black!90,thick},
		s2/.style={rectangle,draw=black!50,thick},scale=0.5]
		
		\node at ( -5,0) { $\mathscr{R}_p$};
		
		\node at (-4,0) { $=$};
		\node(a) at (-2,1)[s1]{};
		\node(b) at (-3,-1){ $\mathscr{R}_{p-1}$};
		\node(c) at (-1,-1){ $\mathscr{R}_{p-1}$};
		\path [draw,-,black!90] (a) -- (b) node{};%
		\path [draw,-,black!90] (a) -- (c) node{};%
		\node at (0.7,0) {$+$};
		%\node at (1.5,0) {$2$};
		\node(d) at (3,1)[s1]{};
		\node(e) at (2,-1){ $\mathscr{R}_{p}$};
		\node(f) at (4,-1.2){ $\sum\limits_{j<p}\mathscr{R}_{j} $};
		\path [draw,-,black!90] (d) -- (e) node{};%
		\path [draw,-,black!90] (d) -- (f) node{};%
		\node at (5+0.7,0) {$+$};
		%5\node at (5+1.5,0) {$2$};
		\node(dd) at (5.5+3,1)[s1]{};
		\node(ee) at (5.5+2,-1.2){ $\sum\limits_{j<p}\mathscr{R}_{j}$};
		\node(ff) at (5.5+4,-1){ $\mathscr{R}_{p}$};
		\path [draw,-,black!90] (dd) -- (ee) node{};%
		\path [draw,-,black!90] (dd) -- (ff) node{};%
		\node(dd) at (13,1)[s1]{};
		\node(ee) at (13,-1){ $\mathscr{R}_{p}$};
		\path [draw,-,black!90] (dd) -- (ee) node{};%
		\node at (11.5,0) {$+\ \ a\cdot$};
	\end{tikzpicture}
\end{center}
In terms of generating functions, these equations read as
\begin{equation*}
	R_p(z)=zR_{p-1}^2(z)+2zR_p(z)\sum_{j<p}R_j(z)+azR_p(z), \quad p\ge1;\quad R_0(z)=1.
\end{equation*}
Amazingly, with the substitution $z=\frac{u}{1+(a+2)u+u^2}$, formally we get the \emph{same} solution as in the binary case:
\begin{equation*}
	R_p(z)=\frac{1-u^2}{u}\frac{u^{2^p}}{1-u^{2^{p+1}}}.
\end{equation*}
The proof by induction is as before. One sees another advantage of the substitution: On a formal level, many manipulations do not need to be repeated. Only when one switches back to the $z$-world, things become different.

Now we move to Hex-trees.
\begin{center}
	\begin{tikzpicture}
		[inner sep=1.3mm,
		s1/.style={circle=10pt,draw=black!90,thick},
		s2/.style={rectangle,draw=black!50,thick},scale=0.5]
		
		\node at ( -12.4,0.1) { ${\mathscr{H}}$};
		
		\node at (-11.0,0) { $=$};
		
		\node(a) at (-1,1)[s1]{};
		\node(b) at (-3,-1){ ${\mathscr{H}\setminus\{\square\}}$};
		\node(c) at (1,-1){ ${\mathscr{H}\setminus\{\square\}}$};
		\path [draw,-,black!90] (a) -- (b) node{};%
		\path [draw,-,black!90] (a) -- (c) node{};%
		
		\begin{scope}[xshift=13cm]
			\node at (-10.0,0) { $+$};
			\node at (-4.5,0) { $+$};
			\node(a1) at (-6.5,1)[s1]{};
			\node(b1) at (-7.5,-1){ $ {\mathscr{H}}\setminus\{\square\}$};
			\path [draw,-,black!90] (a1) -- (b1) node{};%
		\end{scope}
		
		\begin{scope}[xshift=21cm]
			%\node at (-4.5,0) { $+$};
			\node(a1) at (-6.5,1)[s1]{};
			\node(b1) at (-5.5,-1){ $ {\mathscr{H}}\setminus\{\square\}$};
			\path [draw,-,black!90] (a1) -- (b1) node{};%
		\end{scope}

		\begin{scope}[xshift=17cm]
			\node at (-4.5,0) { $+$};
			\node(a1) at (-6.5,1)[s1]{};
			\node(b1) at (-6.5,-1){ $ {\mathscr{H}}\setminus\{\square\}$};
			\path [draw,-,black!90] (a1) -- (b1) node{};%
		\end{scope}

		\node at (-9.0,0) { $\square\ \ +$};
		\node[s1] at (-7.0,0) { };
		\node at (-5.0,0) {$+$ };
	\end{tikzpicture}
\end{center}

\subsection*{Hex trees}

Hex trees either have two non-empty successors, or one of 3 types of unary successors (called left, middle, right).
The author has seen this family first in \cite{KimStanley}, but one can find older literature following the references and the usual
search engines.

The generating function satisfies
\begin{align*}
	H&(z)=1+z(H(z)-1)^2+z+3z(H(z)-1)=\frac{1-z-\sqrt{(1-z)(1-5z)}}{2z}\\
	&=1+z+3{z}^{2}+10{z}^{3}+36{z}^{4}+137{z}^{5}+543{z}^{6}+2219
	{z}^{7}+9285{z}^{8}+39587{z}^{9}+\cdots.
\end{align*}
The same generating function also appears in \cite{HPW}, and it is again sequence A002212 in \cite{OEIS}.
One can rewrite the symbolic equation as
\begin{center}
	\begin{tikzpicture}
		[inner sep=1.3mm,
		s1/.style={circle=10pt,draw=black!90,thick},
		s2/.style={rectangle,draw=black!50,thick},scale=0.5]
		
		\node at ( -12.4,0.1) { ${\mathscr{H}}$};
		
		\node at (-11.0,0) { $=$};
		
		\begin{scope}[xshift=-4cm]
			\node(a) at (-1,1)[s1]{};
			\node(b) at (-3,-1){ $\mathscr{H}$};
			\node(c) at (1,-1){ $\mathscr{H}$};
			\path [draw,-,black!90] (a) -- (b) node{};%
			\path [draw,-,black!90] (a) -- (c) node{};%
		\end{scope}

		\begin{scope}[xshift=7.5cm]
			\node at (-8.5,0) { $+$};
			\node(a1) at (-6.5,1)[s1]{};
			\node(b1) at (-6.5,-1){ $ {\mathscr{H}}\setminus\{\square\}$};
			\path [draw,-,black!90] (a1) -- (b1) node{};%
		\end{scope}

		\node at (-9.0,0) { $\square\ \ \ +$};
		
	\end{tikzpicture}
\end{center}
and sees in this way that the Hex-trees are just unary-binary trees (with parameter $a=1$).

\subsection*{Continuing with enumerations}

First, we will enumerate the number of (generalized) unary-binary trees with $n$ (internal) nodes. For that we need the notion 
of generalized trinomial coefficients, viz.
\begin{equation*}
	\binom{n;1,a,1}{k}:=[z^k](1+az+z^2)^n.
\end{equation*}
Of course, for $a=2$, this simplifies to a binomial coefficient $\binom{2n}{k}$. We will use contour integration to pull out coefficients, and the contour of integer, in whatever variable, is a small circle (or equivalent) around the origin.
The desired number is
\begin{align*}
	[z^n](1+u)&=\frac1{2\pi i}\oint \frac{dz}{z^{n+1}}(1+u)\\
	&=\frac1{2\pi i}\oint \frac{du(1-u^2)(1+(a+2)u+u^2)^{n+1}}{(1+(a+2)u+u^2)^2u^{n+1}}(1+u)\\
	&=[u^{n+1}](1-u)(1+u)^2(1+(a+2)u+u^2)^{n-1}\\
	&=\binom{n-1;1,a+2,1}{n+1}+\binom{n-1;1,a+2,1}{n}\\*
	&\hspace*{4cm}-\binom{n-1;1,a+2,1}{n-1}-\binom{n-1;1,a+2,1}{n-2}.
\end{align*}
Then we introducte $S_p(z)=R_{p}(z)+R_{p+1}(z)+R_{p+2}(z)+\cdots$, the generating function of trees with Horton-Strahler number
$\ge p$. Using the summation formula proved earlier, we get
\begin{equation*}
	S_p(z)=\frac{1-u^2}{u}\frac{u^{2^p}}{1-u^{2^{p}}}=
	\frac{1-u^2}{u}\sum_{k\ge1}u^{k2^p}.
\end{equation*}
Further,
\begin{align*}
	[z^n]S_p(z)&=\sum_{k\ge1}\frac1{2\pi i}\oint \frac{dz}{z^{n+1}}\frac{1-u^2}{u}u^{k2^p}.
\end{align*}

\subsection*{Asymptotics}

We start by deriving asymptotics for the number of (generalized) unary-binary trees with $n$ (internal) nodes. This is
a standard application of singularity analysis of generating functions, as described in \cite{FlOd90} and \cite{FS}.

We start from the generating function
\begin{equation*}
	N(z)=\frac{1-az-\sqrt{1-2(a+2)z+a(a+4)z^2}}{2z}
\end{equation*}
and determine the singularity closest to the origin, which is the value making the square root disappear:
$z=\frac1{a+4}$.
After that, the local expansion of $N(z)$ around this singularity is determined:
\begin{equation*}
	N(z) \sim	2-\sqrt{a+4}\sqrt{1-(a+4)z}.
\end{equation*}
The translation lemmas given in \cite{FlOd90} and \cite{FS} provide the asymptotics:
\begin{align*}
	[z^n]N(z)&\sim [z^n]\Big(2-\sqrt{a+4}\sqrt{1-(a+4)z}\Big)\\&
	=-\sqrt{a+4}(a+4)^n\frac{n^{-3/2}}{\Gamma(-\frac12)}=(a+4)^{n+1/2}\frac{1}{2\sqrt\pi n^{3/2}}.
\end{align*}
Just note that $a=0$ is the well-known formula for binary trees with $n$ nodes.

Now we move to the generating function for the average number of registers. Apart from normalization it is
\begin{align*}
	\sum_{p\ge1}pR_p(z)&=\sum_{p\ge1}S_p(z)=\frac{1-u^2}{u}\sum_{p\ge1}\sum_{k\ge1}u^{k2^p}\\
	&=\frac{1-u^2}{u}\sum_{n\ge1}v_2(n)u^n,
\end{align*}
where $v_2(n)$ is the highest exponent $k$ such $2^k$ divides $n$.

This has to be studied around $u=1$, which, upon setting $u=e^{-t}$, means around $t=0$.
Eventually, and that is the only thing that is different here, this is to be retranslated  into a singular expansion of $z$ around
its singularity, which depends on the parameter $a$. 

For the reader's convenience, we also repeat the steps that were known before.
The first factor is elementary:
\begin{equation*}
	\frac{1-u^2}{u}\sim2t+{\frac {1}{3}}{t}^{3}+\cdots
\end{equation*}
For 
\begin{equation*}
	\sum_{p\ge1}\sum_{k\ge1}e^{-k2^pt},
\end{equation*}
one applies the Mellin transform, with the result
\begin{equation*}
	\frac{\Gamma(s)\zeta(s)}{2^s-1}.
\end{equation*}
Applying the inversion formula, one finds
\begin{equation*}
	\sum_{p\ge1}\sum_{k\ge1}e^{-k2^pt}=\frac1{2\pi i}\int_{2-i\infty}^{2+i\infty}t^{-s}\frac{\Gamma(s)\zeta(s)}{2^s-1}ds.
\end{equation*}
Shifting the line of integration to the left, the residues at the poles $s=1$, $s=0$, $s=\chi_k=\frac{2k\pi i}{\log2}$, $k\neq0$ provide enough terms for our asymptotic expansion.
\begin{equation*}
	\frac1{t}+{\frac {\gamma}{2\log2 }}-\frac14-	\frac {\log   \pi   }{2\log2 }+\frac {\log t }{2\log2}
	+\frac1{\log2}\sum_{k\neq0}\Gamma(\chi_k)\zeta(\chi_k)t^{-\chi_k}.
\end{equation*}
Combined with the elementary factor, this leads to
\begin{equation*}
	2+\Big(\frac {\gamma}{\log2 }-\frac12-\frac {\log   \pi   }{\log2 }+\frac {\log t }{\log2}\Big)t+\frac{2t}{\log2}\sum_{k\neq0}\Gamma(\chi_k)\zeta(\chi_k)t^{-\chi_k}+O(t^2\log t).
\end{equation*}
Now we want to translate into the original $z$-world. Since $z=\frac{u}{1+(a+2)u+u^2}$, $u=1$ translates into the singularity $z=\frac{1}{4+a}$. Further,
\begin{equation*}
	t\sim \sqrt{4+a}\cdot \sqrt{1-z(4+a)},
\end{equation*}
let us abbreviate $A=4+a$, then for singularity analysis we must consider
\begin{align*}
	&\frac {\sqrt{A}\cdot \sqrt{1-zA}\log (1-zA) }{2\log2}\\
	&+	\Big(\frac {\gamma}{\log2 }-\frac12-\frac {\log   \pi   }{\log2 }+\frac{\log A}{2\log 2}\Big)\sqrt{A}\cdot \sqrt{1-zA}\\
	&+\frac{2  }{\log2}\sum_{k\neq0}\Gamma(\chi_k)\zeta(\chi_k) A^{\frac{1-\chi_k}2}(1-zA)^{\frac{1-\chi_k}2}.
\end{align*}
The formula that is perhaps less known and needed here is \cite{FlOd90}
\begin{align*}
	[z^n]\log(1-z)\sqrt{1-z}\sim \frac{n^{-3/2}\log n}{2\sqrt \pi}+\frac{n^{-3/2}}{2\sqrt \pi}(-2+\gamma +2\log2);
\end{align*}
furthermore we need
\begin{equation*}
	[z^n](1-z)^\alpha \sim \frac{n^{-\alpha-1}}{\Gamma(-\alpha)}.
\end{equation*}
We start with the most complicated term:
\begin{align*}
	\frac{[z^n]\frac {\sqrt{A}\cdot \sqrt{1-zA}\log (1-zA) }{2\log2}}{[z^n]N(z)}
	&\sim \frac {\sqrt{A}}{2\log2}\frac{A^n\Big(\frac{n^{-3/2}\log n}{2\sqrt \pi}+\frac{n^{-3/2}}{2\sqrt \pi}(-2+\gamma +2\log2)\Big)}
	{A^{n+1/2}\frac{1}{2\sqrt\pi n^{3/2}}}\\
	&= \log_4 n+1+ \frac{\gamma }{2\log2}- \frac{1}{\log2}.
\end{align*}
The next term we consider is
\begin{align*}
	\Big(\frac {\gamma}{\log2 }-\frac12-\frac {\log   \pi   }{\log2 }+\frac{\log A}{2\log 2}\Big)&\sqrt{A}\frac{[z^n] \sqrt{1-zA}}{[z^n]N(z)}\\*
	&\sim
	\Big(\frac {\gamma}{\log2 }-\frac12-\frac {\log   \pi   }{\log2 }+\frac{\log A}{2\log 2}\Big)\sqrt{A}\frac{[z^n] \sqrt{1-zA}}{-\sqrt{A}[z^n]\sqrt{1-zA}}\\
	&=-\frac {\gamma}{\log2 }+\frac12+\frac {\log   \pi   }{\log2 }-\frac{\log A}{2\log 2}.
\end{align*}
The last term we consider is
\begin{align*}
	\frac{2  }{\log2}&\Gamma(\chi_k)\zeta(\chi_k) A^{\frac{1-\chi_k}2}\frac{[z^n](1-zA)^{\frac{1-\chi_k}2}}{-\sqrt{A}[z^n]\sqrt{1-zA}}\\
	&\sim-\frac{4 \sqrt\pi }{\log2}\frac{\Gamma(\chi_k)\zeta(\chi_k)}{\Gamma\big(\frac{\chi_k-1}{2}\big)} A^{\frac{1-\chi_k}2}n^{\chi_k/2}.
\end{align*}
Eventually we have evaluated the average value of the Horton-Strahler numbers:
\begin{theorem}The average Horton-Strahler of weighted unary-binary trees with $n$ nodes is given by the asymptotic formula
	\begin{align*}
		\log_4 n&- \frac{\gamma }{2\log2}- \frac{1}{\log2}+\frac32+\frac {\log   \pi   }{\log2 }-\frac{\log A}{2\log 2}
		-\frac{4 \sqrt{\pi A} }{\log2}\sum_{k\neq0}\frac{\Gamma(\chi_k)\zeta(\chi_k)}{\Gamma\big(\frac{\chi_k-1}{2}\big)} A^{\frac{-\chi_k}2}n^{\chi_k/2}\\
		&=\log_4 n- \frac{\gamma }{2\log2}- \frac{1}{\log2}+\frac32+\frac {\log   \pi   }{\log2 }-\frac{\log A}{2\log 2}+\psi(\log_4n),
	\end{align*}
	with a tiny periodic function $\psi(x)$ of period 1.
\end{theorem}	

%\clearpage

\section{Marked ordered trees}

In \cite{Deutsch-italy} we find the following variation of ordered trees: Each rightmost edge might be marked or not, if it does not lead to an endnode (leaf).
We depict a marked edge by the red colour and draw all of them of size 4 (4 nodes):
\begin{figure}[h]
	\begin{tikzpicture}[scale=0.7]

		\draw[black,fill=black] (0,0) circle (.5ex);
		\draw[black,fill=black] (0,-1) circle (.5ex);
		\draw[black,fill=black] (0,-2) circle (.5ex);
		\draw[black,fill=black] (0,-3) circle (.5ex);
		\draw [thick] (0,0) -- (0,-1)-- (0,-3) ;%
		
		\draw[black,fill=black,xshift=1cm] (0,0) circle (.5ex);
		\draw[black,fill=black,xshift=1cm] (0,-1) circle (.5ex);
		\draw[black,fill=black,xshift=1cm] (0,-2) circle (.5ex);
		\draw[black,fill=black,xshift=1cm] (0,-3) circle (.5ex);
		\draw [ultra thick,red,xshift=1cm] (0,0) -- (0,-1) ;
		\draw [thick,xshift=1cm] (0,-1) -- (0,-2) ;
		\draw [thick,xshift=1cm] (0,-2) -- (0,-3) ;
		
		\draw[black,fill=black,xshift=2cm] (0,0) circle (.5ex);
		\draw[black,fill=black,xshift=2cm] (0,-1) circle (.5ex);
		\draw[black,fill=black,xshift=2cm] (0,-2) circle (.5ex);
		\draw[black,fill=black,xshift=2cm] (0,-3) circle (.5ex);
		\draw [thick,xshift=2cm] (0,0) -- (0,-1) ;
		\draw [ultra thick,red,xshift=2cm] (0,-1) -- (0,-2) ;
		\draw [thick,xshift=2cm] (0,-2) -- (0,-3) ;
		
		\draw[black,fill=black,xshift=3cm] (0,0) circle (.5ex);
		\draw[black,fill=black,xshift=3cm] (0,-1) circle (.5ex);
		\draw[black,fill=black,xshift=3cm] (0,-2) circle (.5ex);
		\draw[black,fill=black,xshift=3cm] (0,-3) circle (.5ex);
		\draw [ultra thick,red, xshift=3cm] (0,0) -- (0,-1) ;
		\draw [ultra thick,red,xshift=3cm] (0,-1) -- (0,-2) ;
		\draw [thick,xshift=3cm] (0,-2) -- (0,-3) ;
		
		\draw[black,fill=black,xshift=6cm] (0,0) circle (.5ex);
		\draw[black,fill=black,xshift=6cm] (-1,-1) circle (.5ex);
		\draw[black,fill=black,xshift=6cm] (-1,-2) circle (.5ex);
		\draw[black,fill=black,xshift=6cm] (1,-1) circle (.5ex);
		\draw [thick, xshift=6cm] (0,0) -- (-1,-1) ;
		\draw [thick,xshift=6cm] (-1,-1) -- (-1,-2) ;
		\draw [thick,xshift=6cm] (0,0) -- (1,-1) ;
		
		\draw[black,fill=black,xshift=9cm] (0,0) circle (.5ex);
		\draw[black,fill=black,xshift=9cm] (-1,-1) circle (.5ex);
		\draw[black,fill=black,xshift=9cm] (1,-2) circle (.5ex);
		\draw[black,fill=black,xshift=9cm] (1,-1) circle (.5ex);
		\draw [thick, xshift=9cm] (0,0) -- (-1,-1) ;
		\draw [thick,xshift=9cm] (0,0) -- (1,-1) ;
		\draw [thick,xshift=9cm] (1,-1) -- (1,-2) ;
		
		\draw[black,fill=black,xshift=12cm] (0,0) circle (.5ex);
		\draw[black,fill=black,xshift=12cm] (-1,-1) circle (.5ex);
		\draw[black,fill=black,xshift=12cm] (1,-2) circle (.5ex);
		\draw[black,fill=black,xshift=12cm] (1,-1) circle (.5ex);
		\draw [thick, xshift=12cm] (0,0) -- (-1,-1) ;
		\draw [ultra thick,red,xshift=12cm] (0,0) -- (1,-1) ;
		\draw [thick,xshift=12cm] (1,-1) -- (1,-2) ;
		
		\draw[black,fill=black,xshift=15cm] (0,0) circle (.5ex);
		\draw[black,fill=black,xshift=15cm] (0,-1) circle (.5ex);
		\draw[black,fill=black,xshift=15cm] (1,-2) circle (.5ex);
		\draw[black,fill=black,xshift=15cm] (-1,-2) circle (.5ex);
		\draw [thick, xshift=15cm] (0,0) -- (0,-1) ;
		\draw [thick,xshift=15cm] (0,-1) -- (1,-2) ;
		\draw [thick,xshift=15cm] (0,-1) -- (-1,-2) ;
		
		\draw[black,fill=black,xshift=18cm] (0,0) circle (.5ex);
		\draw[black,fill=black,xshift=18cm] (0,-1) circle (.5ex);
		\draw[black,fill=black,xshift=18cm] (1,-2) circle (.5ex);
		\draw[black,fill=black,xshift=18cm] (-1,-2) circle (.5ex);
		\draw [ultra thick, red, xshift=18cm] (0,0) -- (0,-1) ;
		\draw [ thick,xshift=18cm] (0,-1) -- (1,-2) ;
		\draw [thick,xshift=18cm] (0,-1) -- (-1,-2) ;
		
		\draw[black,fill=black,xshift=21cm] (0,0) circle (.5ex);
		\draw[black,fill=black,xshift=21cm] (0,-1) circle (.5ex);
		\draw[black,fill=black,xshift=21cm] (1,-1) circle (.5ex);
		\draw[black,fill=black,xshift=21cm] (-1,-1) circle (.5ex);
		\draw [thick, xshift=21cm] (0,0) -- (0,-1) ;
		\draw [thick,xshift=21cm] (0,0) -- (1,-1) ;
		\draw [thick,xshift=21cm] (0,0) -- (-1,-1) ;
		
	\end{tikzpicture}
	\caption{All 10 marked ordered trees with 4 nodes.}
\end{figure}
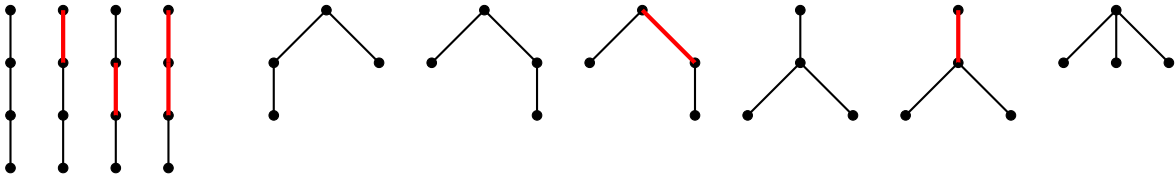

Now we move to a symbolic equation for the marked ordered trees:
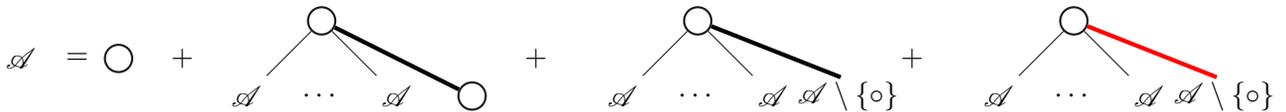
\begin{figure}[h]\small
	\begin{tikzpicture}[scale=1.0,
		s1/.style={circle=10pt,draw=black!90,thick},
		s2/.style={rectangle,draw=black!50,thick},scale=0.5]
		
		\node at ( -3.0,0) { $\mathscr{A}$};
		
		\node at (-1.5,0) { $=$};
		\node(c) at (-0.4,0)[s1]{};
		\node at (1.3,0) {$+$};

		\node(d) at (5,1)[s1]{};
		\node(e) at (3,-1){ $\mathscr{A}$};
		\node(ee) at (5,-1){$\cdots$};
		\node(f) at (7,-1){ $\mathscr{A}$};
		\path [draw,-,black!90] (d) -- (e) node{};%
		\path [draw,-,black!90] (d) -- (f) node{};%
		\node(g) at (9,-1)[s1]{ };
		\path [draw,-,black,ultra thick] (d) -- (g) node {};%
		
		\node[xshift=5cm] at (0.7,0) {$+$};

		\node[xshift=5cm](d) at (5,1)[s1]{};
		\node[xshift=5cm](e) at (3,-1){ $\mathscr{A}$};
		\node[xshift=5cm](ee) at (5,-1){$\cdots$};
		\node[xshift=5cm](f) at (7,-1){ $\mathscr{A}$};
		\path [draw,-,black!90] (d) -- (e) node{};%
		\path [draw,-,black!90] (d) -- (f) node{};%
		\node[xshift=5cm](g) at (9,-1){ $\mathscr{A}\setminus\{\circ\}$};
		\path [draw,-,black,ultra thick] (d) -- (8.8+10,-0.5) node {};%
		
		\node[xshift=10cm] at (0.7,0) {$+$};

		\node[xshift=10cm](d) at (5,1)[s1]{};
		\node[xshift=10cm](e) at (3,-1){ $\mathscr{A}$};
		\node[xshift=10cm](ee) at (5,-1){$\cdots$};
		\node[xshift=10cm](f) at (7,-1){ $\mathscr{A}$};
		\path [draw,-,black!90] (d) -- (e) node{};%
		\path [draw,-,black!90] (d) -- (f) node{};%
		\node[xshift=10cm](g) at (9,-1){ $\mathscr{A}\setminus\{\circ\}$};
		\path [draw,-,black,red,ultra thick] (d) -- (8.8+20,-0.5) node {};%

	\end{tikzpicture}
	\caption{Symbolic equation for marked ordered trees.\\ $\mathscr{A}\cdots\mathscr{A}$ refers to $\ge0$ copies of $\mathscr{A}$.}
\end{figure}

In terms of generating functions,
\begin{equation*}
	A=z+\frac{z}{1-A}z+\frac{z}{1-A}2(A-z),
\end{equation*}
with the solution
\begin{equation*}
	A(z)=\frac{1-z-\sqrt{1-6z+5z^2}}{2}=z+z^2+z^3+3z^3+10z^4+36z^5+\cdots.
\end{equation*}

The importance of this family of trees lies in the bijection to skew Dyck paths, as given in \cite{Deutsch-italy}. One walks around the tree as one usually does and translates it into a Dyck path.
The only difference are the red edges. On the way down, nothing special is to be reported, but on the way up, it is translated into  a skew step $(-1,-1)$. The present author believes that trees are more manageable when it comes to enumeration issues than skew Dyck paths.

The 10 trees of Figure 1 translate as follows:
\begin{equation*}
	\begin{tikzpicture}[scale=0.3]
		\draw (0,0)--(3,3)--(6,0);
	\end{tikzpicture}
	\quad
	\begin{tikzpicture}[scale=0.3]
		\draw (0,0)--(3,3)--(5,1)--(4,0);
	\end{tikzpicture}
	\quad
	\begin{tikzpicture}[scale=0.3]
		\draw (0,0)--(3,3)--(4,2)--(3,1)--(4,0);
	\end{tikzpicture}
	\quad
	\begin{tikzpicture}[scale=0.3]
		\draw (0,0)--(3,3)--(4,2)--(3,1)--(2,0);
	\end{tikzpicture}
	\quad
	\begin{tikzpicture}[scale=0.3]
		\draw (0,0)--(2,2)--(4,0)--(5,1)--(6,0);
	\end{tikzpicture}
\end{equation*}
\begin{equation*}
	\begin{tikzpicture}[scale=0.3]
		\draw (0,0)--(1,1)--(2,0)--(4,2)--(6,0);
	\end{tikzpicture}
	\quad
	\begin{tikzpicture}[scale=0.3]
		\draw (0,0)--(1,1)--(2,0)--(4,2)--(5,1)--(4,0);
	\end{tikzpicture}
	\quad
	\begin{tikzpicture}[scale=0.3]
		\draw (0,0)--(1,1)--(2,2)--(3,1)--(4,2)--(6,0);
	\end{tikzpicture}
	\quad
	\begin{tikzpicture}[scale=0.3]
		\draw (0,0)--(1,1)--(2,2)--(3,1)--(4,2)--(5,1)--(4,0);
	\end{tikzpicture}
	\quad
	\begin{tikzpicture}[scale=0.3]
		\draw (0,0)--(1,1)--(2,0)--(3,1)--(4,0)--(5,1)--(6,0);
	\end{tikzpicture}
\end{equation*}

Skew Dyck paths and a dual version (reading the paths from right-to-left) will be discussed in a later section.

\subsection*{Parameters of marked ordered trees}

There are many parameters, usually considered in the context of ordered trees, that can be considered 
for marked ordered trees. Of course, we cannot be encyclopedic about such parameters. We just consider a few parameters and leave further
analysis to the future.

\subsubsection*{The number of leaves}

To get this, it is most natural to use an additional variable $u$ when translating the symbolic equation, so
that $z^nu^k$ refers to trees with $n$ nodes and $k$ leaves. One obtains
\begin{equation*}
	F=zu+\frac{z}{1-F}\bigl(zu2(F-zu)\bigr),
\end{equation*}
with the solution
\begin{align*}
	F(z,u)&=-z+\frac{zu}2+\frac12-\frac12\sqrt {4{z}^{2}-4z+{z}^{2}{u}^{2}-2zu+1}\\*
	&=zu+{z}^{2}u+ \left( 2u+{u}^{2}\right) {z}^{3}+ \left( 4u+5{u}^{2}+{u}^{3} \right) {z}^{4}+\cdots.
\end{align*}
The factor $4u+5{u}^{2}+{u}^{3}$ corresponds to the 10 trees in Figure 1.

Of interest is also the average number of leaves, when all marked ordered trees of size $n$ are considered to be 
equally likely. For that, we differentiate $F(z,u)$ w. r. t. $u$, followed by $u:=1$, with the result
\begin{equation*}
	\frac{z}{2}+\frac{z-z^2}{2\sqrt{1-6z+5z^2}}=\frac{z}{1-v}, \quad\text{with the usual}\quad z=\frac{v}{1+3v+v^2}.
\end{equation*}
Since $F(z,1)=z(1+v)$, it follows that the average is asymptotic to 
\begin{align*}
	\frac{[z^{n+1}]\frac{z}{1-v}}{[z^{n+1}]z(1+v)}&=\frac{[z^{n}]\frac{1}{1-v}}{[z^{n}](1+v)}=\frac{[z^n]\frac1{\sqrt5}\frac1{\sqrt{1-5z}}}{5^{n+\frac12}\frac1{2\sqrt\pi}n^{3/2}}\\
	&=\frac{\frac{n^{-1/2}}{\Gamma(\frac12)}}{5^{n+\frac12}\frac1{2\sqrt\pi}n^{3/2}}=\frac{n}{10}.
\end{align*}
Note that the corresponding number for ordered trees (unmarked) is $\frac n2$, so we have significantly less leaves here.

\subsubsection*{The height}

As in the seminal paper \cite{BrKnRi72}, we define the height in terms of the longest chain of nodes from the root to a leaf. 
Further, let $p_h=p_h(z)$ be the generating function of marked ordered trees of height $\le h$. From the symbolic equation,
\begin{align*}
	p_{h+1}=z+\frac{z^2}{1-p_h}+\frac{2z(p_h-z)}{1-p_h}=-z+\frac{2z-z^2}{1-p_h},\quad h\ge1,\ p_1=z.
\end{align*}
By some creative guessing, separating numerator and denominator, we find the solution
\begin{equation*}
	p_h=z(1+v)\frac{(1+2v)^{h-1}-v^h(v+2)^{h-1}}{(1+2v)^{h-1}-v^{h+1}(v+2)^{h-1}},
\end{equation*}
which can now be proved  by induction:

We have $p_1=z(1+v)\frac{1-v}{1-v^{2}}=z$. Furthermore, the induction step is best checked using a computer.

The limit of $p_h$ for $h\to\infty$ is $z(1+v)$, the generating function of \emph{all} marked ordered trees, as expected.
Taking differences, we get the generating functions of trees of height $>h$:
\begin{align*}
	z(1+v)&-z(1+v)\frac{(1+2v)^{h-1}-v^h(v+2)^{h-1}}{(1+2v)^{h-1}-v^{h+1}(v+2)^{h-1}}\\
	&=z(1+v)\frac{(1+2v)^{h-1}-v^{h+1}(v+2)^{h-1}-(1+2v)^{h-1}+v^h(v+2)^{h-1}}{(1+2v)^{h-1}-v^{h+1}(v+2)^{h-1}}\\
	&=z(1-v^2)\frac{(v+2)^{h-1}v^h}{(1+2v)^{h-1}-v^{h+1}(v+2)^{h-1}}.
\end{align*}
From this, the average height can be worked out, as in the model paper \cite{HPW}. We sketch the essential steps.
For the average height, one needs
\begin{equation*}
	\sum_{h\ge0}z(1-v^2)\frac{(v+2)^{h-1}v^h}{(1+2v)^{h-1}-v^{h+1}(v+2)^{h-1}}
\end{equation*}
and its behaviour around $v=1$, viz.
\begin{equation*}
	2z(1-v)\sum_{h\ge0}\frac{3^{h-1}v^h}{3^{h-1}-v^{h+1}3^{h-1}}
	\sim2z(1-v)\sum_{h\ge1}\frac{v^h}{1-v^{h}}.
\end{equation*}
The behaviour of the series can be taken straight from \cite{HPW}.

We find there
\begin{equation*}
	\sum_{h\ge1}\frac{v^h}{1-v^{h}}=-\frac{\log(1-v)}{1-v},
\end{equation*}
and further
\begin{equation*}
	\sum_{h\ge0}z(1-v^2)\frac{(v+2)^{h-1}v^h}{(1+2v)^{h-1}-v^{h+1}(v+2)^{h-1}}
	\sim-	2z\log(1-v),
\end{equation*}
so that the coefficient of $z^{n+1}$ is asymptotic to $-2[z^n]\log(1-v)$. Since $1-v\sim \sqrt5\sqrt{1-5z}$,
\begin{equation*}
	-	2z\log(1-v)\sim -2z\log\sqrt{1-5z}= -z\log(1-5z),
\end{equation*}
and the coefficient of $z^{n+1}$ in it is asymptotic to $\frac{5^n}{n}$. This has to be divided (as derived earlier) by
\begin{equation*}
	5^{n+\frac12}\frac{1}{2\sqrt\pi n^{3/2}},
\end{equation*}
with the result
\begin{equation*}
	2\frac{5^n}{n}\frac1{5^{n+\frac12}}\sqrt\pi n^{3/2}=\frac{2}{\sqrt5}\sqrt{\pi n}.
\end{equation*}
Note that the constant in front of $\sqrt{\pi n}$ for ordered trees is $\frac{2}{\sqrt4}=1$, so the average height for marked ordered trees is indeed a bit smaller
thanks to the extra markings.

\section{A bijection between multi-edge trees and 3-coloured Motzkin paths}

Multi-edge trees are like ordered (planar, plane, \dots) trees, but instead of edges there are multiple edges. When drawing such a tree,
instead of drawing, say 5 parallel edges, we just draw one edge and put the number 5 on it as a label. These trees were studied in
\cite{polish, HPW}. For the enumeration, one must count edges. The generating function $F(z)$ satisfies
\begin{equation*}
	F(z)=\sum_{k\ge0}\Big(\frac{z}{1-z}F(z)\Big)^k=\frac1{1-\frac{z}{1-z}F(z)},
\end{equation*}
whence
\begin{equation*}
	F(z)=\frac{1-z-\sqrt{1-6z+5z^2}}{2z}=1+z+3{z}^{2}+10{z}^{3}+36{z}^{4}+137{z}^{5}+543{z}^{6}+\cdots.
\end{equation*}
The coefficients form once again sequence A002212 in \cite{OEIS}.

A Motzkin path consists of up-steps, down-steps, and horizontal steps, see sequence A091965 in \cite{OEIS} and the references given there. As Dyck paths, they start at the origin and end, after $n$ steps again at the $x$-axis, but are not allowed to go below the $x$-axis. A 3-coloured Motzkin path is built as a Motzkin path, but there
are 3 different types of horizontal steps, which we call \emph{red, green, blue}. The generating function $M(z)$ satisfies
\begin{equation*}
	M(z)=1+3zM(z)+z^2M(z)^2=\frac{1-3z-\sqrt{1-6z+5z^2}}{2z^2}, \quad\text{or}\quad F(z)=1+zM(z).
\end{equation*}
So multi-edge trees with $N$ edges (counting the multiplicities) correspond to  3-coloured Motzkin paths of length $N-1$.

The purpose of this note is to describe a bijection. It transforms trees into paths, but all steps are reversible.

\subsection*{The details}

As a first step, the multiplicities will be ignored, and the tree then has only $n$ edges. The standard translation of such tree into the world of Dyck paths,
which is in every book on combinatorics, leads to a Dyck path of length $2n$.
Then the Dyck path will transformed bijectively to a 2-coloured Motzkin path of length $n-1$ (the colours used are red and green).
This transformation plays a prominent role in \cite{Shapiro}, but is most likely much older. I believe that people like Viennot know this for 40 years.
I would be glad to get a proper historic account from the gentle readers.

The last step is then to use the third colour (blue) to deal with the multiplicities.

The first up-step and the last down-step of the Dyck path will be deleted. Then, the remaining $2n-2$ steps are coded pairwise into a 2-Motzkin path of length $n-1$:
\begin{equation*}
	\begin{tikzpicture}[scale=0.3]
		% place nodes
		\path (0,0) node(x1) {\tiny$\bullet$} ;
		\path (1,1) node(x2) {\tiny$\bullet$};
		\path (2,2) node(x3) {\tiny$\bullet$};
		
		\draw (0,0) -- (2,2);

	\end{tikzpicture}
	\raisebox{0.5 em}{$\longrightarrow$}
	\begin{tikzpicture}[scale=0.3]
		% place nodes
		\path (0,0) node(x1) {\tiny$\bullet$} ;
		\path (1,1) node(x2) {\tiny$\bullet$};
		
		\draw (0,0) -- (1,1);

	\end{tikzpicture}
	\qquad
	\begin{tikzpicture}[scale=0.3]
		% place nodes
		\path (0,2) node(x1) {\tiny$\bullet$} ;
		\path (1,1) node(x2) {\tiny$\bullet$};
		\path (2,0) node(x3) {\tiny$\bullet$};
		
		\draw (0,2) -- (2,0);

	\end{tikzpicture}
	\raisebox{0.5 em}{$\longrightarrow$}
	\begin{tikzpicture}[scale=0.3]
		% place nodes
		\path (0,1) node(x1) {\tiny$\bullet$} ;
		\path (1,0) node(x2) {\tiny$\bullet$};
		
		\draw (0,1) -- (1,0);

	\end{tikzpicture}
	\qquad
	\begin{tikzpicture}[scale=0.3]
		% place nodes
		\path (0,0) node(x1) {\tiny$\bullet$} ;
		\path (1,1) node(x2) {\tiny$\bullet$};
		\path (2,0) node(x3) {\tiny$\bullet$};
		
		\draw (0,0) -- (1,1) -- (2,0);

	\end{tikzpicture}
	\raisebox{0.5 em}{$\longrightarrow$}
	\begin{tikzpicture}[scale=0.3]
		% place nodes
		\path (0,0) node[red](x1) {\tiny$\bullet$} ;
		\path (1,0) node[red](x2) {\tiny$\bullet$};
		
		\draw[red, very thick] (0,0) -- (1,0);

	\end{tikzpicture}
	\qquad
	\begin{tikzpicture}[scale=0.3]
		% place nodes
		\path (0,0) node(x1) {\tiny$\bullet$} ;
		\path (1,-1) node(x2) {\tiny$\bullet$};
		\path (2,0) node(x3) {\tiny$\bullet$};
		
		\draw (0,0) -- (1,-1) -- (2,0);

	\end{tikzpicture}
	\raisebox{0.5 em}{$\longrightarrow$}
	\begin{tikzpicture}[scale=0.3]
		% place nodes
		\path (0,0) node[green](x1) {\tiny$\bullet$} ;
		\path (1,0) node[green](x2) {\tiny$\bullet$};
		
		\draw[green, very thick] (0,0) -- (1,0);

	\end{tikzpicture}
\end{equation*}

The last step is to deal with the multiplicities. If an edge is labelled with the number $a$, we will insert $a-1$ horizontal blue steps in the following way:
Since there are currently $n-1$ symbols in the path, we have $n$ possible positions to enter something (in the beginning, in the end, between symbols).
We go through the tree in pre-order, and enter the multiplicities one by one using the blue horizontal steps.

To make this procedure more clear, we prepared a list of 10 multi-edge trees with 3 edges, and the corresponding 3-Motzkin paths of length 2, with intermediate steps 
completely worked out:

%\begin{figure}\label{table1}
\begin{center}
	\begin{table}[h]
		\begin{tabular}{c | c | c  |c}
			\text{Multi-edge tree  }&\text{Dyck path}&\text{2-Motzkin path}&\text{blue edges added}\\
			\hline\hline
			\begin{tikzpicture}[scale=0.5]
				% place nodes
				\path (0,0) node(x1) {\tiny$\bullet$} ;
				\path (0,-1) node(x2) {\tiny$\bullet$};
				\path (0,-2) node(x3) {\tiny$\bullet$};
				\path (0,-3) node(x4) {\tiny$\bullet$};
				\draw (0,0) -- (0,-1)node[pos=0.5,left]{\tiny1} ;
				\draw (0,-1) -- (0,-2)node[pos=0.5,left]{\tiny1} ;
				\draw (0,-2) -- (0,-3)node[pos=0.5,left]{\tiny 1} ;

			\end{tikzpicture}
			& \begin{tikzpicture}[scale=0.45]
				% place nodes
				
				\draw (0,0) -- (3,3) --(6,0);

			\end{tikzpicture}
			& 
			\begin{tikzpicture}[scale=0.45]
				% place nodes
				
				\draw[thick] (0,0) -- (1,1) --(2,0);

			\end{tikzpicture}
			& \begin{tikzpicture}[scale=0.45]
				% place nodes
				
				\draw[thick] (0,0) -- (1,1) --(2,0);

			\end{tikzpicture}\\
			
			\hline

			\begin{tikzpicture}[scale=0.5]
				% place nodes
				\path (0,0) node(x1) {\tiny$\bullet$} ;
				\path (0,-1) node(x2) {\tiny$\bullet$};
				\path (0,-2) node(x3) {\tiny$\bullet$};
				
				\draw (0,0) -- (0,-1)node[pos=0.5,left]{\tiny2} ;
				\draw (0,-1) -- (0,-2)node[pos=0.5,left]{\tiny1} ;

			\end{tikzpicture}
			& \begin{tikzpicture}[scale=0.45]
				% place nodes
				
				\draw (0,0) -- (2,2) --(4,0);

			\end{tikzpicture}
			& \begin{tikzpicture}[scale=0.45]
				% place nodes
				
				\draw [red,thick](0,0) -- (1,0);

			\end{tikzpicture}
			& \begin{tikzpicture}[scale=0.45]
				% place nodes
				\draw [blue,thick](0,0) -- (1,0);
				\draw [red,thick](1,0) -- (2,0);

			\end{tikzpicture}\\
			
			\hline
			\begin{tikzpicture}[scale=0.5]
				% place nodes
				\path (0,0) node(x1) {\tiny$\bullet$} ;
				\path (0,-1) node(x2) {\tiny$\bullet$};
				\path (0,-2) node(x3) {\tiny$\bullet$};
				
				\draw (0,0) -- (0,-1)node[pos=0.5,left]{\tiny1} ;
				\draw (0,-1) -- (0,-2)node[pos=0.5,left]{\tiny2} ;

			\end{tikzpicture}
			& \begin{tikzpicture}[scale=0.45]
				% place nodes
				
				\draw (0,0) -- (2,2) --(4,0);

			\end{tikzpicture}& \begin{tikzpicture}[scale=0.45]
				% place nodes
				
				\draw [red,thick](0,0) -- (1,0);

			\end{tikzpicture}& \begin{tikzpicture}[scale=0.45]
				% place nodes
				\draw [red,thick](0,0) -- (1,0);
				\draw [blue,thick](1,0) -- (2,0);

			\end{tikzpicture}\\
			
			\hline
			\begin{tikzpicture}[scale=0.5]
				% place nodes
				\path (0,0) node(x1) {\tiny$\bullet$} ;
				\path (0,-1) node(x2) {\tiny$\bullet$};

				\draw (0,0) -- (0,-1)node[pos=0.5,left]{\tiny3} ;

			\end{tikzpicture}
			& \begin{tikzpicture}[scale=0.45]
				% place nodes
				
				\draw (0,0) -- (1,1) --(2,0);

			\end{tikzpicture} & & \begin{tikzpicture}[scale=0.45]
				% place nodes
				\draw [blue,thick](0,0) -- (2,0);

			\end{tikzpicture}\\
			
			\hline 
			\begin{tikzpicture}[scale=0.5]
				% place nodes
				\path (0,0) node(x1) {\tiny$\bullet$} ;
				\path (-1,-1) node(x2) {\tiny$\bullet$};
				\path (-1,-2) node(x3) {\tiny$\bullet$};
				\path (1,-1) node(x4) {\tiny$\bullet$};
				\draw (0,0) -- (-1,-1)node[pos=0.3,left]{\tiny1} ;
				\draw (-1,-1) -- (-1,-2)node[pos=0.3,left]{\tiny1} ;
				\draw (0,0) -- (1,-1)node[pos=0.3,right]{\tiny1} ;

			\end{tikzpicture}
			& \begin{tikzpicture}[scale=0.45]
				% place nodes
				
				\draw (0,0) -- (2,2) --(4,0)--(5,1)--(6,0);

			\end{tikzpicture}& \begin{tikzpicture}[scale=0.45]
				% place nodes
				
				\draw [red,thick](0,0) -- (1,0);
				\draw [green,thick](1,0) -- (2,0);
				
			\end{tikzpicture}& \begin{tikzpicture}[scale=0.45]
				% place nodes
				
				\draw [red,thick](0,0) -- (1,0);
				\draw [green,thick](1,0) -- (2,0);
				
			\end{tikzpicture}\\
			
			\hline
			\begin{tikzpicture}[scale=0.5]
				% place nodes
				\path (0,0) node(x1) {\tiny$\bullet$} ;
				\path (-1,-1) node(x2) {\tiny$\bullet$};
				
				\path (1,-1) node(x4) {\tiny$\bullet$};
				\draw (0,0) -- (-1,-1)node[pos=0.3,left]{\tiny2} ;
				
				\draw (0,0) -- (1,-1)node[pos=0.3,right]{\tiny1} ;

			\end{tikzpicture}
			& \begin{tikzpicture}[scale=0.45]
				% place nodes
				
				\draw (0,0) -- (1,1) --(2,0)--(3,1)--(4,0);

			\end{tikzpicture}& \begin{tikzpicture}[scale=0.45]
				% place nodes

				\draw [green,thick](0,0) -- (1,0);
				
			\end{tikzpicture}			& \begin{tikzpicture}[scale=0.45]
				% place nodes
				
				\draw [blue,thick](0,0) -- (1,0);
				\draw [green,thick](1,0) -- (2,0);
				
			\end{tikzpicture}\\
			
			\hline
			\begin{tikzpicture}[scale=0.5]
				% place nodes
				\path (0,0) node(x1) {\tiny$\bullet$} ;
				\path (-1,-1) node(x2) {\tiny$\bullet$};
				
				\path (1,-1) node(x4) {\tiny$\bullet$};
				\draw (0,0) -- (-1,-1)node[pos=0.3,left]{\tiny1} ;
				
				\draw (0,0) -- (1,-1)node[pos=0.3,right]{\tiny2} ;

			\end{tikzpicture}
			& \begin{tikzpicture}[scale=0.45]
				% place nodes
				
				\draw (0,0) -- (1,1) --(2,0)--(3,1)--(4,0);

			\end{tikzpicture}& \begin{tikzpicture}[scale=0.45]
				% place nodes

				\draw [green,thick](0,0) -- (1,0);
				
			\end{tikzpicture}& \begin{tikzpicture}[scale=0.45]
				% place nodes
				
				\draw [green,thick](0,0) -- (1,0);
				\draw [blue,thick](1,0) -- (2,0);
				
			\end{tikzpicture}\\
			
			\hline 
			\begin{tikzpicture}[scale=0.5]
				% place nodes
				\path (-1,0) node(x1) {\tiny$\bullet $} ;
				\path (-1,1) node(x2) {\tiny$\bullet$};
				\path (-2,2) node(x3) {\tiny$\bullet$};
				\path (-3,1) node(x4) {\tiny$\bullet$};
				\draw (-1,0) -- (-1,1)node[pos=0.7,right]{\tiny1} ;
				\draw (-1,1) -- (-2,2)node[pos=0.7,right]{\tiny1} ;
				\draw (-2,2) -- (-3,1)node[pos=0.3,left]{\tiny1} ;

			\end{tikzpicture}
			& \begin{tikzpicture}[scale=0.45]
				% place nodes
				
				\draw (0,0) -- (1,1) --(2,0)--(4,2)--(6,0);

			\end{tikzpicture}& \begin{tikzpicture}[scale=0.45]
				% place nodes

				\draw [green,thick](0,0) -- (1,0);
				\draw[red,thick](1,0)--(2,0);			
			\end{tikzpicture}& \begin{tikzpicture}[scale=0.45]
				% place nodes

				\draw [green,thick](0,0) -- (1,0);
				\draw[red,thick](1,0)--(2,0);			
			\end{tikzpicture}\\
			
			\hline
			\begin{tikzpicture}[scale=0.5]
				% place nodes
				\path (0,0) node(x1) {\tiny$\bullet$} ;
				\path (-1,-1) node(x2) {\tiny$\bullet$};
				\path (0,-1) node(x3) {\tiny$\bullet$};
				\path (1,-1) node(x4) {\tiny$\bullet$};
				\draw (0,0) -- (-1,-1)node[pos=0.3,left]{\tiny1} ;
				\draw (0,0) -- (0,-1)node[pos=0.6]{\tiny\;\;1} ;
				\draw (0,0) -- (1,-1)node[pos=0.3,right]{\tiny1} ;

			\end{tikzpicture}
			& \begin{tikzpicture}[scale=0.45]
				% place nodes
				
				\draw (0,0) -- (1,1) --(2,0)--(3,1)--(4,0)--(5,1)--(6,0);

			\end{tikzpicture}& \begin{tikzpicture}[scale=0.45]
				% place nodes

				\draw [green,thick](0,0) -- (1,0);
				\draw [green,thick](1,0) -- (2,0);
				
			\end{tikzpicture}& \begin{tikzpicture}[scale=0.45]
				% place nodes

				\draw [green,thick](0,0) -- (1,0);
				\draw [green,thick](1,0) -- (2,0);
				
			\end{tikzpicture}\\
			
			\hline 
			\begin{tikzpicture}[scale=0.5]
				% place nodes
				\path (0,0) node(x1) {\tiny$\bullet$} ;
				\path (0,-1) node(x2) {\tiny$\bullet$};
				\path (-1,-2) node(x3) {\tiny$\bullet$};
				\path (1,-2) node(x4) {\tiny$\bullet$};
				\draw (0,0) -- (0,-1)node[pos=0.5,left]{\tiny1} ;
				\draw (0,-1) -- (1,-2)node[pos=0.5,right]{\tiny1} ;
				\draw (0,-1) -- (-1,-2)node[pos=0.5,left]{\tiny 1} ;

			\end{tikzpicture}
			& \begin{tikzpicture}[scale=0.45]
				% place nodes
				
				\draw (0,0) -- (2,2) --(3,1)--(4,2)--(6,0);

			\end{tikzpicture}& \begin{tikzpicture}[scale=0.45]
				% place nodes

				\draw [red,thick](0,0) -- (1,0);
				\draw [red,thick](1,0) -- (2,0);
				
			\end{tikzpicture}& \begin{tikzpicture}[scale=0.45]
				% place nodes

				\draw [red,thick](0,0) -- (1,0);
				\draw [red,thick](1,0) -- (2,0);
				
			\end{tikzpicture}\\

		\end{tabular}
		
		\caption{First row is a multi-edge tree with 3 edges, second row is the standard Dyck path (multiplicities ignored), third row is cutting off first and last step, and then translated pairs of steps, fourth row is inserting blued horizontal edges, according to multiplicities.}
	\end{table}	
\end{center}
%\end{figure}

\subsection*{Connecting unary-binary trees with multi-edge trees} 

This is not too difficult: We start from multi-edge trees, and ignore the multiplicities at the moment. Then we apply the classical rotation correspondence (also called: natural correspondence).
Then we add vertical edges, if the multiplicity is higher than 1. To be precise, if there is a node, and an edge with multiplicity $a$ leads to it from the top, we insert $a-1$ extra nodes in a chain on the top, and connect them with unary branches. The following example with 10 objects will help to understand this procedure.
After that, all the structures studied in this paper are connected with bijections.

\begin{center}
	\begin{table}[h]
		\begin{tabular}{c | c | c  |c}
			\text{Multi-edge tree  }&\text{Binary tree (rotation)}&\text{vertical edges added}\\
			\hline\hline
			\begin{tikzpicture}[scale=0.5]
				% place nodes
				\path (0,0) node(x1) {\tiny$\bullet$} ;
				\path (0,-1) node(x2) {\tiny$\bullet$};
				\path (0,-2) node(x3) {\tiny$\bullet$};
				\path (0,-3) node(x4) {\tiny$\bullet$};
				\draw (0,0) -- (0,-1)node[pos=0.5,left]{\tiny1} ;
				\draw (0,-1) -- (0,-2)node[pos=0.5,left]{\tiny1} ;
				\draw (0,-2) -- (0,-3)node[pos=0.5,left]{\tiny 1} ;

			\end{tikzpicture}

			& 			\begin{tikzpicture}[scale=0.5]
				% place nodes
				\path (0,0) node(x1) {\tiny$\bullet$} ;
				\path (-1,-1) node(x2) {\tiny$\bullet$};
				\path (-2,-2) node(x3) {\tiny$\bullet$};
				
				\draw (0,0) -- (-1,-1);
				\draw (-1,-1) -- (-2,-2) ;

			\end{tikzpicture}
			
			& \begin{tikzpicture}[scale=0.5]
				% place nodes
				\path (0,0) node(x1) {\tiny$\bullet$} ;
				\path (-1,-1) node(x2) {\tiny$\bullet$};
				\path (-2,-2) node(x3) {\tiny$\bullet$};
				
				\draw (0,0) -- (-1,-1);
				\draw (-1,-1) -- (-2,-2) ;

			\end{tikzpicture}
			\\
			\hline 
			\begin{tikzpicture}[scale=0.5]
				% place nodes
				\path (0,0) node(x1) {\tiny$\bullet$} ;
				\path (0,-1) node(x2) {\tiny$\bullet$};
				\path (0,-2) node(x3) {\tiny$\bullet$};
				
				\draw (0,0) -- (0,-1)node[pos=0.5,left]{\tiny2} ;
				\draw (0,-1) -- (0,-2)node[pos=0.5,left]{\tiny1} ;

			\end{tikzpicture}
			& \begin{tikzpicture}[scale=0.5]
				% place nodes
				\path (0,0) node(x1) {\tiny$\bullet$} ;
				\path (-1,-1) node(x2) {\tiny$\bullet$};
				%\path (-2,-2) node(x3) {\tiny$\bullet$};
				
				\draw (0,0) -- (-1,-1);
				%\draw (-1,-1) -- (-2,-2) ;

			\end{tikzpicture}
			&\begin{tikzpicture}[scale=0.5]
				% place nodes
				\path (0,0) node(x1) {\tiny$\bullet$} ;
				\path (-1,-1) node(x2) {\tiny$\bullet$};
				\path (0,1) node(x3) {\tiny$\bullet$};
				
				\draw (0,0) -- (-1,-1);
				\draw (0,0) -- (0,1) ;

			\end{tikzpicture}
			\\
			\hline 
			\begin{tikzpicture}[scale=0.5]
				% place nodes
				\path (0,0) node(x1) {\tiny$\bullet$} ;
				\path (0,-1) node(x2) {\tiny$\bullet$};
				\path (0,-2) node(x3) {\tiny$\bullet$};
				
				\draw (0,0) -- (0,-1)node[pos=0.5,left]{\tiny1} ;
				\draw (0,-1) -- (0,-2)node[pos=0.5,left]{\tiny2} ;

			\end{tikzpicture}
			& \begin{tikzpicture}[scale=0.5]
				% place nodes
				\path (0,0) node(x1) {\tiny$\bullet$} ;
				\path (-1,-1) node(x2) {\tiny$\bullet$};
				%\path (-2,-2) node(x3) {\tiny$\bullet$};
				
				\draw (0,0) -- (-1,-1);
				%\draw (-1,-1) -- (-2,-2) ;

			\end{tikzpicture}
			&\begin{tikzpicture}[scale=0.5]
				% place nodes
				\path (-1,0) node(x1) {\tiny$\bullet$} ;
				\path (-1,-1) node(x2) {\tiny$\bullet$};
				\path (0,1) node(x3) {\tiny$\bullet$};
				
				\draw (-1,0) -- (-1,-1);
				\draw (-1,0) -- (0,1) ;

			\end{tikzpicture}
			\\
			\hline
			\begin{tikzpicture}[scale=0.5]
				% place nodes
				\path (0,0) node(x1) {\tiny$\bullet$} ;
				\path (0,-1) node(x2) {\tiny$\bullet$};

				\draw (0,0) -- (0,-1)node[pos=0.5,left]{\tiny3} ;

			\end{tikzpicture}
			& \begin{tikzpicture}[scale=0.5]
				% place nodes
				\path (0,0) node(x1) {\tiny$\bullet$} ;

			\end{tikzpicture}
			&\begin{tikzpicture}[scale=0.5]
				% place nodes
				\path (0,0) node(x1) {\tiny$\bullet$} ;
				\path (0,-1) node(x2) {\tiny$\bullet$};
				\path (0,1) node(x3) {\tiny$\bullet$};
				
				\draw (0,0) -- (0,-1);
				\draw (0,0) -- (0,1) ;

			\end{tikzpicture}
			\\
			
			\hline 
			\begin{tikzpicture}[scale=0.5]
				% place nodes
				\path (0,0) node(x1) {\tiny$\bullet$} ;
				\path (-1,-1) node(x2) {\tiny$\bullet$};
				\path (-1,-2) node(x3) {\tiny$\bullet$};
				\path (1,-1) node(x4) {\tiny$\bullet$};
				\draw (0,0) -- (-1,-1)node[pos=0.3,left]{\tiny1} ;
				\draw (-1,-1) -- (-1,-2)node[pos=0.3,left]{\tiny1} ;
				\draw (0,0) -- (1,-1)node[pos=0.3,right]{\tiny1} ;

			\end{tikzpicture}
			& \begin{tikzpicture}[scale=0.5]
				% place nodes
				\path (0,-2) node(x1) {\tiny$\bullet$} ;
				\path (-1,-1) node(x2) {\tiny$\bullet$};
				\path (-2,-2) node(x3) {\tiny$\bullet$};
				
				\draw (0,-2) -- (-1,-1);
				\draw (-1,-1) -- (-2,-2) ;

			\end{tikzpicture}
			& \begin{tikzpicture}[scale=0.5]
				% place nodes
				\path (0,-2) node(x1) {\tiny$\bullet$} ;
				\path (-1,-1) node(x2) {\tiny$\bullet$};
				\path (-2,-2) node(x3) {\tiny$\bullet$};
				
				\draw (0,-2) -- (-1,-1);
				\draw (-1,-1) -- (-2,-2) ;

			\end{tikzpicture}
			\\
			\hline
			\begin{tikzpicture}[scale=0.5]
				% place nodes
				\path (0,0) node(x1) {\tiny$\bullet$} ;
				\path (-1,-1) node(x2) {\tiny$\bullet$};
				
				\path (1,-1) node(x4) {\tiny$\bullet$};
				\draw (0,0) -- (-1,-1)node[pos=0.3,left]{\tiny2} ;
				
				\draw (0,0) -- (1,-1)node[pos=0.3,right]{\tiny1} ;

			\end{tikzpicture}
			& \begin{tikzpicture}[scale=0.5]
				% place nodes
				\path (0,-2) node(x1) {\tiny$\bullet$} ;
				\path (-1,-1) node(x2) {\tiny$\bullet$};
				%\path (-2,-2) node(x3) {\tiny$\bullet$};
				
				\draw (0,-2) -- (-1,-1);
				%\draw (-1,-1) -- (-2,-2) ;

			\end{tikzpicture}
			&  	\begin{tikzpicture}[scale=0.5]
				% place nodes
				\path (0,0) node(x1) {\tiny$\bullet$} ;
				\path (1,-1) node(x2) {\tiny$\bullet$};
				\path (0,1) node(x3) {\tiny$\bullet$};
				
				\draw (0,0) -- (1,-1);
				\draw (0,0) -- (0,1) ;

			\end{tikzpicture}
			\\
			\hline
			\begin{tikzpicture}[scale=0.5]
				% place nodes
				\path (0,0) node(x1) {\tiny$\bullet$} ;
				\path (-1,-1) node(x2) {\tiny$\bullet$};
				
				\path (1,-1) node(x4) {\tiny$\bullet$};
				\draw (0,0) -- (-1,-1)node[pos=0.3,left]{\tiny1} ;
				
				\draw (0,0) -- (1,-1)node[pos=0.3,right]{\tiny2} ;

			\end{tikzpicture}
			& \begin{tikzpicture}[scale=0.5]
				% place nodes
				\path (0,-2) node(x1) {\tiny$\bullet$} ;
				\path (-1,-1) node(x2) {\tiny$\bullet$};
				%\path (-2,-2) node(x3) {\tiny$\bullet$};
				
				\draw (0,-2) -- (-1,-1);
				%\draw (-1,-1) -- (-2,-2) ;

			\end{tikzpicture}
			
			&  \begin{tikzpicture}[scale=0.5]
				% place nodes
				\path (0,0) node(x1) {\tiny$\bullet$} ;
				\path (0,-1) node(x2) {\tiny$\bullet$};
				\path (-1,1) node(x3) {\tiny$\bullet$};
				
				\draw (0,0) -- (-1,1);
				\draw (0,0) -- (0,-1) ;

			\end{tikzpicture}
			
			\\
			\hline 
			\begin{tikzpicture}[scale=0.5]
				% place nodes
				\path (-1,0) node(x1) {\tiny$\bullet $} ;
				\path (-1,1) node(x2) {\tiny$\bullet$};
				\path (-2,2) node(x3) {\tiny$\bullet$};
				\path (-3,1) node(x4) {\tiny$\bullet$};
				\draw (-1,0) -- (-1,1)node[pos=0.7,right]{\tiny1} ;
				\draw (-1,1) -- (-2,2)node[pos=0.7,right]{\tiny1} ;
				\draw (-2,2) -- (-3,1)node[pos=0.3,left]{\tiny1} ;

			\end{tikzpicture}
			& \begin{tikzpicture}[scale=0.5]
				% place nodes
				\path (-2,0) node(x1) {\tiny$\bullet$} ;
				\path (-1,-1) node(x2) {\tiny$\bullet$};
				\path (-2,-2) node(x3) {\tiny$\bullet$};
				
				\draw (-2,0) -- (-1,-1);
				\draw (-1,-1) -- (-2,-2) ;

			\end{tikzpicture}
			& \begin{tikzpicture}[scale=0.5]
				% place nodes
				\path (-2,0) node(x1) {\tiny$\bullet$} ;
				\path (-1,-1) node(x2) {\tiny$\bullet$};
				\path (-2,-2) node(x3) {\tiny$\bullet$};
				
				\draw (-2,0) -- (-1,-1);
				\draw (-1,-1) -- (-2,-2) ;

			\end{tikzpicture}
			\\
			\hline
			\begin{tikzpicture}[scale=0.5]
				% place nodes
				\path (0,0) node(x1) {\tiny$\bullet$} ;
				\path (-1,-1) node(x2) {\tiny$\bullet$};
				\path (0,-1) node(x3) {\tiny$\bullet$};
				\path (1,-1) node(x4) {\tiny$\bullet$};
				\draw (0,0) -- (-1,-1)node[pos=0.3,left]{\tiny1} ;
				\draw (0,0) -- (0,-1)node[pos=0.6]{\tiny\;\;1} ;
				\draw (0,0) -- (1,-1)node[pos=0.3,right]{\tiny1} ;

			\end{tikzpicture}&
			\begin{tikzpicture}[scale=0.5]
				% place nodes
				\path (-3,3) node(x1) {\tiny$\bullet$} ;
				\path (-1,1) node(x2) {\tiny$\bullet$};
				\path (-2,2) node(x3) {\tiny$\bullet$};
				
				\draw (-2,2) -- (-1,1);
				\draw (-3,3) -- (-2,2) ;

			\end{tikzpicture}
			&\begin{tikzpicture}[scale=0.5]
				% place nodes
				\path (-3,3) node(x1) {\tiny$\bullet$} ;
				\path (-1,1) node(x2) {\tiny$\bullet$};
				\path (-2,2) node(x3) {\tiny$\bullet$};
				
				\draw (-2,2) -- (-1,1);
				\draw (-3,3) -- (-2,2) ;

			\end{tikzpicture}
			\\
			
			\hline 
			\begin{tikzpicture}[scale=0.5]
				% place nodes
				\path (0,0) node(x1) {\tiny$\bullet$} ;
				\path (0,-1) node(x2) {\tiny$\bullet$};
				\path (-1,-2) node(x3) {\tiny$\bullet$};
				\path (1,-2) node(x4) {\tiny$\bullet$};
				\draw (0,0) -- (0,-1)node[pos=0.5,left]{\tiny1} ;
				\draw (0,-1) -- (1,-2)node[pos=0.5,right]{\tiny1} ;
				\draw (0,-1) -- (-1,-2)node[pos=0.5,left]{\tiny 1} ;

			\end{tikzpicture}
			& \begin{tikzpicture}[scale=0.5]
				% place nodes
				\path (-2,0) node(x1) {\tiny$\bullet$} ;
				\path (-1,1) node(x2) {\tiny$\bullet$};
				\path (-1,-1) node(x3) {\tiny$\bullet$};
				
				\draw (-2,0) -- (-1,-1);
				\draw (-2,0) -- (-1,1) ;

			\end{tikzpicture}
			& \begin{tikzpicture}[scale=0.5]
				% place nodes
				\path (-2,0) node(x1) {\tiny$\bullet$} ;
				\path (-1,1) node(x2) {\tiny$\bullet$};
				\path (-1,-1) node(x3) {\tiny$\bullet$};
				
				\draw (-2,0) -- (-1,-1);
				\draw (-2,0) -- (-1,1) ;

			\end{tikzpicture}\\
			
			%\hline

		\end{tabular}
		
		\caption{First row is a multi-edge tree with 3 edges, second row the corresponding binary tree, according to the classical rotation correspondence, ignoring the unary branches.
			Third row is inserting extra horizontal edges when the multiplicities are higher than 1.}
	\end{table}	
\end{center}

\section{The combinatorics of skew Dyck paths}

Skew Dyck are a variation of Dyck paths, where additionally to steps $(1,1)$ and $(1,-1)$ a south-west step $(-1,-1)$ is also allowed, provided that the path
does not intersect itself. Otherwise, like for Dyck path, it must never go below the $x$-axis and end eventually (after $2n$ steps) on the $x$-axis.
Here are a few references: \cite{Deutsch-italy, KimStanley, Baril-neu, Prodinger-hex}. The enumerating sequence is
\begin{equation*}
	1, 1, 3, 10, 36, 137, 543, 2219, 9285, 39587, 171369, 751236, 3328218, 14878455,\dots,
\end{equation*}
which is A002212 in \cite{OEIS}.

Skew Dyck appeared very briefly in a previous section; here we want to give a more thorough analysis of them, using generating functions and
the kernel method. 
Here is a list of the 10 skew paths consisting of 6 steps:

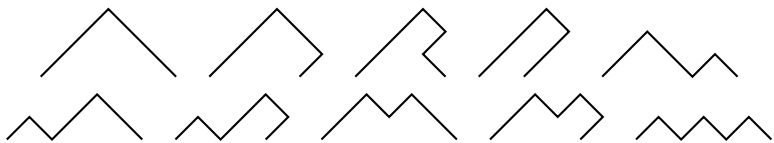
\begin{figure}[h]
	\begin{equation*}
		\begin{tikzpicture}[scale=0.3]
			\draw [thick](0,0)--(3,3)--(6,0);
		\end{tikzpicture}
		\quad
		\begin{tikzpicture}[scale=0.3]
			\draw [thick](0,0)--(3,3)--(5,1)--(4,0);
		\end{tikzpicture}
		\quad
		\begin{tikzpicture}[scale=0.3]
			\draw [thick](0,0)--(3,3)--(4,2)--(3,1)--(4,0);
		\end{tikzpicture}
		\quad
		\begin{tikzpicture}[scale=0.3]
			\draw [thick](0,0)--(3,3)--(4,2)--(3,1)--(2,0);
		\end{tikzpicture}
		\quad
		\begin{tikzpicture}[scale=0.3]
			\draw [thick](0,0)--(2,2)--(4,0)--(5,1)--(6,0);
		\end{tikzpicture}
	\end{equation*}
	\begin{equation*}
		\begin{tikzpicture}[scale=0.3]
			\draw [thick](0,0)--(1,1)--(2,0)--(4,2)--(6,0);
		\end{tikzpicture}
		\quad
		\begin{tikzpicture}[scale=0.3]
			\draw [thick](0,0)--(1,1)--(2,0)--(4,2)--(5,1)--(4,0);
		\end{tikzpicture}
		\quad
		\begin{tikzpicture}[scale=0.3]
			\draw [thick](0,0)--(1,1)--(2,2)--(3,1)--(4,2)--(6,0);
		\end{tikzpicture}
		\quad
		\begin{tikzpicture}[scale=0.3]
			\draw [thick](0,0)--(1,1)--(2,2)--(3,1)--(4,2)--(5,1)--(4,0);
		\end{tikzpicture}
		\quad
		\begin{tikzpicture}[scale=0.3]
			\draw [thick](0,0)--(1,1)--(2,0)--(3,1)--(4,0)--(5,1)--(6,0);
		\end{tikzpicture}
	\end{equation*}
	\caption{All 10 skew Dyck paths of length 6 (consisting of 6 steps).}
\end{figure}

We prefer to work with the equivalent model (resembling more traditional Dyck paths) where
we replace each step $(-1,-1)$ by $(1,-1)$ but label it red. Here is the list of the 10 paths again (Figure 2):

\begin{figure}[h]
	\begin{equation*}
		\begin{tikzpicture}[scale=0.3]
			\draw [thick](0,0)--(3,3)--(6,0);
		\end{tikzpicture}
		\quad
		\begin{tikzpicture}[scale=0.3]
			\draw [thick](0,0)--(3,3)--(5,1);
			\draw [thick,red](5,1)--(6,0);
		\end{tikzpicture}
		\quad
		\begin{tikzpicture}[scale=0.3]
			\draw [thick](0,0)--(3,3)--(4,2);
			\draw[red,thick] (4,2)--(5,1);
			\draw [thick](5,1)--(6,0);
		\end{tikzpicture}
		\quad
		\begin{tikzpicture}[scale=0.3]
			\draw [thick](0,0)--(3,3)--(4,2);
			\draw[red,thick](4,2)--(6,0);
		\end{tikzpicture}
		\quad
		\begin{tikzpicture}[scale=0.3]
			\draw [thick](0,0)--(2,2)--(4,0)--(5,1)--(6,0);
		\end{tikzpicture}
	\end{equation*}
	\begin{equation*}
		\begin{tikzpicture}[scale=0.3]
			\draw [thick](0,0)--(1,1)--(2,0)--(4,2)--(6,0);
		\end{tikzpicture}
		\quad
		\begin{tikzpicture}[scale=0.3]
			\draw [thick](0,0)--(1,1)--(2,0)--(4,2)--(5,1);
			\draw[red,thick] (5,1)--(6,0);
		\end{tikzpicture}
		\quad
		\begin{tikzpicture}[scale=0.3]
			\draw [thick](0,0)--(1,1)--(2,2)--(3,1)--(4,2)--(6,0);
		\end{tikzpicture}
		\quad
		\begin{tikzpicture}[scale=0.3]
			\draw [thick](0,0)--(1,1)--(2,2)--(3,1)--(4,2)--(5,1);
			\draw[red,thick] (5,1)--(6,0);
		\end{tikzpicture}
		\quad
		\begin{tikzpicture}[scale=0.3]
			\draw [thick](0,0)--(1,1)--(2,0)--(3,1)--(4,0)--(5,1)--(6,0);
		\end{tikzpicture}
	\end{equation*}
	\caption{The 10 paths redrawn, with red south-east edges instead of south-west edges.}
\end{figure}
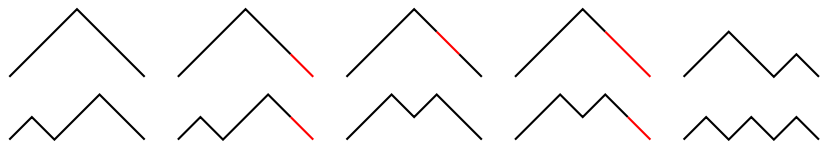

The rules to generate such decorated Dyck paths are: each edge $(1,-1)$ may be  black or red, but
\begin{tikzpicture}[scale=0.3]\draw [thick](0,0)--(1,1); \draw [red,thick] (1,1)--(2,0);\end{tikzpicture}
and
\begin{tikzpicture}[scale=0.3] \draw [red,thick] (0,1)--(1,0);\draw [thick](1,0)--(2,1);\end{tikzpicture}
are forbidden.

Our interest is in particular in \emph{partial} decorated Dyck paths, ending at level $j$, for fixed $j\ge0$;
the instance $j=0$ is the classical case.

The analysis of partial skew Dyck paths was recently  started in \cite{Baril-neu} (using the notion `prefix of a skew Dyck path') using
Riordan arrays instead of our kernel method. The latter  gives us \emph{bivariate} generating functions,
from which it is easier to draw conclusions. Two variables, $z$ and $u$, are used, where $z$ marks the length
of the path and $j$ marks the end-level. We briefly mention that one can, using a third variable $w$, also
count the number of red edges.

Again, once all generating functions are explicitly known, many corollories can be derived in a standard fashion.
We only do this in a few instances. But we would like to emphasize that the substitution
\begin{equation*}
	x=\frac{v}{1+3v+v^2},
\end{equation*}
which was used in \cite{HPW, Prodinger-hex} allows to write \emph{explicit enumerations}, using
the notion of a (weighted) trinomial coefficient:
\begin{equation*}
	\binom{n;1,3,1}{k}:=[t^k](1+3t+t^2)^n.
\end{equation*}

The second part of this section  deals with a dual version, where the paths are read from right to left.

\subsection*{Generating functions and the kernel method}

We catch the essence of a decorated Dyck path using a state-diagram:

\begin{figure}[h]

	\begin{center}
		\begin{tikzpicture}[scale=1.5]
			\draw (0,0) circle (0.1cm);
			\fill (0,0) circle (0.1cm);
			
			\foreach \x in {0,1,2,3,4,5,6,7,8}
			{
				\draw (\x,0) circle (0.05cm);
				\fill (\x,0) circle (0.05cm);
			}
			
			\foreach \x in {0,1,2,3,4,5,6,7,8}
			{
				\draw (\x,-1) circle (0.05cm);
				\fill (\x,-1) circle (0.05cm);
			}
			
			\foreach \x in {0,1,2,3,4,5,6,7,8}
			{
				\draw (\x,-2) circle (0.05cm);
				\fill (\x,-2) circle (0.05cm);
			}
			
			\foreach \x in {0,1,2,3,4,5,6,7}
			{
				\draw[ thick,-stealth] (\x,0) -- (\x+1,0);
				
			}

			\foreach \x in {0,1,2,3,4,5,6,7}
			{
				\draw[ thick,<->] (\x+1,0) -- (\x,-1);
				
			}
			
			\foreach \x in {0,1,2,3,4,5,6,7}
			{
				\draw[ thick,-stealth] (\x+1,-1) -- (\x,-1);
				
			}
			\foreach \x in {0,1,2,3,4,5,6,7}
			{
				\draw[ thick,-stealth,red] (\x+1,-1) -- (\x,-2);
				
			}
			
			\foreach \x in {0,1,2,3,4,5,6,7}
			{
				\draw[ thick,-stealth,red] (\x+1,-2) -- (\x,-2);
				
			}
			
			\foreach \x in {0,1,2,3,4,5,6,7}
			{
				\draw[ thick,-stealth] (\x+1,-2) -- (\x,-1);
				
			}

		\end{tikzpicture}
	\end{center}
	\caption{Three layers of states according to the type of steps leading to them (up, down-black, down-red).}
\end{figure}
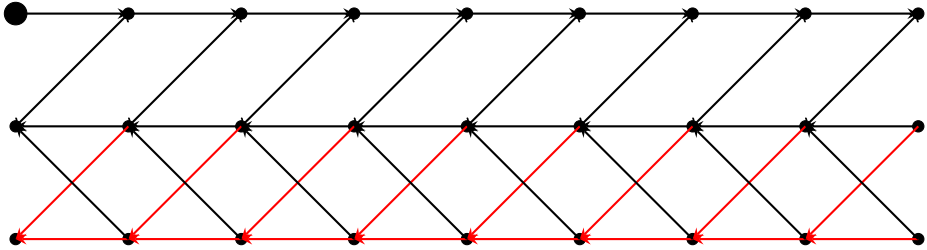
It has three types of states, with $j$ ranging from 0 to infinity; in the drawing, only $j=0..8$ is shown. The first
layer of states refers to an up-step leading to a state, the second layer refers to a black down-step leading to a state
and the third layer refers to a red down-step leading to a state. We will work out generating functions describing all paths
leading to a particular state. We will use the notations $f_j,g_j,h_j$ for the three respective layers, from top to bottom.
Note that the syntactic rules of forbidden patterns
\begin{tikzpicture}[scale=0.3]\draw [thick](0,0)--(1,1); \draw [red,thick] (1,1)--(2,0);\end{tikzpicture}
and
\begin{tikzpicture}[scale=0.3] \draw [red,thick] (0,1)--(1,0);\draw [thick](1,0)--(2,1);\end{tikzpicture}
can be clearly seen from the picture. The functions depend on the variable $z$ (marking the number of steps), but mostly
we just write $f_j$ instead of $f_j(z)$, etc.

The following recursions can be read off immediately from the diagram:
\begin{gather*}
	f_0=1,\quad f_{i+1}=zf_i+zg_i,\quad i\ge0,\\
	g_i=zf_{i+1}+zg_{i+1}+zh_{i+1},\quad i\ge0,\\
	h_i=zh_{i+1}+zg_{i+1},\quad i\ge0.
\end{gather*}
And now it is time to introduce the promised \emph{bivariate} generating functions: 
\begin{equation*}
	F(z,u)=\sum_{i\ge0}f_i(z)u^i,\quad
	G(z,u)=\sum_{i\ge0}g_i(z)u^i,\quad
	H(z,u)=\sum_{i\ge0}h_i(z)u^i.
\end{equation*}
Again, often we just write $F(u)$ instead of $F(z,u)$ and treat $z$ as a `silent' variable. Summing the recursions leads to
\begin{align*}
	\sum_{i\ge0}u^if_{i+1}&=\sum_{i\ge0}u^izf_i+\sum_{i\ge0}u^izg_i,\\
	\sum_{i\ge0}u^ig_i&=\sum_{i\ge0}u^izf_{i+1}+\sum_{i\ge0}u^izg_{i+1}+\sum_{i\ge0}u^izh_{i+1},\\
	\sum_{i\ge0}u^ih_i&=\sum_{i\ge0}u^izh_{i+1}+\sum_{i\ge0}u^izg_{i+1}.
\end{align*}
This can be rewritten as
\begin{align*}
	\frac1u(F(u)-1)&=zF(u)+zG(u),\\*
	G(u)&=\frac zu(F(u)-1)+\frac zu(G(u)-G(0))+\frac zu(H(u)-H(0)),\\*
	H(u)&=	\frac zu(G(u)-G(0))+\frac zu(H(u)-H(0)).
\end{align*}
This is a typical application of the kernel method. For a gentle example-driven introduction to the kernel method, see \cite{Prodinger-kernel}. First,
\begin{align*}
	F(u)&=\frac{z^2uG(0)+z^2uH(0)+z^2u-u-z^3+2z}{-{z}^{3}-u+2z+z{u}^{2}-{z}^{2}u},\\
	G(u)&=\frac{z(H(0)-uzH(0)+z^2+G(0)-zuG(0)-zu)}{-{z}^{3}-u+2z+z{u}^{2}-{z}^{2}u},\\
	H(u)&=\frac{z(-uzH(0)-z^2-zuG(0)+G(0)-z^2H(0)+H(0)-z^2G(0))}{-{z}^{3}-u+2z+z{u}^{2}-{z}^{2}u}.
\end{align*}
The denominator factors as $z(u-r_1)(u-r_2)$, with
\begin{equation*}
	r_1=\frac{1+z^2+\sqrt{1-6z^2+5z^4}}{2z},\quad r_2=\frac{1+z^2-\sqrt{1-6z^2+5z^4}}{2z}.
\end{equation*}
Note that $r_1r_2=2-z^2$. 
Since the factor $u-r_2$ in the denominator is ``bad,'' it must also cancel in the numerators. From this
we conclude as a first step
\begin{equation*}
	G(0) =  \frac{1-2z^2H(0)-3z^2-\sqrt{1-6z^2+5z^4}}{2z^2},
\end{equation*}
and by further simplification
\begin{equation*}
	H(0)=\frac{1-4z^2+z^4+(z^2-1)\sqrt{1-6z^2+5z^4}}{2-z^2}.
\end{equation*}
Thus (with $W=\sqrt{1-6z^2+5z^4}=\sqrt{(1-z^2)(1-5z^2)}$\,)
\begin{align*}
	F(u)&=\frac{-1-z^2-W}{2z(u-r_1)}=\frac{1+z^2+W}{2zr_1(1-u/r_1)},\\
	G(u)&=\frac{-1+z^2+W}{2z(u-r_1)}=\frac{1-z^2-W}{2zr_1(1-u/r_1)},\\
	H(u)&=\frac{-1+3z^2+W}{2z(u-r_1)}=\frac{1-3z^2-W}{2zr_1(1-u/r_1)}.
\end{align*}
The total generating function is
\begin{equation*}
	S(u)=F(u)+G(u)+H(u)=\frac{3-3z^2-W}{2zr_1(1-u/r_1)}.
\end{equation*}
The coefficient of $u^jz^n$ in $S(u)$ counts the partial paths of length $n$, ending at level $j$. 
We will write $s_j=[u^j]S(u)$.
Furthermore
\begin{align*}
	f_j=[u^j]	F(u)&=[u^j]\frac{1+z^2+W}{2zr_1(1-u/r_1)},\\
	g_j=[u^j]	G(u)&=[u^j]\frac{1-z^2-W}{2zr_1(1-u/r_1)},\\
	h_j=[u^j]	H(u)&=[u^j]\frac{1-3z^2-W}{2zr_1(1-u/r_1)}.
\end{align*}
At this stage, we are only interested in
\begin{equation*}
	s_j=f_j+g_j+h_j=[u^j]\frac{3-3z^2-W}{2zr_1(1-u/r_1)}=\frac{3-3z^2-W}{2zr_1^{j+1}},
\end{equation*}
which is the generating function of all (partial) paths ending at level $j$. Parity considerations give us that only coefficients $[z^n]s_j$ are non-zero if $n\equiv j\bmod2$.
To make this more transparent, we set 
\begin{equation*}
	P(z)=zr_1=\frac{1+z^2+\sqrt{1-6z^2+5z^4}}{2}, 
\end{equation*}
and then
\begin{equation*}
	s_j=f_j+g_j+h_j=z^j\frac{3-3z^2-W}{2P^{j+1}}.
\end{equation*}
Now we read off coefficients. We do this using residues and contour integration. The path of integration, in both variables $x$ resp.\ $v$ is a
small circle or an equivalent contour.
\begin{align*}
	[z^{2m+j}]s_j&=[z^{2m}]\frac{3-3z^2-W}{2P^{j+1}}=
	[x^m]\frac{3-3x-\sqrt{1-6x+5x^2}}{2\Big(\frac{1+x-\sqrt{1-6x+5x^2}}{2}\Big)^{j+1}}\\
	&=[x^m]\frac{3-3\frac v{1+3v+v^2}-\frac{1-v^2}{1+3v+v^2}}{2\big(\frac{v(v+2)}{1+3v+v^2}\big)^{j+1}}\\
	&=[x^m]\frac{(1+v)(1+2v)}{v^{j+1}(v+2)^{j+1}}(1+3v+v^2)^j\\
	&=\frac1{2\pi i}\oint\frac{dx}{x^{m+1}}\frac{(1+v)(1+2v)}{v^{j+1}(v+2)^{j+1}}(1+3v+v^2)^j\\
	&=\frac1{2\pi i}\oint\frac{dv}{v^{m+1}}\frac{(1+v)(1+2v)(1-v^2)}{v^{j+1}(v+2)^{j+1}}(1+3v+v^2)^{m-1+j}\\
	&=[v^{m+j+1}]\frac{(1+v)^2(1+2v)(1-v)}{(v+2)^{j+1}}(1+3v+v^2)^{m-1+j}.
\end{align*}
Note that
\begin{equation*}(1+v)^2(1+2v)(1-v)=
	-9+27( v+2 ) -29( v+2 ) ^{2}+13( v+2) ^{3}-2( v+2 ) ^{4};
\end{equation*}
consequently
\begin{align*}
	[v^k]&\frac{(1+v)^2(1+2v)(1-v)}{(v+2)^{j+1}}\\
	&=-9\frac1{2^{j+1+k}}\binom{-j-1}{k}
	+27\frac1{2^{j+k}}\binom{-j}{k}
	-29\frac1{2^{j-1+k}}\binom{-j+1}{k}\\&
	+13\frac1{2^{j-2+k}}\binom{-j+2}{k}
	-2\frac1{2^{j-3+k}}\binom{-j+3}{k}=:\lambda_{j;k}.
\end{align*}
With this abbreviation we find
\begin{equation*}
	[v^{m+j+1}]\frac{(1+v)^2(1+2v)(1-v)}{(v+2)^{j+1}}(1+3v+v^2)^{m-1+j}
	=\sum_{k=0}^{m+j+1}\lambda_{j;k}\binom{m-1+j;1,3,1}{m+j+1-k}.
\end{equation*}
This is not extremely pretty but it is \emph{explicit} and as good as it gets.
Here are the first few generating functions:
\begin{align*}
	s_0&=1+z^2+3z^4+10z^6+36z^8+137z^{10}+543z^{12}+\cdots\\*
	s_1&=z+2z^3+6z^5+21z^7+79z^9+311z^{11}+1265z^{13}+\cdots\\
	s_2&=z^2+3z^4+10z^6+37z^8+145z^{10}+589z^{12}+2455z^{14}+\cdots\\
	s_3&=z^3+4z^5+15z^7+59z^9+241z^{11}+1010z^{13}+4314^{15}+\cdots\\
\end{align*}
We could also give such lists for the functions $f_j$, $g_j$, $h_j$, if desired. We summarize the essential findings of this section:
\begin{theorem} The generating function of decorated (partial) Dyck paths, consisting of $n$ steps, ending on level $j$, is given by
	\begin{equation*}
		S(z,u)=\frac{3-3z^2-\sqrt{1-6z^2+5z^4}}{2zr_1(1-u/r_1)},
	\end{equation*}
	with
	\begin{equation*}
		r_1=\frac{1+z^2+\sqrt{1-6z^2+5z^4}}{2z}.
	\end{equation*}
	Furthermore
	\begin{equation*}
		[u^j]S(z,u)=\frac{3-3z^2-\sqrt{1-6z^2+5z^4}}{2zr_1^{j+1}}.
	\end{equation*}
\end{theorem}

\subsection*{Open ended paths}

If we do not specify the end of the paths, in other words we sum over all $j\ge0$, then at the level of generating functions
this is very easy, since we only have to set $u:=1$.
We find
\begin{align*}
	S(1)&=-\frac{(z+1)(z^2+3z-2)+(z+2)\sqrt{1-6z^2+5z^4}}{2z(z^2+2z-1)}\\
	&=1+z+2z^2+3z^3+7z^4+11z^5+26z^6+43z^7+102z^8+175z^9+416z^{10}+\cdots.
\end{align*}

\subsection*{Counting red edges}

We can use an extra variable, $w$, to count additionally the red edges that occur in a path. We use the same
letters for generating functions. Eventually, the coefficient $[z^nu^jw^k]S$ is the number of (partial) paths consisting of $n$ steps, leading
to level $j$, and having passed $k$ red edges. The endpoint of the original skew path has then coordinates $(n-2k,j)$. The computations are very similar, and we only
sketch the key steps.

\begin{equation*}
	f_0=1,\quad f_{i+1}=zf_i+zg_i,\quad i\ge0,
\end{equation*}
\begin{equation*}
	g_i=zf_{i+1}+zg_{i+1}+zh_{i+1},\quad i\ge0,
\end{equation*}
\begin{equation*}
	h_i=wzh_{i+1}+wzg_{i+1},\quad i\ge0;
\end{equation*}
\begin{align*}
	\frac1u(F(u)-1)&=zF(u)+zG(u),\\*
	G(u)&=\frac zu(F(u)-1)+\frac zu(G(u)-G(0))+\frac zu(H(u)-H(0)),\\*
	H(u)&=	\frac {wz}u(G(u)-G(0))+\frac {wz}u(H(u)-G(0));
\end{align*}
\begin{align*}
	F(u)&=\frac{z^2uG(0)+z^2uH(0)+z^2u-u-wz^3+z+wz}{-w{z}^{3}-u+z+wz+z{u}^{2}-w{z}^{2}u},\\
	G(u)&=\frac{z(H(0)-uzH(0)+wz^2+G(0)-zuG(0)-zu)}{-w{z}^{3}-u+z+wz+z{u}^{2}-w{z}^{2}u},\\
	H(u)&=\frac{wz(-uzH(0)-z^2-zuG(0)+G(0)-z^2H(0)+H(0)-z^2G(0))}{-w{z}^{3}-u+z+wz+z{u}^{2}-w{z}^{2}u}.
\end{align*}
The denominator factors as $z(u-r_1)(u-r_2)$, with
\begin{align*} 
	r_1&=\frac{1+wz^2+\sqrt{1-(4+2w)z^2+(4w+w^2)z^4}}{2z},\\*
	r_2&=\frac{1+wz^2-\sqrt{1-(4+2w)z^2+(4w+w^2)z^4}}{2z}.
\end{align*}
Note the factorization $1-(4+2w)z^2+(4w+w^2)z^4=(1-z^2w)(1-(4+w)z^2)$. 
Since the factor $u-r_2$ in the denominator is ``bad,'' it must also cancel in the numerators. From this
we eventually find, with the abbreviation
$W=\sqrt{1-(4+2w)z^2+(4w+w^2)z^4}\,$)
\begin{align*}
	F(u)&=\frac{-1-wz^2-W}{2z(u-r_1)},\\
	G(u)&=\frac{-1+wz^2+W}{2z(u-r_1)},\\
	H(u)&=\frac{-1+(2+w)z^2+W}{2z(u-r_1)}.
\end{align*}
The total generating function is
\begin{equation*}
	S(u)=F(u)+G(u)+H(u)=\frac{-2-w+z^2(w+w^2)+ wW}{2z(u-r_1)}.
\end{equation*}
The special case $u=0$ (return to the $x$-axis) is to be noted:
\begin{equation*}
	S(0)=\frac{-2-w+z^2(w+w^2)+ wW}{-2zr_1}=\frac{1-wz^2-W}{2z^2}.
\end{equation*}
Since there are only even powers of $z$ in this function, we replace $x=z^2$ and get
\begin{align*}
	S(0)&=\frac{1-wx-\sqrt{1-(4+2w)x+(4w+w^2)x^2}}{2x}\\
	&=1+x+(w+2)x^2+(w^2+4w+5)x^3+(w^3+6w^2+15w+14)x^4+\cdots.
\end{align*}
Compare the factor $(w^2+4w+5)$ with the earlier drawing of the 10 paths.

There is again a substitution that allows for better results:
\begin{equation*}
	z=\frac{v}{1+(2+w)v+v^2}, \quad\text{then}\quad S(0)=1+v.
\end{equation*}
Reading off coefficients can now be done using modified trinomial coefficients:
\begin{equation*}
	\binom{n;1,2+w,1}{k}=[t^k]\bigl(1+(2+w)t+t^2\bigr)^n.
\end{equation*}
Again, we use contour integration to extract coefficients:
\begin{align*}
	[x^n](1+v)&=\frac1{2\pi i}\oint \frac{dx}{x^{n+1}}(1+v)\\
	&=\frac1{2\pi i}\oint \frac{dx}{v^{n+1}}\frac{1-v^2}{(1+(2+w)v+v^2)^2}(1+(2+w)v+v^2)^{n+1}(1+v)\\
	&=[v^n](1-v)(1+v)^2(1+(2+w)v+v^2)^{n-1}\\
	&=\binom{n-1;1,2+w,1}{n}+\binom{n-1;1,2+w,1}{n-1}\\*
	&\qquad-\binom{n-1;1,2+w,1}{n-2}-\binom{n-1;1,2+w,1}{n-3}.
\end{align*}

Now we want to count the average number of red edges. For that, we differentiate $S(0)$ w.r.t.\ $w$, followed by $w:=1$.
This leads to
\begin{equation*}
	\frac{-1+6x-5x^2+(1+3x)\sqrt{1-6x+5x^2}}{2(1-x)(1-5x)}.
\end{equation*}

A simple application of singularity analysis leads to
\begin{equation*}
	\frac{\frac1{2\sqrt5}[x^n]\frac1{\sqrt{1-5x}}}{-\sqrt5[x^n]\sqrt{1-5x}}\sim \frac {n}{5}.
\end{equation*}
So, a random path consisting of $2n$ steps has about $n/5$ red steps, on average. 

For readers who are not familiar with singularity analysis of generating functions \cite{FlOd90, FS}, we just
mention that one determines the local expansion around the dominating singularity, which is at $z=\frac15$ in our instance.
In the denominator, we just have the total number of skew Dyck paths, according to the sequence A002212 in \cite{OEIS}.

In the example of Figure~2, the exact average is $6/10$, which curiously is exactly the same as $3/5$.

We finish the discussion by considering fixed powers of $w$ in $S(0)$, counting skew Dyck paths consisting of zero, one, two, three, \dots red edges. We find
\begin{align*}
	[w^0]S(0)&=\frac{1-\sqrt{1-4x}}{2x},\\
	[w^1]S(0)&=\frac{1-2x-\sqrt{1-4x}}{2\sqrt{1-4x}},\\
	[w^2]S(0)&=\frac{x^3}{(1-4x)^{3/2}},\\
	[w^3]S(0)&=\frac{x^4(1-2x)}{(1-4x)^{5/2}},\\
	[w^4]S(0)&=\frac{x^5(1-4x+5x^2)}{(1-4x)^{7/2}}, \quad\&\text{c}.
\end{align*}
The generating function $[w^0]S(0)$ is of course the generating function of Catalan numbers, since no red edges just means: ordinary Dyck paths.
We can also conclude that the asymptotic behaviour is of the form $n^{k-3/2}4^n$, where the polynomial contribution gets higher, but the exponential growth
stays the same: $4^n$. This is compared to the scenario of an \emph{arbitrary} number of red edges, when we get an exponential growth of the form $5^n$.

\subsection*{Dual skew Dyck paths}

The mirrored version of skew Dyck paths with two types of up-steps, $(1,1)$ and $(-1,1)$ are also cited among the objects in A002212 in \cite{OEIS}.
We call them dual skew paths and drop the `dual' when it isn't necessary. When the paths come back to the $x$-axis, no new enumeration is necessary, but this
is no longer true for paths ending at level $j$.

Here is a list of the 10 skew paths consisting of 6 steps:

\begin{figure}[h]
	\begin{equation*}
		\begin{tikzpicture}[scale=0.3]
			\draw [thick](0,0)--(3,3)--(6,0);
		\end{tikzpicture}
		\quad
		\begin{tikzpicture}[scale=0.3]
			\draw [thick](0,0)--(-1,1)--(1,3)--(4,0);
		\end{tikzpicture}
		\quad
		\begin{tikzpicture}[scale=0.3]
			\draw [thick](0,0)--(1,1)--(0,2)--(1,3)--(4,0);
		\end{tikzpicture}
		\quad
		\begin{tikzpicture}[scale=0.3]
			\draw [thick](0,0)--(-2,2)--(-1,3)--(2,0);
		\end{tikzpicture}
		\quad
		\begin{tikzpicture}[scale=0.3]
			\draw [thick](0,0)--(2,2)--(4,0)--(5,1)--(6,0);
		\end{tikzpicture}
	\end{equation*}
	\begin{equation*}
		\begin{tikzpicture}[scale=0.3]
			\draw [thick](0,0)--(1,1)--(2,0)--(4,2)--(6,0);
		\end{tikzpicture}
		\quad
		\begin{tikzpicture}[scale=0.3]
			\draw [thick](0,0)--(-1,1)--(0,2)--(2,0)--(3,1)--(4,0);
		\end{tikzpicture}
		\quad
		\begin{tikzpicture}[scale=0.3]
			\draw [thick](0,0)--(1,1)--(2,2)--(3,1)--(4,2)--(6,0);
		\end{tikzpicture}
		\quad
		\begin{tikzpicture}[scale=0.3]
			\draw [thick](0,0)--(-1,1)--(0,2)--(1,1)--(2,2)--(4,0);
		\end{tikzpicture}
		\quad
		\begin{tikzpicture}[scale=0.3]
			\draw [thick](0,0)--(1,1)--(2,0)--(3,1)--(4,0)--(5,1)--(6,0);
		\end{tikzpicture}
	\end{equation*}
	\caption{All 10 dual skew Dyck paths of length 6 (consisting of 6 steps).}
\end{figure}

We prefer to work with the equivalent model (resembling more traditional Dyck paths) where
we replace each step $(-1,-1)$ by $(1,-1)$ but label it blue. Here is the list of the 10 paths again (Figure 2):

\begin{figure}[h]
	\begin{equation*}
		\begin{tikzpicture}[scale=0.3]
			\draw [thick](0,0)--(3,3)--(6,0);
		\end{tikzpicture}
		\quad
		\begin{tikzpicture}[scale=0.3]
			\draw [thick,cyan](0,0)--(1,1);
			\draw [thick](1,1)--(3,3)--(6,0);
		\end{tikzpicture}
		\quad
		\begin{tikzpicture}[scale=0.3]
			\draw [thick](0,0)--(1,1);
			\draw [thick,cyan](1,1)--(2,2);
			\draw [thick](2,2)--(3,3)--(6,0);
		\end{tikzpicture}
		\quad
		\begin{tikzpicture}[scale=0.3]
			\draw [thick,cyan](0,0)--(2,2);
			\draw [thick](2,2)--(3,3)--(6,0);
		\end{tikzpicture}
		\quad
		\begin{tikzpicture}[scale=0.3]
			\draw [thick](0,0)--(2,2)--(4,0)--(5,1)--(6,0);
		\end{tikzpicture}
	\end{equation*}
	\begin{equation*}
		\begin{tikzpicture}[scale=0.3]
			\draw [thick](0,0)--(1,1)--(2,0)--(4,2)--(6,0);
		\end{tikzpicture}
		\quad
		\begin{tikzpicture}[scale=0.3]
			\draw [thick,cyan](0,0)--(1,1);
			\draw [thick](1,1)--(2,2)--(4,0)--(5,1)--(6,0);
		\end{tikzpicture}
		\quad
		\begin{tikzpicture}[scale=0.3]
			\draw [thick](0,0)--(1,1)--(2,2)--(3,1)--(4,2)--(6,0);
		\end{tikzpicture}
		\quad
		\begin{tikzpicture}[scale=0.3]
			\draw [thick,cyan](0,0)--(1,1);
			\draw [thick](1,1)--(2,2)--(3,1)--(4,2)--(6,0);
		\end{tikzpicture}
		\quad
		\begin{tikzpicture}[scale=0.3]
			\draw [thick](0,0)--(1,1)--(2,0)--(3,1)--(4,0)--(5,1)--(6,0);
		\end{tikzpicture}
	\end{equation*}
	\caption{All 10 dual skew Dyck paths of length 6 (consisting of 6 steps).}
\end{figure}

The rules to generate such decorated Dyck paths are: each edge $(1,-1)$ may be  black or blue, but
\begin{tikzpicture}[scale=0.3]\draw [thick](0,1)--(1,0); \draw [cyan,thick] (1,0)--(2,1);\end{tikzpicture}
and
\begin{tikzpicture}[scale=0.3] \draw [cyan,thick] (0,0)--(1,1);\draw [thick](1,1)--(2,0);\end{tikzpicture}
are forbidden.

Our interest is in particular in \emph{partial} decorated Dyck paths, ending at level $j$, for fixed $j\ge0$;
the instance $j=0$ is the classical case. 

As before, two variables, $z$ and $u$, are used, where $z$ marks the length
of the path and $j$ marks the end-level. We briefly mention that one can, using a third variable $w$, also
count the number of blue edges.

The substitution
\begin{equation*}
	x=\frac{v}{1+3v+v^2},
\end{equation*}
  is again the key to the success.

\subsection*{Generating functions and the kernel method}

We catch the essence of a decorated (dual skew) Dyck path using a state-diagram:

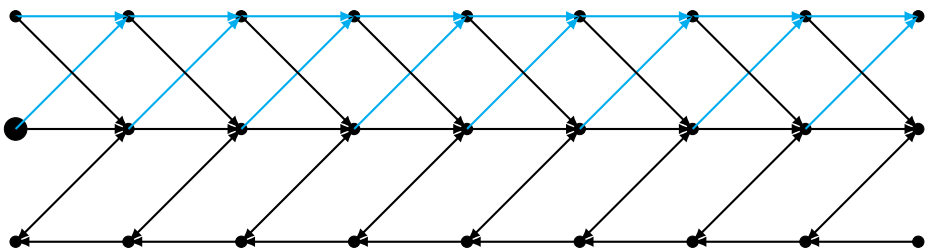
\begin{figure}[h]
	
	%\draw[thick,red, -latex]   (0.5,1) [out=-90,in=145] to (1,0);
	\begin{center}
		\begin{tikzpicture}[scale=1.5]
			\draw (0,0) circle (0.1cm);
			\fill (0,0) circle (0.1cm);
			
			\foreach \x in {0,1,2,3,4,5,6,7,8}
			{
				\draw (\x,0) circle (0.05cm);
				\fill (\x,0) circle (0.05cm);
			}
			
			\foreach \x in {0,1,2,3,4,5,6,7,8}
			{
				\draw (\x,-1) circle (0.05cm);
				\fill (\x,-1) circle (0.05cm);
			}
			
			\foreach \x in {0,1,2,3,4,5,6,7,8}
			{
				\draw (\x,1) circle (0.05cm);
				\fill (\x,1) circle (0.05cm);
			}
			
			\foreach \x in {0,1,2,3,4,5,6,7}
			{
				\draw[ thick,-latex] (\x,0) -- (\x+1,0);
				
			}
			\foreach \x in {0,1,2,3,4,5,6,7}
			{
				\draw[ thick,-latex,cyan] (\x,1) -- (\x+1,1);
				
			}
			
			\foreach \x in {0,1,2,3,4,5,6,7}
			{
				\draw[ thick,-latex,cyan] (\x,0) to (\x+1,1);
				
			}
			
			\foreach \x in {0,1,2,3,4,5,6,7}
			{
				\draw[ thick,-latex] (\x,1)  to (\x+1,0);
				
			}
			
			\foreach \x in {0,1,2,3,4,5,6,7}
			{
				\draw[ thick,latex-latex] (\x+1,0)  to (\x,-1);
				
			}
			
			\foreach \x in {0,1,2,3,4,5,6,7}
			{
				\draw[ thick,-latex] (\x+1,-1)  to (\x,-1);
				
			}

		\end{tikzpicture}
	\end{center}
	\caption{Three layers of states according to the type of steps leading to them (down, up-black, up-blue).}
\end{figure}
It has three types of states, with $j$ ranging from 0 to infinity; in the drawing, only $j=0..8$ is shown. The first
layer of states refers to an up-step leading to a state, the second layer refers to a black down-step leading to a state
and the third layer refers to a blue down-step leading to a state. We will work out generating functions describing all paths
leading to a particular state. We will use the notations $c_j,a_j,b_j$ for the three respective layers, from top to bottom.
Note that the syntactic rules of forbidden patterns
\begin{tikzpicture}[scale=0.3]\draw [thick,cyan](0,0)--(1,1); \draw [thick] (1,1)--(2,0);\end{tikzpicture}
and
\begin{tikzpicture}[scale=0.3] \draw [thick] (0,1)--(1,0);\draw [thick,cyan](1,0)--(2,1);\end{tikzpicture}
can be clearly seen from the picture. The functions depend on the variable $z$ (marking the number of steps), but mostly
we just write $a_j$ instead of $a_j(z)$, etc.

The following recursions can be read off immediately from the diagram:
\begin{gather*}
	a_0=1,\quad a_{i+1}=za_i+zb_i+zc_i,\quad i\ge0,\\
	b_i=za_{i+1}+zb_{i+1},\quad i\ge0,\\
	c_{i+1}=za_{i}+zc_{i},\quad i\ge0.
\end{gather*}
And now it is time to introduce the promised \emph{bivariate} generating functions: 
\begin{equation*}
	A(z,u)=\sum_{i\ge0}a_i(z)u^i,\quad
	B(z,u)=\sum_{i\ge0}b_i(z)u^i,\quad
	C(z,u)=\sum_{i\ge0}c_i(z)u^i.
\end{equation*}
Again, often we just write $A(u)$ instead of $A(z,u)$ and treat $z$ as a `silent' variable. Summing the recursions leads to
\begin{align*}
	\sum_{i\ge0}u^ia_i &=1+u\sum_{i\ge0}u^i(za_i+zb_i+zc_i)\\
	&=1+uzA(u)+uzB(u)+uzC(u),\\
	\sum_{i\ge0}u^ib_i &= \sum_{i\ge0}u^i(za_{i+1}+zb_{i+1})\\
	&=\frac zu\sum_{i\ge1}u^ia_i+\frac zu\sum_{i\ge1}u^ib_i,\\
	\sum_{i\ge1}u^ic_i &=uz\sum_{i\ge0}u^ia_i+uz\sum_{i\ge0}u^ic_i.
\end{align*}
This can be rewritten as
\begin{align*}
	A(u)&=1+uzA(u)+uzB(u)+uzC(u),\\
	B(u)&=\frac zu(A(u)-a_0)+\frac zu(B(u)-b_0),\\
	C(u)&=c_0+uzA(u)+uzC(u).
\end{align*}
Note that $a_0=1$, $c_0=0$.
Simplification leads to
\begin{equation*}
	C(u)=\frac{uzA(u)}{1-uz}
\end{equation*}
and
\begin{equation*}
	B(u)=\frac{z(A(u)-1-B(0))}{u-z},
\end{equation*}
leaving us with just one equation
\begin{equation*}
	A(u)={\frac { \left( z-u+u{z}^{2}+u{z}^{2}B(0) \right)  \left( uz-1
			\right) }{{u}^{2}{z}^{3}+u{z}^{2}-2{u}^{2}z-z+u}}.
\end{equation*}
This is again a typical application of the kernel method. 
\begin{equation*}
	{u}^{2}{z}^{3}+u{z}^{2}-2{u}^{2}z-z+u=z(z^2-2)(u-s_1)(u-s_2)
\end{equation*}
The denominator factors as $2z(z^2-2)(u-s_1)(u-s_2)$, with
\begin{equation*}
	s_1=\frac{1+z^2+\sqrt{1-6z^2+5z^4}}{2z(2-z^2)},\quad s_2=\frac{1+z^2-\sqrt{1-6z^2+5z^4}}{2z(2-z^2)}.
\end{equation*}
Note that $s_1s_2=\frac{1}{2-z^2}$. 
Since the factor $u-s_2$ in the denominator is ``bad,'' it must also cancel in the numerators. From this
we conclude (again with the abbreviation $W=\sqrt{1-6z^2+5z^4}\,$)
\begin{equation*}
	B(0) =  \frac{zs_2}{1-2zs_2},
\end{equation*}
and further
\begin{equation*}
	A(u)
	=\frac{(1-uz)(1+z^2+W)}{2z(z^2-2)(u-s_1)},
\end{equation*}
\begin{equation*}
	B(u)=\frac{1-2z^2-W}{z(2-z^2)(u-s_1)},
\end{equation*}
\begin{equation*}
	C(u)=\frac{1+z^2+W}{2(z^2-2)}\frac{u}{u-s_1},
\end{equation*}
and for the function of main interest
\begin{equation*}
	G(u)=A(u)+B(u)+C(u)=\frac{3z^2-3+W}{2z(2-z^2)(u-s_1)}.
\end{equation*}
Note that
\begin{align*}
	\frac1{s_1}&=\frac{1+z^2-\sqrt{1-6z^2+5z^4}}{2z}=zS,\\
	\frac1{s_2}&=\frac{1+z^2+\sqrt{1-6z^2+5z^4}}{2z}.
\end{align*}
Then
\begin{align*}
	[u^j]G(u)&=[u^j]\frac{3z^2-3+W}{2z(z^2-2)s_1(1-u/s_1)}\\
	&=\frac{3z^2-3+W}{2z(z^2-2)s_1^{j+1}}
	=\frac{3z^2-3+W}{2(z^2-2)}z^{j}S^{j+1}.
\end{align*}
So $[u^j]G(u)$ contains only powers of the form $z^{j+2N}$. Now we continue
\begin{align*}
	[z^{j+2N}u^j]G(u)&
	=[z^{2N}]\frac{3z^2-3+W}{2(z^2-2)}S^{j+1}
	\\&=[x^{N}]\frac{3x-3+\sqrt{1-6x+5x^2}}{2(x-2)}\bigg(\frac{1+x-\sqrt{1-6x+5x^2}}
	{2x}\bigg)^{j+1}\\
	&=[x^{N}](v+1)(v+2)^{j}
\end{align*}
which is the generating function of all (partial) paths ending at level $j$.

Now we read off coefficients. We do this using residues and contour integration. The path of integration, in both variables $x$ resp.\ $v$ is a
small circle or an equivalent contour;
\begin{align*}
	[z^{j+2N}u^j]G(u)&=[x^{N}](v+1)(v+2)^{j}\\
	&=\frac1{2\pi i}\oint \frac{dx}{x^{N+1}}(v+1)(v+2)^{j}\\
	&=\frac1{2\pi i}\oint \frac{dv}{v^{N+1}}(1+3v+v^2)^{N+1}\frac{(1-v^2)}{(1+3v+v^2)^2}(v+1)(v+2)^{j}\\
	&=[v^{N}](1+3v+v^2)^{N-1}(1-v)(1+v)^2(v+2)^{j}.
\end{align*}
Note that
\begin{equation*}(1-v)(1+v)^2=
	3-7( v+2 ) +5( v+2 ) ^{2}- ( v+2) ^{3};
\end{equation*}
consequently
\begin{align*}
	[z^{j+2N}u^j]G(u)&=[v^{N}](1+3v+v^2)^{N-1}\Big[3-7( v+2 ) +5( v+2 ) ^{2}- ( v+2) ^{3}
	\Big](v+2)^{j}.
\end{align*}
We abbreviate:
\begin{align*}
	\mu_{j;k}&=[v^{k}]\Big[3(v+2)^{j}-7(v+2)^{j+1} +5(v+2)^{j+2}- (v+2)^{j+3}\Big]\\
	&=3\binom{j}{k}2^{j-k}-7\binom{j+1}{k}2^{j+1-k}+5\binom{j+2}{k}2^{j+2-k}-\binom{j+3}{k}2^{j+3-k}.
\end{align*}
With this notation we get
\begin{equation*}
	[z^{j+2N}u^j]G(u)
	=\sum_{0\le k\le N-1}\mu_{j;k}\binom{N-1;1,3,1}{N-k}.
\end{equation*}
Here are the first few generating functions:
\begin{align*}
	G_0&=1+{z}^{2}+3{z}^{4}+10{z}^{6}+36{z}^{8}+137{z}^{10}+543{z}^{
		12}+2219{z}^{14}
	+\cdots\\*
	G_1&=2z+3{z}^{3}+10{z}^{5}+36{z}^{7}+137{z}^{9}+543{z}^{11}+
	2219{z}^{13}+9285{z}^{15}
	+\cdots\\
	G_2&=4{z}^{2}+8{z}^{4}+29{z}^{6}+111{z}^{8}+442{z}^{10}+1813{z
	}^{12}+7609{z}^{14}+32521{z}^{16}
	+\cdots\\
	G_3&=8{z}^{3}+20{z}^{5}+78{z}^{7}+315{z}^{9}+1306{z}^{11}+5527
	{z}^{13}+23779{z}^{15}+103699{z}^{17}
	+\cdots\\
\end{align*}
We could also give such lists for the functions $a_j$, $b_j$, $c_j$, if desired. We summarize the essential findings of the rest of this section:
\begin{theorem} The generating function of decorated (partial) dual skew Dyck paths, consisting of $n$ steps, ending on level $j$, is given by
	\begin{equation*}
		G(z,u)=\frac{3z^2-3+\sqrt{1-6z^2+5z^4}}{2z(2-z^2)(u-s_1)},
	\end{equation*}
	with
	\begin{equation*}
		s_1=\frac{2z}{1+z^2-\sqrt{1-6z^2+5z^4}}.
	\end{equation*}
	Furthermore
	\begin{equation*}
		[u^j]G(z,u)=\frac{3z^2-3+\sqrt{1-6z^2+5z^4}}{2(z^2-2)}z^jS^{j+1},
	\end{equation*}
	with
	\begin{equation*}
		S=\frac{1+z^2-\sqrt{1-6z^2+5z^4}}{2z^2}.
	\end{equation*}
\end{theorem}

\subsection*{Open ended paths}

If we do not specify the end of the paths, in other words we sum over all $j\ge0$, then at the level of generating functions
this is very easy, since we only have to set $u:=1$.
We find
\begin{align*}
	G(1)&=\frac{(1+z)(1-3z)}{2z(z^2+2z-1)-\sqrt{1-6z^2+5z^4}}\\
	&=1+2z+5{z}^{2}+11{z}^{3}+27{z}^{4}+62{z}^{5}+151{z}^{6}+
	354{z}^{7}+859{z}^{8}+2036{z}^{9}+\cdots.
\end{align*}

\subsection*{Counting blue edges}

We can use an extra variable, $w$, to count additionally the blue edges that occur in a path. We use the same
letters for generating functions. Eventually, the coefficient $[z^nu^jw^k]S$ is the number of (partial) paths consisting of $n$ steps, leading
to level $j$, and having passed $k$ blue edges. The endpoint of the original skew path has then coordinates $(n-2k,j)$. The computations are very similar, and we only
sketch the key steps.

\begin{gather*}
	a_0=1,\quad a_{i+1}=za_i+zb_i+zc_i,\quad i\ge0,\\
	b_i=za_{i+1}+zb_{i+1},\quad i\ge0,\\
	c_{i+1}=wza_{i}+wzc_{i},\quad i\ge0.
\end{gather*}
This leads to
\begin{align*}
	A(u)&=1+uzA(u)+uzB(u)+uzC(u),\\
	B(u)&=\frac zu(A(u)-a_0)+\frac zu(B(u)-b_0),\\
	C(u)&=c_0+wuzA(u)+wuzC(u).
\end{align*}
Solving,
\begin{equation*}
	S(u)=A(u)+B(u)+C(u)={\frac {u-wu{z}^{2}-zA(0)-zB(0)+uw{z}^{2}A(0)+uw{z}^{2}B(0)}{{u}^{2}{z}^{3}w+u-w{u}^{2}z-{u}^{2}z-z+wu{z}^{2}}}.
\end{equation*}
The denominator factors as $-z(1+w-z^2w)(u-s_1)(u-s_2)$, with
\begin{align*} 
	s_1&={\frac {1+{z}^{2}w+\sqrt {1-2\,{z}^{2}w+{z}^{4}{w}^{2}-4\,{z}^{2
				}+4{z}^{4}w}}{2z \left( 1+w-{z}^{2}w \right) }},\\*
	s_2&= {\frac {1+{z}^{2}w-\sqrt {1-2\,{z}^{2}w+{z}^{4}{w}^{2}-4\,{z}^{2
				}+4{z}^{4}w}}{2z \left( 1+w-{z}^{2}w \right) }}	.
\end{align*}
Note the factorization $1-(4+2w)z^2+(4w+w^2)z^4=(1-z^2w)(1-(4+w)z^2)$. 
Since the factor $u-r_2$ in the denominator is ``bad,'' it must also cancel in the numerators. From this
we eventually find, with the abbreviation
$W=\sqrt{1-(4+2w)z^2+(4w+w^2)z^4}\,$)
\begin{equation*}
	G(0)={\frac {1-{z}^{2}w-W }{2{z}^{2}}},
\end{equation*}
and further
\begin{equation*}
	G(u)=\frac {w-{z}^{2}{w}^{2}-wW+2-2{z}^{2}w}
	{2z \left( -w	-1+{z}^{2}w \right) (u-s_1)}.
\end{equation*}
The special case $u=0$ (return to the $x$-axis) is to be noted:
\begin{equation*}
	G(0)=1+{z}^{2}+ \left( w+2 \right) {z}^{4}+ \left( {w}^{2}+4w+5 \right) 
	{z}^{6}+ \left( w+2 \right)  \left( {w}^{2}+4w+7 \right) {z}^{8}+\cdots.
\end{equation*}
Compare the factor $(w^2+4w+5)$ with the earlier drawing of the 10 paths.
There is again a substitution that allows for better results:
\begin{equation*}
	z=\frac{v}{1+(2+w)v+v^2}, \quad\text{then}\quad G(0)=1+v.
\end{equation*}
Since $S(u)=G(u)$ with $S(u)$ from the first part of the paper, as it means the same objects, read from left to right resp.\ from right to left, no 
new analysis is required.

\section{Counting ternary trees according to the number of middle edges}

	The recent preprint \cite{Burstein} triggered my interest in the sequence A120986 in \cite{OEIS}. The double-indexed sequence
enumerates ternary trees according to the number of edges and the number of middle edges. We consider here $T(n,k)$, the
number of ternary trees with $n$ nodes and $k$ middle edges. The difference is marginal, but we want to compare/relate our analysis with \cite{christmas}, and there it is also the number of nodes that is considered. Let
$G=G(x,u)=\sum_{n,k\ge0}T(n,k)x^nu^k$. Then it is easy to see (decomposition at the root) that
\begin{equation*}
	G=1+xG^2(1-u+uG).
\end{equation*}
The substitution
\begin{equation*}
	x=\frac{t(1-t)^2}{(1-t+ut)}
\end{equation*}
makes the cubic equation manageable and also allows, as in \cite{christmas}, to introduce a (refined) version of the $(3/2)$-ary trees.

Here is a small table of these numbers and a ternary tree:
\begin{center}
	
	\begin{tabular}{c|ccccccccccc}
		
		$n\backslash k$& 0&1&2&3&4&5\\
		\hline
		0&1&&&&&&&\\
		1&1&&&&&&&\\
		2&2&1&&&&&&\\
		3&5&6&1&&&&&\\
		4&14&28&12&1&&&&\\
		5&42&120&90&20&1&&&\\
		6&132&495&550&220&30&1&&\\
	\end{tabular}
	
\end{center}
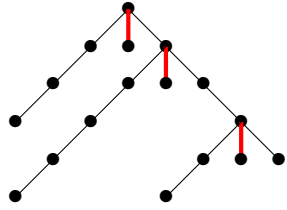
\begin{figure}[h]
	\begin{tikzpicture}[scale=0.5]
		
		\draw[](0,0)to(4,-4);
		
		\node at (0,0) {$\bullet$};\node at (1,-1) {$\bullet$};\node at (2,-2) {$\bullet$};\node at (3,-3) {$\bullet$};
		\node at (4,-4) {$\bullet$};
		
		\draw(0,0)to(-3,-3);\node at (-1,-1) {$\bullet$};\node at (-2,-2) {$\bullet$};\node at (-3,-3) {$\bullet$};

		\draw[ultra thick,red](3,-3)to(3,-4);
		\draw[ultra thick,red](1,-1)to(1,-2);\draw(1,-1)to(-3,-5);
		\draw[ultra thick,red](0,0)to(0,-1);\node at (0,-1) {$\bullet$};
		
		\node at (3,-4) {$\bullet$};
		\node at (1,-2) {$\bullet$};\node at (0,-2) {$\bullet$};\node at (-1,-3) {$\bullet$};
		\node at (-2,-4) {$\bullet$};\node at (-3,-5) {$\bullet$};
		
		\draw(3,-3)to(1,-5);\node at (1,-5) {$\bullet$};\node at (2,-4) {$\bullet$};

	\end{tikzpicture}
	
	\caption{Ternary tree with 17 nodes and 3 middle edges}
	
\end{figure}

\subsection*{Analysis of the cubic equation}

The cubic equation has the following solutions:
\begin{align*}
	r_1&=\frac{1}{1-t},\\
	r_2&=\frac{-t+t^2-t^2u-\sqrt{t(1-t+ut)(4u+t-4ut-t^2+t^2u)}}{2ut(1-t)},\\
	r_3&=\frac{-t+t^2-t^2u+\sqrt{t(1-t+ut)(4u+t-4ut-t^2+t^2u)}}{2ut(1-t)}.
\end{align*}
Note that
\begin{equation*}
	r_2r_3=-\frac{1-t+ut}{ut(1-t)}.
\end{equation*}
The root with the combinatorial significance is $r_1$. But it is the explicit form of the two other roots that makes everything here interesting and challenging. 

We extract coefficients of $r_1$ using contour integration, which is closely related to the Lagrange inversion formula. The path of integration is a small circle in the $x$-plane which is then transformed into a small circle in the $t$-plane.
\begin{align*}
	[x^n]r_1&=\frac1{2\pi i}\oint \frac{dx}{x^{n+1}}\frac1{1-t}\\
	&=\frac1{2\pi i}\oint \frac{dt(1-t)(1-3t+2t^2-2t^2u)}{(1-t+tu)^2}\frac{(1-t+tu)^{n+1}} {t^{n+1}(1-t)^{2n+2}}\frac1{1-t}\\
	&=[t^n](1-3t+2t^2-2t^2u)\frac{(1-t+tu)^{n-1}} {(1-t)^{2n+2}}.
\end{align*}
Furthermore
\begin{align*}
	[x^nu^k]r_1
	&=[t^n][u^k](1-3t+2t^2-2t^2u)\frac{(1-t+tu)^{n-1}} {(1-t)^{2n+2}}\\
	&=[t^n]\binom{n-1}{k}\frac{t^k(1-2t)}{(1-t)^{n+k+2}}	-2[t^n]\binom{n-1}{k-1}\frac{t^{k+1}}{(1-t)^{n+k+2}}\\
	&=\binom{n-1}{k}[t^{n-k}]\frac{(1-2t)}{(1-t)^{n+k+2}}	-2\binom{n-1}{k-1}[t^{n-k-1}]\frac{1}{(1-t)^{n+k+2}}\\
	&=\binom{n-1}{k}\binom{2n+1}{n-k}-2\binom{n-1}{k}\binom{2n}{n-k-1}
	-2\binom{n-1}{k-1}\binom{2n}{n-k-1}\\
	&=\frac1n\binom nk\binom{2n}{n-1-k}.
\end{align*}
For $u=1$, which means that the middle edges are not especially counted, we get
\begin{equation*}
	\sum_k\frac1n\binom nk\binom{2n}{n-1-k}=\frac1n\binom{3n}{n-1},
\end{equation*}
the number of ternary trees with $n$ nodes.

\subsection*{Factorizing the solution of the cubic equation}

For $u=1$, Knuth \cite{christmas} was able to factor the generating function $r_1$ into two factors, for which he coined the catchy name $(3/2)$-ary trees. For this factorization, see also \cite{naimi-paper,BM-P}. The goal in this section is to perform this factorization in the context of counting middle edges, i. e., for the generating function with the additional variable $u$. In Knuth's instance, the generating function was expressible as a generalized binomial series (in the sense of Lambert \cite{GKP}), but that does not seem to be an option here.

Note that
\begin{equation*}
	\frac1{r_2}=\frac t2-\frac{\sqrt t\sqrt{t(1-t)+u(2-t)^2}}{2\sqrt{1-t+tu}}
\end{equation*}
and
\begin{equation*}
	\frac1{r_3}=\frac t2+\frac{\sqrt t\sqrt{t(1-t)+u(2-t)^2}}{2\sqrt{1-t+tu}}.
\end{equation*}
From the cubic equation we deduce that
\begin{equation*}
	r_1=-\frac1{uxr_2r_3},
\end{equation*}
which is the desired factorization. The factor $ux$ will be fairly split as $\sqrt{ux}\cdot\sqrt{ux}$, whereas the minus sign goes to the
factor $1/r_2$. In the following we work out how this factorization can be obtained. To say it again, it is not as appealing as in the original case.

Let us write
\begin{equation*}
	t=x\Phi(t), \quad\text{with}\quad  \Phi(t)=\frac{1-t+tu}{(1-t)^2},
\end{equation*}
so that we can use the Lagrange inversion formula to get
\begin{equation*}
	[x^n]t^{\ell}=\frac{\ell}{n}[t^{n-\ell}]\Phi(t)^n
\end{equation*}
and
\begin{align*}
	[x^nu^k]t^{\ell}&=\frac{\ell}{n}[t^{n-\ell}][u^k]\frac{(1-t+tu)^n}{(1-t)^{2n}}\\
	&=\frac{\ell}{n}[t^{n-\ell-k}]\binom{n}{k}\frac1{(1-t)^{n+k}}
	=\frac{\ell}{n}\binom{n}{k}\binom{2n-\ell-1}{n-\ell-k}.
\end{align*}
In particular,
\begin{equation*}
	t=\sum_{n\ge1}x^n\sum_{0\le k\le n}\frac{1}{n}\binom{n}{k}\binom{2n-2}{n-1-k}u^k;
\end{equation*}
this series expansion may be used in the following developments whenever needed.

To proceed further, we set $u=1+U$ and $\tau=t/u$:
\begin{gather*}
	\frac1{r_2}=\frac t2-\frac{\sqrt{x}}{2(1-t)}\sqrt{4-3t+U(2-t)^2},\\
	\frac1{r_3}=\frac t2+\frac{\sqrt{x}}{2(1-t)}\sqrt{4-3t+U(2-t)^2}
\end{gather*}
Since the first term is well understood, we concentrate on the second:
\begin{multline*}
	\frac{\sqrt{x}}{2(1-t)}\sqrt{4-3t+U(2-t)^2}
	\\=\sqrt{ux}\bigg(	1+\frac18 \left( 5+4U \right) \tau+\frac {1}{128}  
	( 71+136U	+64{U}^{2} ) {\tau}^{2}\\+\frac {1}{1024} ( 541+1596U	+1568{U}^{2}+512{U}^{3} ) {\tau}^{3}+\cdots\bigg)	=:\sqrt{ux}\cdot\Xi.
\end{multline*} 
With this expanded form $\Xi$, we have now our final formula, the expansion of $r_1$ into two factors:
\begin{align*}
	r_1=\frac1{1-t}=-\frac{1}{ux}\frac1{r_2}\frac1{r_3}
	=\Big(-\frac{1}{\sqrt{ux}}\frac t2+\Xi\Big)\Big(\frac{1}{\sqrt{ux}}\frac t2+\Xi\Big).
\end{align*}
These two factors do not have a combinatorial meaning, as it seems, but we can still stick to the $(3/2)$-ary tree notation, with the additional counting of middle edges.

\section{More about Motzkin paths}

\subsection*{Retakh's Motzkin paths}
V.~Retakh~\cite{EZ} introduced the following restricted class of Dyck paths: Peaks are only allowed on level 1 and on even-numbered levels. Here is an example, and the corresponding plane tree using the standard bijection. 
\begin{center}
	\begin{tikzpicture}[scale=0.6]
		
		\draw[step=1.cm,black,dotted] (-0.0,-0.0) grid (20.0,6.0);

		%	\draw[ thick] (0,0) to (20,0);

		\draw[thick] (0,0) -- (1,1) -- (2,0)-- (3,1)-- (4,2)-- (5,3)-- (6,4)-- (7,3)-- (8,4)-- (9,5)-- (10,6)-- (11,5)-- (12,4)-- (13,3)-- (14,2)-- (15,1)-- (16,0)-- (17,1)-- (18,0)-- (19,1)-- (20,0);
		
		\node at (1,1){$\bullet$};
		\node at (6,4){$\bullet$};
		\node at (10,6){$\bullet$};
		\node at (17,1){$\bullet$};
		\node at (19,1){$\bullet$};
	\end{tikzpicture}
	\end	{center}
	
	\begin{center}
		\begin{tikzpicture}[scale=0.6]
			\node (1) at (0, 0){$\bullet$};
			\node (2) at (-3, -1){$\bullet$};
			\node (3) at (-1, -1){$\bullet$};
			\node (4) at (1, -1){$\bullet$};
			\node (5) at (3, -1){$\bullet$};
			\node (6) at (-1, -2){$\bullet$};
			\node (7) at (-1, -3){$\bullet$};
			\node (8) at (-2, -4){$\bullet$};
			\node (9) at (0, -4){$\bullet$};	
			\node (10) at (0, -5){$\bullet$};
			\node (11) at (0, -6){$\bullet$};
			
			\node  at (-5, -1){level $1$};
			\node  at (-5, -4){level $4$};	
			\node  at (-5, -6){level $6$};		
			
			\draw[-] (1.center) to (2.center);
			\draw[-] (1.center) to (3.center);
			\draw[-] (1.center) to (4.center);
			\draw[-] (1.center) to (5.center);
			\draw[-] (3.center) to (6.center);
			\draw[-] (6.center) to (7.center);
			\draw[-] (8.center) to (7.center);
			\draw[-] (9.center) to (7.center);
			\draw[-] (9.center) to (10.center);
			\draw[-] (10.center) to (11.center);
		\end{tikzpicture}
		
	\end{center}
	
	\begin{center}
		\begin{tikzpicture}[scale=0.6]
			\node (1) at (0, 0){$\bullet$};
			\node (2) at (-3, -1){$\bullet$};
			\node (3) at (-1, -1){$\bullet$};
			\node (4) at (1, -1){$\bullet$};
			\node (5) at (3, -1){$\bullet$};
			
			\node (8) at (-2, -4){  };
			\node (9) at (0, -4){ };

			\draw[-] (1.center) to (2.center);
			\draw[-] (1.center) to (3.center);
			\draw[-] (1.center) to (4.center);
			\draw[-] (1.center) to (5.center);
			\draw[-] (3.center) to (8.center);
			\draw[-] (3.center) to (9.center);
			\draw[-] (8.center) to (9.center);
		\end{tikzpicture}
		
	\end{center}
	Ekhad and Zeilberger \cite{EZ} proved  recently that these restricted paths are enumerated by Motzkin numbers. 
	Recall that the generating function of the Motzkin numbers $M(z)$ according to length satisfies $M=1+zM+z^2M^2$ and thus
	\begin{equation*}
		M(z)=\frac{1-z-\sqrt{1-2z-3z^2}}{2z^2}.
	\end{equation*}
	
	Here, I want to present a few additional observations, also including the height of the paths (or the associated plane trees).
	First, we are going to confirm the connection to Motzkin paths:	
	
	Since the level 1 is somewhat special, we only consider  trees as symbolized by the triangle. We will use two generating functions, to deal with the odd/even situation. We have
	\begin{equation*}
		F=\frac{zG}{1-G}\quad\text{and}\quad G=\frac{z}{1-F}.
	\end{equation*}
	This is not too difficult to see, since the family $\mathscr{F}$ of trees as symbolized by the triangle does not contain a single node, so $F=zG+zG^2+\cdots$. However, the next generation $\mathscr{G}$ can contain a single node, and thus $G=z+zF+zF^2+\cdots$.
	
	Solving this (best by a computer) we find $F(z)=z^2M(z)$ and the total generating function (allowing sequences of single nodes between copies of $\mathscr{F}$) is
	\begin{equation*}
		\frac z{1-z}\sum_{r\ge0}\Big(\frac{F}{1-z}\Big)^r=zM(z),
	\end{equation*}
	as predicted. Recall that the number of nodes in trees is always one more than the half-length of the corresponding Dyck path.
	
	We will       compute the average height of such restricted paths, using singularity analysis of generating functions, as in \cite{FlOd90, FS}. Whether we define the height in terms of the maximal chain of edges resp.\ nodes only makes a difference of one,
	and we will only compute the average height according to the leading term of order $\sqrt n$. For readers who wish to see how more terms could be computed, at least in principle, we cite~\cite{Prodinger-ars}.
	
	The average height of planted plane trees (and subclasses of them) has been of central interest to combinatorialists and theoretical computer scientists alike since the seminal paper \cite{BrKnRi72}. The number of leaves (endnodes) is also one of the key parameters since Narayana \cite{Narayana}. We investigate it in the last section of this paper to see how the restrictions according to Retakh influence this parameter.

	\subsection*{The height}	
	
	Now we will use the substitution $z=\frac{v}{1+v+v^2}$, which occurred for the first time in \cite{Prodinger-three}, but has been used more recently in different models where Motzkin numbers are involved \cite{Prodinger-Deutsch, HHP, Baril-neu}.
	Motzkin numbers appear in \cite{OEIS} as sequence A001006.
	For example,  the generating function $M(z)$ of Motzkin numbers has the very simple form  $M(z)=1+v+v^2$ using this auxiliary variable.
	We define 
	\begin{equation*}
		G_{k+1}=\dfrac{z}{1-\dfrac{zG_k}{1-G_k}},\quad\text{with} \quad G_1=z.
	\end{equation*}
	There is a simple formula, viz.
	\begin{equation*}
		G_k=\frac v{1+v}\frac{1-v^{2k}}{1-v^{2k+1}}.
	\end{equation*}
	This is easy to prove by induction, which we will do for the convenience of the reader.
	The start is
	\begin{equation*}
		G_1=\frac v{1+v}\frac{1-v^{2}}{1-v^{3}}=\frac v{1+v}\frac{1+v}{1+v+v^{2}}=\frac{v}{1+v+v^{2}}=z.
	\end{equation*}
	And now
	\begin{align*}
		G_{k+1}&=\dfrac{z}{1-\dfrac{zG_k}{1-G_k}}=\dfrac{z(1-G_k)}{1-(1+z)G_k}=\frac{v}{1+v+v^{2}}\dfrac{1-\frac v{1+v}\frac{1-v^{2k}}{1-v^{2k+1}}}{1-\frac{(1+v)^2}{1+v+v^2}\frac v{1+v}\frac{1-v^{2k}}{1-v^{2k+1}}}\\
		&=\frac{v}{1+v}\dfrac{1+v-v\frac{1-v^{2k}}{1-v^{2k+1}}}{1+v+v^{2}-v(1+v)\frac{1-v^{2k}}{1-v^{2k+1}}}\\&
		=\frac{v}{1+v}\dfrac{(1+v)(1-v^{2k+1})-v(1-v^{2k})}{(1+v+v^{2})(1-v^{2k+1})-v(1+v)(1-v^{2k})}\\
		&=\frac{v}{1+v}\dfrac{1-v^{2k+2}}{1-v^{2k+3}},
	\end{align*}
	as claimed. From this we also get
	\begin{equation*}
		F_k=\frac{zG_k}{1-G_k}=\frac{v^2}{1+v+v^2}\frac{1-v^{2k}	}{1-v^{2k+2}}.
	\end{equation*}
	For $k\ge1$, $F_k$ is the generating function of trees (like in the triangle) of height $\le 2k$. 
	
	Note that the height is currently counted in terms of nodes;
	\begin{equation*}
		F_1=\frac{z^2}{1-z},
	\end{equation*} 
	which describes a root with $\ell\ge1$ leaves attached to the root.
	
	Now we incorporate the irregular beginning of the tree and compute
	\begin{equation*}
		\frac z{1-z}\sum_{r\ge0}\Big(\frac{F_h}{1-z}\Big)^r=\frac z{1-z}\dfrac1{1-\dfrac{F_h}{1-z}}=v\frac{1-v^{2h+2}}{1-v^{2h+4}}.
	\end{equation*}
	From here onwards it seems to be more natural to define the height of the whole tree in terms of the number of \emph{edges}, 
	and then the quantity we just derived is the generating function of all trees with height $\le 2h$, for $h\ge1$. Note that the limit
	$h\to\infty$ gives us simply $v=zM(z)$, which is consistent. There is also a contribution of trees of height $\le 1$, namely 
	$\frac{z}{1-z}=\frac{v}{1+v^2}$, but this term is, when we compute the average height, irrelevant and only contributes to the error term, as we only compute  the leading term, which is of order $\sqrt n$.
	
	So, apart from normalization, we are led to investigate
	\begin{align*}
		\sum_{h\ge1}&2h\bigg[v\frac{1-v^{2h+2}}{1-v^{2h+4}}-v\frac{1-v^{2h}}{1-v^{2h+2}}\bigg]
		=2v(1-v^{-2})\sum_{h\ge1}h\bigg[\frac{v^{2h+4}}{1-v^{2h+4}}-\frac{v^{2h+2}}{1-v^{2h+2}}\bigg]\\
		&=2v(1-v^{-2})\sum_{h\ge0}h\frac{v^{2h+4}}{1-v^{2h+4}}
		-2v(1-v^{-2})\sum_{h\ge0}(h+1)\frac{v^{2h+4}}{1-v^{2h+4}}\\
		&=-2v+\frac{2(1-v^{2})}{v}\sum_{h\ge1}\frac{v^{2h}}{1-v^{2h}}.
	\end{align*}
	
	Note that we could get explicit coefficients from here, using trinomial coefficients, $\binom{n,3}{k}=[v^k](1+v+v^2)^n$ (notation from \cite{Comtet-book}). To show the reader how this works, we compute
	\begin{align*}
		[z^{n+1}]&\frac{1-v^{2}}{v}\sum_{h\ge1}\frac{v^{2h}}{1-v^{2h}}
		=\frac1{2\pi i}\oint\frac{dz}{z^{n+2}}\frac{1-v^{2}}{v}\sum_{h\ge1}\frac{v^{2h}}{1-v^{2h}}\\
		&=\frac1{2\pi i}\oint dv(1-v^2)^2\frac{(1+v+v^2)^{n}}{v^{n+3}}\sum_{h\ge1}\sum_{h,k\ge1}v^{2hk}\\
		&=[v^{n+2}](1-2v^2+v^4)\sum_{h\ge1}d(h)v^{2h}(1+v+v^2)^{n}\\
		&=\sum_{h\ge1}d(h)\bigg[\binom{n,3}{n+2-2h}-2\binom{n,3}{n-2h}+\binom{n,3}{n-2-2h}\bigg].
	\end{align*}
	Note that $d(h)$ is the number of divisors of $h$. We will, however, not use this explicit form.
	The expression as derived before,
	\begin{equation*}
		-2v+\frac{2(1-v^{2})}{v}\sum_{h\ge1}\frac{v^{2h}}{1-v^{2h}},
	\end{equation*}
	has to be expanded around $v=1$, which is a standard application of the Mellin transform. Details are worked out in
	\cite{HPW}, for example:
	\begin{equation*}
		\sum_{h\ge1}\frac{v^{2h}}{1-v^{2h}}=\sum_{k\ge1}d(k)v^{2k}\sim-\frac{\log(1-v^2)}{1-v^2}
		\sim-\frac{\log(1-v)}{2(1-v)}.
	\end{equation*}
	Note again that $d(k)$ is the number of divisors of $k$. Consequently
	\begin{equation*}
		-2v+\frac{2(1-v^{2})}{v}\sum_{h\ge1}\frac{v^{2h}}{1-v^{2h}}\sim
		-2\log(1-v).
	\end{equation*}
	We have $1-v\sim \sqrt 3\sqrt{1-3z}$, and $z=\frac13$ is the relevant singularity when discussing Motzkin numbers. We can continue
	\begin{equation*}
		-2v+\frac{2(1-v^{2})}{v}\sum_{h\ge1}\frac{v^{2h}}{1-v^{2h}}\sim -\log(1-3z).
	\end{equation*}
	The coefficient of $z^n$ in this is $\frac{3^n}{n}$. This has to be divided by
	\begin{equation*}
		[z^n]zM(z)=[z^{n-1}]M(z)\sim \frac{3^{n+\frac12}}{2\sqrt{\pi }n^{3/2}},
	\end{equation*}
	with the final result for the average height of  restricted Dyck paths (\`a la Retakh):
	\begin{equation*}
		\sim 2\sqrt{\frac{\pi n}{3}}.
	\end{equation*}
	Recall \cite{Prodinger-three} that the average height of Motzkin paths of length $n$ is asymptotic to
	\begin{equation*}
		\sqrt{\frac{\pi n}{3}}.
	\end{equation*}
	
	\subsection*{The number of leaves}
	We can use a second variable, $u$, to count the number of leaves. Then we have
	\begin{equation*}
		F(z,u)=\frac{zG(z,u)}{1-G(z,u)}\quad\text{and}\quad G(z,u)=zu+\frac{zF(z,u)}{1-F(z,u)},
	\end{equation*}
	which leads to
	\begin{multline*}
		F(z,u)= \\{\frac {1-zu-{z}^{2}+{z}^{2}u-\sqrt {1-2zu-2{z}^{2}-2{z}^{
						2}u+{z}^{2}{u}^{2}-2{z}^{3}u+2{z}^{3}{u}^{2}+{z}^{4}-2{z}^{4}u+{
						z}^{4}{u}^{2}}}{2(1-zu+z)}}.
	\end{multline*}
	Bringing the irregular beginning also into the game leads to
	\begin{align*}
		\frac{z}{1-zu}\sum_{r\ge0}\Big(\frac{F}{1-zu}\Big)^r+zu-z.
	\end{align*}
	This is an ugly expression that we do not display here. However, we can compute the average number of leaves, by differentiation w.r.t. $u$, followed by setting $u=1$:
	\begin{equation*}
		R:={\frac {v \left( 1+v \right)   ( 1-v+2{v}^{2}-{v}^{3}) }{
				(1-v)  \left( 1+v+{v}^{2} \right) }}.
	\end{equation*}
	The coefficient of $z^n$ in this can be expressed in terms of trinomial coefficients, if needed. However, we only compute an asymptotic formula, to keep this section short. Expanding around $v=1$, we find
	\begin{equation*}
		R\sim\frac23\frac{1}{1-v}\sim\frac23\frac{1}{\sqrt3\sqrt{1-3z}},
	\end{equation*}
	and thus
	\begin{equation*}
		[z^n]R\sim\frac23\frac{1}{\sqrt3}3^n\frac1{\sqrt{\pi n}}.
	\end{equation*}
	We divide this again by
	\begin{equation*}
		\frac{3^{n+\frac12}}{2\sqrt{\pi} n^{3/2}}
	\end{equation*}
	with the result
	\begin{equation*}
		\frac49n,
	\end{equation*}
	which is the asymptotic number of leaves in a Retakh tree of size $n$. Recall that for unrestricted planar trees, the result is
	$\frac n2$, which is a folklore result using Narayana numbers. So the constant in the restricted case, $\frac49$, is a bit smaller than 
	$\frac12$.
	
	With some effort, more precise approximations could be obtained, as well as the variance. This might be a good project for a student.

\section{The amplitude of Motzkin paths}

A Motzkin path consists of up-steps, down-steps, and horizontal steps, see sequence A091965 in \cite{OEIS} and the references given there. As Dyck paths, they start at the origin and end, after $n$ steps 
again at the $x$-axis, but are not allowed to go below the $x$-axis. The height of a Motzkin path is the highest $y$-coordinate that the path
reaches. The average height of such random paths of length $n$ was considered in an early paper \cite{Prodinger-three}, it is asymptotically given by
$\sqrt{\frac{\pi n}{3}}$. 

In the recent paper \cite{irene} an interesting new concept was introduced: the \emph{amplitude}. It is basically twice the height, but with a twist.
If there exists a horizontal  step on level $h$, which is the height, the amplitude is $2h+1$, otherwise it is $2h$. To clarify the concept, we  
created a list of all 9 Motzkin paths of length 4 with height and amplitude given.

\begin{center}
	
	\begin{table}[h]
		\begin{tabular}{c | c | c  |c}
			\text{Motzkin path  }&\text{horizontal   on maximal level}&\text{height}&\text{amplitude}\\
			\hline\hline
			\begin{tikzpicture}[scale=0.4]
				
				\draw[ultra thick] (0,0) to (1,0) to (2,0) to (3,0) to (4,0)  ;
			\end{tikzpicture}
			& \text{Yes}& 0& 1\\
			
			\hline
			
			\begin{tikzpicture}[scale=0.4]
				
				\draw[ ultra thick] (0,0) to (1,0) to (2,0) to (3,1) to (4,0)  ;
			\end{tikzpicture}
			& \text{No}& 1& 2\\
			
			\hline
			\begin{tikzpicture}[scale=0.4]
				
				\draw[ ultra thick] (0,0) to (1,1) to (2,0) to (3,1) to (4,0)  ;
			\end{tikzpicture}
			& \text{No}& 1& 2\\
			
			\hline
			\begin{tikzpicture}[scale=0.4]
				
				\draw[ ultra thick] (0,0) to (1,0) to (2,1) to (3,1) to (4,0)  ;
			\end{tikzpicture}
			& \text{Yes}& 1& 3\\
			
			\hline
			\begin{tikzpicture}[scale=0.4]
				
				\draw[ ultra thick] (0,0) to (1,1) to (2,1) to (3,1) to (4,0)  ;
			\end{tikzpicture}
			& \text{Yes}& 1& 3\\
			\hline
			\begin{tikzpicture}[scale=0.4]
				
				\draw[ ultra thick] (0,0) to (1,1) to (2,1) to (3,0) to (4,0)  ;
			\end{tikzpicture}
			& \text{Yes}& 1& 3\\
			\hline
			\begin{tikzpicture}[scale=0.4]
				
				\draw[ ultra thick] (0,0) to (1,1) to (2,0) to (3,0) to (4,0)  ;
			\end{tikzpicture}
			& \text{No}& 1& 2\\
			\hline 
			\begin{tikzpicture}[scale=0.4]
				
				\draw[ ultra thick] (0,0) to (1,0) to (2,1) to (3,0) to (4,0)  ;
			\end{tikzpicture}
			& \text{No}& 1& 2\\
			
			\hline
			\begin{tikzpicture}[scale=0.4]
				
				\draw[ ultra thick] (0,0) to (1,1) to (2,2) to (3,1) to (4,0)  ;
			\end{tikzpicture}
			& \text{No}& 2& 4\\
			\hline
		\end{tabular}

	\end{table}	
\end{center}

The goal of this paper is to investigate this new parameter; in the next section, generating functions will be given, in the following section explicit enumerations,
involving trinomial coefficients  $\binom{n,3}{k}=[t^k](1+t+t^2)^n$ (notation following Comtet's book \cite{Comtet-book}).
In the last section, the intuitive result that the average amplitude is about twice the average height, is confirmed, and then
it will be shown, that, asymptotically, there are about as many Motzkin paths with/without horizontal steps on the maximal level.

\subsection*{Generating functions}

Let 
\begin{equation*}
	M^{\le h}(z)=\sum_{n\ge0}[\text{number of Motzkin paths of length $n$ and height $\le h$}]z^n.
\end{equation*}
For the computation, let $f_i=f_i(z)$ be the generating function of Motzkin-like paths, bounded by height $h$, but ending at level $i$.
Distinguishing the last step, we get
\begin{equation*}
	f_0=1+zf_0+zf_1,\quad f_i=zf_{i-1}+zf_i+zf_{i+1}\ \ \text{for}\ 0<i<h, \quad f_h=zf_{h-1}+zf_h.
\end{equation*}
This system is best written in matrix form:
\begin{equation*}
	\begin{pmatrix}
		1-z&-z&0&\dots\\
		-z&1-z&-z&0&\dots\\
		0&-z&1-z&-z&0&\dots\\
		\ddots&\ddots&\ddots\\
		\phantom{0}&\phantom{0}&\phantom{0}&\phantom{0}&-z&1-z
	\end{pmatrix}
	\begin{pmatrix}
		f_0\\f_1\\f_2\\ \vdots\\ f_h
	\end{pmatrix}
	=\begin{pmatrix}
		1\\0\\0\\ \vdots\\ 0
	\end{pmatrix}
\end{equation*}
Let $D_n$ be the determinant of the system matrix, with $n$ rows and columns. Using Cramer's rule to solve a linear system, one finds
\begin{equation*}
	M^{\le h}(z)=f_0=\frac{D_h}{D_{h+1}}.
\end{equation*}
Expanding the determinant along the first row, we get the recursion $D_n=(1-z)D_{n-1}-z^2D_{n-2}$, and $D_1=1-z$, $D_0=1$. Solving,
\begin{equation*}
	D_n=\frac1{\sqrt{1-2z-3z^2}}\bigg[\biggl(\frac{1-z+\sqrt{1-2z-3z^2}}{2}\biggr)^{n+1}-\biggl(\frac{1-z-\sqrt{1-2z-3z^2}}{2}\biggr)^{n+1}\bigg].
\end{equation*}
If one deals with enumeration of Motzkin-like objects, the substitution $z=\frac{v}{1+v+v^2}$ makes the expressions prettier:
\begin{equation*}
	D_n=\frac1{1-v^2}\frac{1-v^{2n+2}}{(1+v+v^2)^n}
\end{equation*}
and further
\begin{equation*}
	M^{\le h}(z)=f_0=\frac{D_h}{D_{h+1}}=(1+v+v^2)\frac{1-v^{2h+2}}{1-v^{2h+4}}.
\end{equation*}

Now let $N^{\le h}(z)$ be the generating function of Motzkin paths of height $\le h$, but where horizontal steps on level $h$ are forbidden. The system to
compute this is quite similar:
\begin{equation*}
	\begin{pmatrix}
		1-z&-z&0&\dots\\
		-z&1-z&-z&0&\dots\\
		0&-z&1-z&-z&0&\dots\\
		\ddots&\ddots&\ddots\\
		\phantom{0}&\phantom{0}&\phantom{0}&\phantom{0}&-z&\boldsymbol{1}
	\end{pmatrix}
	\begin{pmatrix}
		g_0\\g_1\\g_2\\ \vdots\\ g_h
	\end{pmatrix}
	=\begin{pmatrix}
		1\\0\\0\\ \vdots\\ 0
	\end{pmatrix}
\end{equation*}
The only difference in the matrix is the entry in the last row, written in boldface. Let $D_n^*$ be the determinant of this system matrix, with $n$ rows and columns. Again,
\begin{equation*}
	N^{\le h}(z)=g_0=\frac{D_h^*}{D_{h+1}^*}.
\end{equation*}
Expanding the matrix along the last row, we find
\begin{equation*}
	D_n^*=D_{n-1}-z^2D_{n-2}=\frac1{1-v}\frac{1-v^{2n+1}}{(1+v+v^2)^n}
\end{equation*}
and
\begin{equation*}
	N^{\le h}(z)=g_0=\frac{D_h^*}{D_{h+1}^*}=(1+v+v^2)\frac{1-v^{2h+1}}{1-v^{2h+3}}.
\end{equation*}
Now let us consider
\begin{equation*}
	M^{\le h}(z)-N^{\le h}(z).
\end{equation*}
There is obviously a lot of cancellation going on. The objects which are still counted, have height $=h$, and \emph{have} horizontal steps on level $h$. That is one of the two quantities that we wanted to compute, and we get
\begin{align*}
	\text{Horiz}_h(z)&=(1+v+v^2)\frac{1-v^{2h+2}}{1-v^{2h+4}}-(1+v+v^2)\frac{1-v^{2h+1}}{1-v^{2h+3}}\\&=
	(1+v+v^2)(1-v^{-2})\bigg[\frac{v^{2h+4}}{1-v^{2h+4}}-\frac{v^{2h+3}}{1-v^{2h+3}}\bigg].
\end{align*}

Similarly, considering $N^{\le h}(z)-M^{\le h-1}(z)$, we find that only objects are counted that have height $=h$, and \emph{no} horizontal steps on level $h$.
Thus
\begin{align*}
	\text{No-Horiz}_h(z)&=(1+v+v^2)\bigg[\frac{1-v^{2h+1}}{1-v^{2h+3}}-\frac{1-v^{2h}}{1-v^{2h+2}}\bigg]\\&=
	(1+v+v^2)(1-v^{-2})\bigg[\frac{v^{2h+3}}{1-v^{2h+3}}-\frac{v^{2h+2}}{1-v^{2h+2}}\bigg].
\end{align*}
As a check, we get
\begin{align*}
	\text{Horiz}_h(z)+\text{No-Horiz}_h(z)&=(1+v+v^2)(1-v^{-2})\bigg[\frac{v^{2h+4}}{1-v^{2h+4}}-\frac{v^{2h+2}}{1-v^{2h+2}}\bigg]
	\\*&=M^{\le h}(z)-M^{\le h-1}(z)=M^{= h}(z).
\end{align*}

\subsection*{Explicit enumerations}

All our generating functions contain the term 
\begin{equation*}
	(1+v+v^2)(1-v^{-2})\frac{v^{2h+a}}{1-v^{2h+a}}=(1+v+v^2)(1-v^{-2})\sum_{k\ge1}v^{k(2h+a)}
\end{equation*}
for various values of $a$. We show how to compute the coefficient of $z^n$ in this. It will be done using contour integration. The contour is a small circle in the $z$-plane or
$v$-plane.
\begin{align*}
	[z^n]&(1+v+v^2)(1-v^{-2})\frac{v^{2h+a}}{1-v^{2h+a}}
	=\frac1{2\pi i}\oint \frac{dz}{z^{n+1}}(1+v+v^2)(1-v^{-2})\sum_{k\ge1}v^{k(2h+a)}\\
	&=\frac1{2\pi i}\oint \frac{dv(1-v^2)}{(1+v+v^2)^2}\frac{(1+v+v^2)^{n+1}}{v^{n+1}}(1+v+v^2)(1-v^{-2})\sum_{k\ge1}v^{k(2h+a)}\\
	&=-\sum_{k\ge1}[v^{n+2-k(2h+a)}](1-v^2)^2(1+v+v^2)^{n}\\
	&=-\sum_{k\ge1}\bigg[\binom{n,3}{n+2-k(2h+a)}-2\binom{n,3}{n-k(2h+a)}+\binom{n,3}{n-2-k(2h+a)}\bigg].
\end{align*}
In this computation, we used again the notion of \emph{trinomial} coefficients:
\begin{equation*}
	\binom{n,3}{k}=[t^k](1+t+t^2)^n.
\end{equation*}

\subsection*{The average amplitude}

Here we compute:
\begin{align*}
	\sum_{h\ge0}&(2h+1)\text{Horiz}_h(z)+\sum_{h\ge1}(2h)\text{No-Horiz}_h(z)\\
	&=(1+v+v^2)(1-v^{-2})\sum_{h\ge0}(2h+1)\bigg[\frac{v^{2h+4}}{1-v^{2h+4}}-\frac{v^{2h+3}}{1-v^{2h+3}}\bigg]
	\\&+(1+v+v^2)(1-v^{-2})\sum_{h\ge0}(2h)\bigg[\frac{v^{2h+3}}{1-v^{2h+3}}-\frac{v^{2h+2}}{1-v^{2h+2}}\bigg]\\
	&=-(1+v+v^2)(1-v^{-2})\sum_{h\ge0}\frac{v^{2h+3}}{1-v^{2h+3}}\\
	&+(1+v+v^2)(1-v^{-2})\sum_{h\ge0}(2h+1)\frac{v^{2h+4}}{1-v^{2h+4}}\\&-
	(1+v+v^2)(1-v^{-2})\sum_{h\ge0}(2h+2)\frac{v^{2h+4}}{1-v^{2h+4}}\\
	&=-(1+v+v^2)(1-v^{-2})\sum_{h\ge0}\frac{v^{2h+3}}{1-v^{2h+3}}-(1+v+v^2)(1-v^{-2})\sum_{h\ge1}\frac{v^{2h+2}}{1-v^{2h+2}}\\
	&=-(1+v+v^2)(1-v^{-2})\sum_{h\ge0}\frac{v^{2h+1}}{1-v^{2h+1}}+(1+v+v^2)(1-v^{-2})\frac{v}{1-v}
	\\&-(1+v+v^2)(1-v^{-2})\sum_{h\ge1}\frac{v^{2h}}{1-v^{2h}}+(1+v+v^2)(1-v^{-2})\frac{v^{2}}{1-v^{2}}\\
	&=-(1+v+v^2)(1-v^{-2})\sum_{h\ge0}\frac{v^{2h+1}}{1-v^{2h+1}}
	\\&-(1+v+v^2)(1-v^{-2})\sum_{h\ge1}\frac{v^{2h}}{1-v^{2h}}-\frac{(1+2v)(1+v+v^2)}{v}\\
	\\&=-(1+v+v^2)(1-v^{-2})\sum_{h\ge1}\frac{v^{h}}{1-v^{h}}-\frac{(1+2v)(1+v+v^2)}{v}.
\end{align*}

To find asymptotics from here, we need the local expansion of this generating function around $v\sim1$, or, equivalently, $z\sim \frac13$.
See \cite{FGD} for more explanations of how this method works; it also involves singularity analysis of generating functions \cite{FlOd90}.
While this combined approach of Mellin transform and singularity analysis has been used for more than 30 years, we would like to cite
\cite{HPW}, where many technical details have been worked out in a similar instance. For example, the expansion of the sum that we need is 
derived there. We combine all the expansions:
\begin{equation*}
	\Big(6(1-v)+\cdots\Big)\Big(-\frac{\log(1-v)}{1-v}+\frac{\gamma}{1-v}+\cdots\Big)-9+6(1-v)+\cdots=-6\log(1-v)+\cdots
\end{equation*}
Translating this into the $z$-world means $1-v\sim \sqrt3\sqrt{1-3z}$, 
and then
\begin{equation*}
	[z^n]\Big(-6\log(1-v)\Big)\sim [z^n]\Big(-6\log\sqrt{1-3z}\Big)=3\frac{3^n}{n}=\frac{3^{n+1}}{n}.
\end{equation*}
The total number of Motzkin paths is
\begin{equation*}
	[z^n]\frac{1-z-\sqrt{1-2z-3z^2}}{2z^2}\sim [z^n]\Big(3-3\sqrt3\sqrt{1-3z}\Big)\sim 3\sqrt3\frac{3^{n}}{2\sqrt\pi n^{3/2}}.
\end{equation*}
Taking quotients, we get the average amplitude of random Motzkin paths of length $n$, which is asymptotic to
\begin{equation*}
	2\sqrt{\frac{\pi n}{3}}.
\end{equation*}
This is about twice as much as the average height of random Motzkin paths, which is what we expected.

Now we consider
\begin{align*}
	\sum_{h\ge0}&\text{No-Horiz}_h(z)=
	(1+v+v^2)(1-v^{-2})\sum_{h\ge0}\bigg[\frac{v^{2h+3}}{1-v^{2h+3}}-\frac{v^{2h+2}}{1-v^{2h+2}}\bigg]\\
	&=(1+v+v^2)(1-v^{-2})\sum_{h\ge1}\bigg[\frac{v^{2h-1}}{1-v^{2h-1}}-\frac{v^{2h}}{1-v^{2h}}\bigg]
	+\frac{(1+v)(1+v+v^2)}{v}\\
	&=(1+v+v^2)(1-v^{-2})\sum_{h\ge1}\frac{v^{h}}{1-v^{h}}(-1)^{h-1}
	+\frac{(1+v)(1+v+v^2)}{v}.
\end{align*}
This needs to be expanded about $v=1$:
\begin{equation*}
	(1+v+v^2)(1-v^{-2})\sim-6(1-v)-3(1-v)^2+\cdots;
\end{equation*}
\begin{equation*}
	\frac{(1+v)(1+v+v^2)}{v}\sim6-3(1-v)+\cdots;
\end{equation*}
for the remaining sum we need the Mellin transform \cite{FGD}. Set $v=e^{-t}$, and transform
\begin{equation*}
	\mathscr{M}\sum_{h\ge1}\frac{v^{h}}{1-v^{h}}(-1)^{h-1}
	=\mathscr{M}\sum_{h,k\ge1}e^{-thk}(-1)^{h-1}=
	\Gamma(s)\zeta(s)^2(1-2^{1-s}).
\end{equation*}
By the Mellin  inversion formula:
\begin{equation*}
	\sum_{h\ge1}\frac{v^{h}}{1-v^{h}}(-1)^{h-1}=
	\frac1{2\pi i}\int_{2-i\infty}^{2+i\infty}\Gamma(s)\zeta(s)^2(1-2^{1-s})t^{-s}ds.
\end{equation*}
The line of integration will be shifted to the left, and the collected residues constitute the expansion about $t=0$:
\begin{align*}
	\sum_{h\ge1}\frac{v^{h}}{1-v^{h}}(-1)^{h-1}\sim \frac{\log 2}{t}-\frac14+\frac{t}{48}\sim\frac{\log2}{1-v}-\frac{\log2}{2}-\frac14+\cdots.
\end{align*}
Combining,
\begin{equation*}
	\sum_{h\ge0}\text{No-Horiz}_h(z)\sim6-6\log 2-\frac32(1-v)+\cdots.
\end{equation*}
But
\begin{equation*}
	\frac{1-z-\sqrt{1-2z-3z^2}}{2z^2}=1+v+v^2\sim 3-3(1-v)+\cdots.
\end{equation*}
Comparing the coefficients of $1-v$, we find that asymptotically about half of the Motzkin paths belong to the `No-horizontal' and
about half of the Motzkin paths belong to the `horizontal' class. Again, this is intuitively clear, once one sees the generating functions of both families.

%\section{Oscillations of Dyck paths}

%To be written

\section{Oscillations in Dyck paths revisited}

Rainer Kemp's paper \cite{Kemp-oscillations} was unfortunately largely overlooked. An extension was published quickly \cite{KP-hyper},
and then it fell into oblivion. We want to come back to this gem, with modern methods, in particular, the kernel method
and singularity analysis. Kemp was interested in peaks and valleys of Dyck paths, which he called \textsc{max}-turns and
\textsc{min}-turns, probably motivated by Computer Science applications. The peaks/valleys are enumerated from left to right, and
the interest is the height of them, for a fixed number and the length of the Dyck path going to infinity. In the corresponding
ordered (plane) tree, the peaks correspond to the leaves.

Very precise information is available for leaves of binary trees \cite{Kirschen-leaves, Gutjahr, PP1, PP2} but the situation is
a bit different for Dyck paths since the number of peaks/valleys  isn't directly related to the length of the Dyck path. (Narayana numbers enumerate them.)
Kemp's results in a nutshell are: The average height of the $m$-th peak/valley tends to a constant for $n$ going to infinity (the length of the random Dyck path).
This constant is $\sim 4\sqrt{2m/\pi }$, and the difference between the height of the peak and the valley is about 2, with more terms being available in principle.

\subsection*{A trivariate generating function for heights of valleys}

The goal is to derive an expression for $\Phi(u)=\Phi(u;z,w)$, where $z$ is used for the length of the path, $w$ for the enumeration of the
$valleys$ ($w^m$ corresponds to the the $m$-th valley), and $u$ is used to record the height of the $m$-th (and last) valley of a partial Dyck path
(the path does not need to return to the $x$-axis). We could think about it continued in any possible fashion, as in the following figure.

\begin{figure}[h]
	\begin{center}
		\begin{tikzpicture}[scale=0.5]
			\draw (0,0) -- (25,0);
			\draw (0,0) -- (0,10);
			\draw[thick](0,0)--(3,5)--(6,2)--(10,9)--(12,7)--(14,8)--(16,3);
			%\draw[thick](5,6)--(6,9);
			\draw [decorate,decoration=snake,thick]   (16,3) to [out=50,in=130] (24,0);
			%\draw[draw=blue, decorate,decoration=snake,thick] (16,3) arc (-30:150:5cm);
			\draw[thick, red,latex-latex](16,3)--(16,0);

			\fill (16,3) circle (0.20cm);
			
			%\draw (-0.3,6) --(0.3,6);
			%\draw (-0.3,9) --(0.3,9);
			%\node[thick] at (-1+0.4,9){$j$};

			%		\draw (-0.3,9) --(0.3,9);
		\end{tikzpicture}
	\end{center}
	\caption{The third valley at level $j$.}
\end{figure}
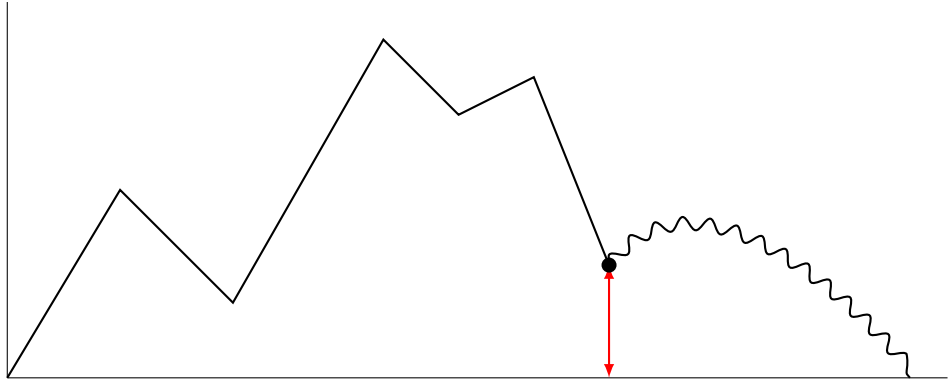

Our goal is, as often, to use the adding-a-new-slice technique, namely adding another `mountain' (a maximal sequence of up-steps, followed by a maximal
sequence of down-steps), or going from the $m$-th valley to the $(m+1)$-st valley.

We investigate what can happen to $u^j$:
\begin{equation*}
	\sum_{l\ge1}\sum_{i=1}^{j+l}z^lu^{j+l}z^iu^{-i}.
\end{equation*}
Working this out, the following substitution is essential for our problem:
\begin{equation*}
	u^j\longrightarrow \frac{z^{2}u^k}{(u-z)(1-zu^k)}u^j-\frac{z^{k+2}}{(u-z)(1-z^{k+1})}z^j.
\end{equation*}
Working this into a generating function of the type
\begin{equation*}
	\Phi(u)=\sum_{m\ge0}w^m\varphi_m(u),
\end{equation*}
where the variable $w$ is used to keep track of the number of mountains, we find from the substitution
\begin{equation*}
	\Phi(u)=1+\frac{wz^{2}u}{(u-z)(1-zu)}\Phi(u)-\frac{wz^{3}}{(u-z)(1-z^{2})}\Phi(z),
\end{equation*}
where 1 stands for the empty path having no mountains.
Rearranging, 
\begin{equation*}
	\Phi(u)\frac{z(u-s_1)(u-s_2)}{(u-z)(zu-1)}=1-\frac{wz^{3}}{(u-z)(1-z^{2})}\Phi(z),
\end{equation*}
and
\begin{equation*}
	\Phi(u)=\frac{(zu-1)}{z(u-s_1)(u-s_2)}\bigg[u-z-\frac{wz^{3}}{(1-z^{2})}\Phi(z)\bigg].
\end{equation*}
Here,
\begin{equation*}
	s_2={\frac {{z}^{2}+1-w{z}^{2}-\sqrt {{z}^{4}-2\,{z}^{2}-2\,{z}^{4}
				w+1-2\,w{z}^{2}+{w}^{2}{z}^{4}}}{2z}}\quad\text{and}\quad s_1=\frac1{s_2}.
\end{equation*}
In the spirit of the kernel method, the factor $u-s_2$ is `bad' and must cancel out. That leads first to
\begin{equation*}
	\Phi(z)=\frac{(1-z^2)(s_2-z)}{wz^3}
\end{equation*}
and further to
\begin{equation*}
	\Phi(u)=\frac{(zu-1)}{z(u-s_1)}=\frac{s_2(1-zu)}{z(1-us_2)}.
\end{equation*}
From this it is easy to read off coefficients:
\begin{equation*}
	[u^j]\Phi(u)=[u^j]\frac{s_2(1-zu)}{z(1-us_2)}=\frac1zs_2^{j+1}-s_2^{j}.
\end{equation*}
Note that setting $w=1$ ignores the number of mountains, and the generating function would then
be enumerating partial Dyck paths ending on level $j$ with a down-step. The answer could then be derived by combinatorial means as well.

For Kemp's problem,  we need
\begin{equation*}
	S=\sum_{j\ge0}j\Big(\frac1zs_2^{j+1}-s_2^{j}\Big)\cdot \Big(\frac1zs_2^{j+1}-s_2^{j}\Big)\Big|_{w=1}.
\end{equation*}
The factor $j$ comes in because of the average value of the height of the valley, the first factor is what we just worked out,
and the third factor is the rest, which, when read from right to left, is just what we discussed, since the number of mountains in the rest is
irrelevant. Thanks to Computer Algebra (not available when Kemp worked on the oscillations), we get
\begin{equation*}
	S=4 {\frac { \left( -3 z+W_1 z-W_1+1 \right)  \left( -W_2+wzW_2+1+{z}^{2}{w}^{2}-w{z}^{2}-2 wz-z \right) }{z
			\left( -3 z-W_1 z+1-W_1-wz+wzW_1-W_2+W_2 W_1 \right) ^{2}}}
\end{equation*}
with
\begin{align*}
	W_1&=\sqrt{1-4z}\quad\text{and}\quad
	W_2=\sqrt{z^2-2z-2z^2w+1-2wz+w^2z^2}.
\end{align*}
Note carefully that $z^2$ was replaced by $z$, since Dyck paths (returning to the $x$-axis) must have an even number of steps.
Their enumeration is classical:
\begin{equation*}
	D(z)=\frac{1-\sqrt{1-4z}}{2z}\sim 2-2\sqrt{1-4z},
\end{equation*}
for $z$ close to the (dominant) singularity $z=\frac14$. We are in the regime of the subcritical case (\cite{FS}, Section IX-3).
The function $S$ has a similar local expansion:
\begin{equation*}
	S\sim \mathsf{constant}_1-\mathsf{constant}_2\sqrt{1-4z},
\end{equation*}
and the function $\frac{-\mathsf{constant}_2}{-2}$ is the resulting generating function. Working out the details,
\begin{align*}
	S&\sim{\frac {w+\sqrt { \left(1-w \right)  \left( 9-w \right) }-3}{-1+w}}\\&
	-\sqrt{1-4z}\bigg({\frac {{w}^{2}+2w-3+(1+w)\sqrt { \left( 1-w \right)  \left( 9-w \right) }}{ \left( 1-w		\right) ^{2}}}\bigg)+\cdots
\end{align*}
Eventually we are led to
\begin{equation*}
	\mathsf{Valley}(w):={\frac {{w}^{2}+2w-3+(1+w)\sqrt { \left( 1-w \right)  \left(9-w \right) }}{ 2( 1-w	) ^{2}}}.
\end{equation*}
To say it again, the coefficient of $w^m$ in this is the average value of the $m$-th valley in `very long' Dyck path. To say more about it, we
can use singularity analysis again and expand (this time around $w=1$, which is dominant):
\begin{equation*}
	\mathsf{Valley}(w)\sim{\frac {2\sqrt {2}}{ \left( 1-w \right) ^{3/2}}}-\frac{2}{1-w}-{\frac {7}{8}} {\frac {\sqrt {2}}{\sqrt {1-w}}}.
\end{equation*}
The traditional translation theorems \cite{FlOd90} lead to the average value of the height of the $m$-th valley:
\begin{equation*}
	4\sqrt2\sqrt\frac{m}{\pi}-2+\frac{5\sqrt2}{8\sqrt{\pi m}}+\cdots.
\end{equation*}

\subsection*{From valleys to peaks.} 

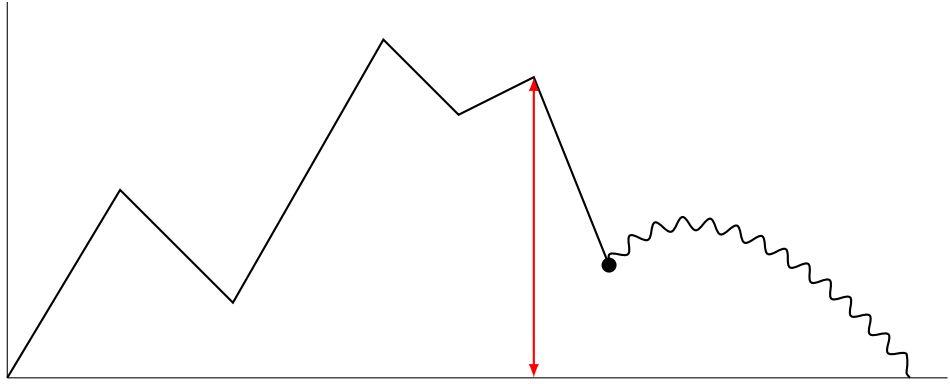
\begin{figure}[h]
	\begin{center}
		\begin{tikzpicture}[scale=0.5]
			\draw (0,0) -- (25,0);
			\draw (0,0) -- (0,10);
			\draw[thick](0,0)--(3,5)--(6,2)--(10,9)--(12,7)--(14,8)--(16,3);
			%\draw[thick](5,6)--(6,9);
			\draw [decorate,decoration=snake,thick]   (16,3) to [out=50,in=130] (24,0);
			%\draw[draw=blue, decorate,decoration=snake,thick] (16,3) arc (-30:150:5cm);
			\draw[thick, red,latex-latex](14,8)--(14,0);

			\fill (16,3) circle (0.20cm);
			
			%\draw (-0.3,6) --(0.3,6);
			%\draw (-0.3,9) --(0.3,9);
			%\node[thick] at (-1+0.4,9){$j$};

			%		\draw (-0.3,9) --(0.3,9);
		\end{tikzpicture}
	\end{center}
	\caption{The third peak at level $j$.}
\end{figure}	

We don't need too many new computations, as we can modify the previous results. If one adds an arbitrary non-empty number
of up-steps after the $m$-th valley, one has reached the $(m+1)$-st peak! This is basically a substitution!
Start from
\begin{equation*}
	\Phi(u)=\frac{s_2(1-zu)}{z(1-us_2)}
\end{equation*}
and attach a sequence of up-steps: $u^j\to \frac{zu}{1-zu}u^j$. A factor $w$ is also important, since the $m$-th valley corresponds to the
$(m+1)$-st peak.
Now
\begin{equation*}
	\frac{zuw}{1-zu}\frac{s_2(1-zu)}{z(1-us_2)}=\frac{us_2w}{1-us_2}=w\sum_{j\ge1}u^js_2^j.
\end{equation*}
The computation
\begin{equation*}
	w\sum_{j\ge0}js_2^j\cdot s_2^j\Big|_{w=1}
\end{equation*}
was basically done before, and the local expansion leads to
\begin{equation*}
	\frac{2w}{1-w}-\frac{2w\sqrt{(1-w)(9-w)}}{(1-w)^2}\sqrt{1-4z},
\end{equation*}
and the generating function of the average values of the $m$-th peak is
\begin{equation*}
	\mathsf{Peak}(w)=\frac{w\sqrt{(1-w)(9-w)}}{(1-w)^2}.
\end{equation*}
A local expansion of this results in
\begin{equation*}
	\mathsf{Peak}(w)\sim{\frac {2\sqrt {2}}{ \left( 1-w \right) ^{3/2}}}-{\frac {15}{8}}{\frac {\sqrt {2}}{\sqrt {1-w}}}.
\end{equation*}
Taking differences:
\begin{equation*}
	\mathsf{Peak}(w)-\mathsf{Valley}(w)\sim\frac{2}{1-w}-{\frac {\sqrt {2}}{\sqrt {1-w}}},
\end{equation*}
and translating into asymptotics:
\begin{equation*}
	2-\frac{\sqrt2}{\sqrt{\pi m}}.
\end{equation*}
The formula $2+O(m^{-1/2})$ was already known to Kemp \cite{Kemp-oscillations}. As Kemp stated in
\cite{Kemp-oscillations}, which was confirmed in  \cite{KP-hyper}, the generating functions $\mathsf{Peak}(w)$ and $\mathsf{Valley}(w)$
can be expressed by Legendre polynomials at special values. This is a bit artificial and not too useful in itself.

\section{Deutsch-paths in a strip}

Emeric Deutsch \cite{Deutsch} had the idea to consider a variation of ordinary Dyck paths, by augmenting the usual up-steps and down-steps by one unit each,
by down-steps of size $3,5,7,\dots$. This leads to ternary equations, as can be seen for instance from \cite{Deutsch-ternary}.

The present author started to investigate a related but simpler model of down-steps $1,2,3,4,\dots$ and investigated it (named Deutsch paths in honour of Emeric Deutsch)
in a series of papers, \cite{Deutsch1,Deutsch-slice,Prodinger-fibo}.

This section is a  further member of this series: The condition that (as with Dyck paths) the paths cannot enter negative territory, is relaxed, by introducing
a negative boundary $-t$. Here are two recent publications about such a negative boundary: \cite{Selkirk-master} and \cite{jcmcc}.

Instead of allowing negative altitudes, we think about the whole system shifted up by $t$ units, and start at the point $(0,t)$ instead. This
is much better for the generating functions that we are going to investigate. Eventually, the results can be re-interpreted as results about 
enumerations with respect to a negative boundary.

The setting with flexible initial level $t$ and final level $j$ allows us to consider the Deutsch paths also from  right to left (they are not symmetric!),
without any new computations.

The next sections achieves this, using the celebrated kernel-method,  one of the tools that is dear to our heart \cite{Prodinger-kernel}.

In the following section, an additional upper bound is introduced, so that the Deutsch paths live now in a strip. The way to attack this
is linear algebra. Once everything has been computated, one can relax the conditions and let lower/upper boundary go to $\mp\infty$.

\subsection*{Generating functions and the kernel method}

As discussed, we consider Deutsch paths starting at $(0,t)$ and ending at $(n,j)$, for $n,t,j\ge0$. First we consider univariate generating functions
$f_j(z)$, where $z^n$ stays for $n$ steps done, and $j$ is the final destination. The recursion is immediate:
\begin{equation*}
	f_j(z)=[\![t=j]\!]+zf_{j-1}(z)+z\sum_{k>j}f_k(z),
\end{equation*}
where $f_{-1}(z)=0$. Next, we consider
\begin{equation*}
	F(z,u):=\sum_{j\ge0}f_j(z)u^j,
\end{equation*}
and get
\begin{align*}
	F(z,u)&=u^t+zuF(z,u)+z\sum_{j\ge0}u^j\sum_{k>j}f_k(z)\\
	&=u^t+zuF(z,u)+z\sum_{k>0}f_k(z)\sum_{0\le j<k}u^j\\
	&=u^t+zuF(z,u)+z\sum_{k\ge0}f_k(z)\frac{1-u^k}{1-u}\\
	&=u^t+zuF(z,u)+\frac z{1-u}[ F(z,1)-F(z,u)]\\
	&=\frac{u^t(1-u)+zF(z,1)}{z-zu+zu^2+1-u}.
\end{align*}
Since the critical value is around $u=1$, we write the denominator as
\begin{equation*}
	z(u-1)^2+(u-1)(z-1)+z=z(u-1-r_1)(u-1-r_2),
\end{equation*}
with
\begin{align*}
	r_1=\frac{1-z+\sqrt{1-2z-3z^2}}{2z},\quad r_2=\frac{1-z-\sqrt{1-2z-3z^2}}{2z}.
\end{align*}
The factor $(u-1-r_2)$ is bad, so the numerator must vanish for $[u^t(1-u)+zF(z,1)]|_{u=1+r_2}$, therefore
\begin{equation*}
	zF(z,1)=(1+r_2)^tr_2.
\end{equation*} 
Furthermore
\begin{equation*}
	F(z,u)=
	\frac{\frac{u^t(1-u)+zF(z,1)}{u-r_2}}{z(u-r_1)}.
\end{equation*}
The expressions become prettier using the substitution $z=\frac{v}{1+v+v^2}$; then
\begin{equation*}
	r_1=\frac{1}{v},\quad r_2=v.
\end{equation*}
It can be proved by induction (or computer algebra) that 
\begin{equation*}
	\frac{u^t(1-u)+v(1+v)^t}{u-1-v}=-v\sum_{k=0}^{t-1}(1+v)^{t-1-k}-u^t.
\end{equation*}
Furthermore
\begin{equation*}
	\frac1{z(u-1-r_1)}=-\frac1{z(1+r_1)(1-\frac{u}{1+r_1})},
\end{equation*}
and so 
\begin{equation*}
	f_j(z)=[u^j]F(z,u)=[u^j]\biggl[v\sum_{k=0}^{t-1}(1+v)^{t-1-k}u^k+u^t\biggr]\sum_{\ell\ge0}\frac{u^{\ell}}{z(1+r_1)^{\ell+1}}.
\end{equation*}
Of interest are two special cases: 
The case that was studied before \cite{Deutsch1} is $t=0$:
\begin{equation*}
	f_j=\frac{(1+v+v^2)v^{j}}{(1+v)^{j+1}}.
\end{equation*}
The other special case is $j=0$ for general $t$, as it may be interpreted as Deutsch paths read from right to left, starting at level $0$ and 
ending at level $t\ge1$ (for $t=0$, the previous formula applies):
\begin{align*}
	f_0(z)&=[u^0]\biggl[v\sum_{k=0}^{t-1}(1+v)^{t-1-k}u^k+u^t\biggr]\sum_{\ell\ge0}\frac{u^{\ell}}{z(1+r_1)^{\ell+1}}\\
	&=v(1+v)^{t-1}\frac{1}{z(1+r_1)}=v(1+v+v^2)(1+v)^{t-2}.
\end{align*}

The next section will present a simplification of the expression for $f_j(z)$, which could be obtained directly by distinguishing cases and summing some geometric series.

\subsection*{Refined analysis: lower and upper boundary}

Now we consider Deutsch paths bounded from below by zero and bounded from above by $m-1$; they start at level $t$ and end at level $j$ after $n$ steps.
For that, we use generating functions $\varphi_j(z)$ (the quantity $t$ is a silent parameter here). The recursions that are straight-forwarded are best organized in a
matrix, as the following example shows.
\begin{equation*}
	\left(\begin{matrix}
		1&-z&-z&-z&-z&-z&-z&-z\\
		-z&	1&-z&-z&-z&-z&-z&-z\\
		0&	-z&	1&-z&-z&-z&-z&-z\\
		0& 0&	-z&	1&-z&-z&-z&-z\\
		0& 0& 0&	-z&	1&-z&-z&-z\\
		0& 0& 0&0&	-z&	1&-z&-z\\
		0& 0& 0&0&0&	-z&	1&-z\\
		0& 0& 0&0&0&0&	-z&	1\\
	\end{matrix}\right)
	\left(\begin{matrix}
		\varphi_0\\
		\varphi_1\\
		\varphi_2\\
		\varphi_3\\
		\varphi_4\\
		\varphi_5\\
		\varphi_6\\
		\varphi_7\\
	\end{matrix}\right)=
	\left(\begin{matrix}
		0\\
		0\\
		0\\
		1\\
		0\\
		0\\
		0\\
		0\\
	\end{matrix}\right)
	\begin{tikzpicture} 
		\draw [](0,0)--(0,0);
		\node at (0.5,1.3){\Bigg\}$t$};
	\end{tikzpicture}
\end{equation*}	

The goal is now to solve this system. For that the substitution $z=\frac{v}{1+v+v^2}$ is used throughout. The method is to use Cramer's rule, which means that the right-hand side has to replace various columns of the matrix, and determinants have to be computed. At the end, one has to divide by the determinant of the system.

Let $D_m$ be the determinant of the matrix with $m$ rows and columns. 
The recursion
\begin{equation*}
	(1+v+v^2)^2m_{n+2}-(1+v+v^2)(1+v)^2D_{m+1}+v(1+v)^2D_{m}=0
\end{equation*}	
appeared already in \cite{Deutsch1} and is not difficult to derive and to solve:
\begin{equation*}
	D_m=\frac{(1+v)^{m-1}}{(1+v+v^2)^m}\frac{1-v^{m+2}}{1-v}.
\end{equation*}

To solve the system with Cramer's rule, we must compute a determinant of the following type,
\begin{center}\small
	\begin{tikzpicture}
		[scale=0.4]
		\draw (0,0)--(5,0)--(5,-3)--(0,-3)--(0,0);
		\node at (2.5,0.9){$j$};
		\draw[<->](0,0.5)--(5,0.5);
		\draw[<->](-0.5,0)--(-0.5,-3);
		\node at (-0.9,-1.5){$t$};
		\node at (5.5,-3.5){$\boldsymbol{1}$};	
		\newcommand\x{6};\draw (0+\x,0)--(8+\x,0)--(8+\x,-3)--(0+\x,-3)--(0+\x,0);	
		\newcommand\y{-4};				\draw (0,0+\y)--(5,0+\y)--(5,-7+\y)--(0,-7+\y)--(0,0+\y);
		\draw (0+\x,0+\y)--(8+\x,0+\y)--(8+\x,-7+\y)--(0+\x,-7+\y)--(0+\x,0+\y);
		\draw[<->](0,-11.5)--(14,-11.5);
		\draw[<->](14.5,0)--(14.5,-11);
		\node at (7,-11.9){$m$};
		\node at (15.0,-5.5){$m$};
		\foreach \x in {0,1,2,3,4}
		{
			\node at (5.5,-\x/1.5){$\tiny\boldsymbol{0}$};
		}
		
		\foreach \x in {6.5,7.5,8.5,9.5,10.5,11.5,12.5,13.5,14.5,15.5,16.5}
		{
			\node at (5.5,-\x/1.5){$\tiny\boldsymbol{0}$};
		}
		
	\end{tikzpicture}
\end{center}
where the various rows are replaced by the right-hand side. While it is not impossible to solve this recursion by hand,
it is very easy to make mistakes, so it is best to employ a computer. Let $D(m;t,j)$ the determinant according to the
drawing.

It is not unexpected that the results are different for $j<t$ resp.\ $j\ge t$. Here is what we found:
\begin{equation*}
	D(m;t,j)=\frac{(1+v)^{t-j-3+m}(1-v^{j+1})v(1-v^{m-t})}{(1-v)^2(1+v+v^2)^{m-1}},\quad\text{for}
	\ j<t,
\end{equation*}
\begin{equation*}
	D(m;t,j)=\frac{v^{j-t}(1-v^{t+2})(1-v^{1-j+m})}{(1-v)^2(1+v+v^2)^{m-1}(1+v)^{j-t+3-m}},\quad\text{for}
	\ j\ge t.
\end{equation*}
To solve the system, we have to divide by the determinant $D_m$, with the result
\begin{equation*}
	\varphi_j=	\frac{D(m;t,j)}{D_m}=
	\frac{(1+v)^{t-j-2}(1-v^{j+1})v(1-v^{m-t})(1+v+v^2)}{(1-v)(1-v^{m+2})}
	,\quad\text{for}
	\ j<t,
\end{equation*}
\begin{equation*}
	\varphi_j=	\frac{D(m;t,j)}{D_m}=\frac{v^{j-t}(1-v^{t+2})(1-v^{1-j+m})(1+v+v^2)}{(1-v)(1+v)^{j-t+2}(1-v^{m+2})}
	,\quad\text{for}
	\ j\ge t.
\end{equation*}
We found all this using Computer algebra. Some critical minds may argue that this is only experimental. One way of rectifying this
would be to show that indeed the functions $\varphi_j$ solve the system, which consists of summing various geometric series; again,
a computer could be helpful for such an enterprise.

Of interest are also the limits for $m\to\infty$, i.e., no upper boundary:
\begin{equation*}
	\varphi_j=\lim_{m\to\infty}\frac{D(m;t,j)}{D_m}=
	\frac{(1+v)^{t-j-2}(1-v^{j+1})v(1+v+v^2)}{(1-v)}
	,\quad\text{for}
	\ j<t,
\end{equation*}
\begin{equation*}
	\varphi_j=\frac{v^{j-t}(1-v^{t+2})(1+v+v^2)}{(1-v)(1+v)^{j-t+2}	}
	,\quad\text{for}
	\ j\ge t.
\end{equation*}
The special case $t=0$ appeared already in the previous section:
\begin{equation*}
	\varphi_j=\frac{v^{j}(1+v+v^2)}{(1+v)^{j+1}	}.
\end{equation*}
Likewise, for $t\ge1$,
\begin{equation*}
	\varphi_0=v(1+v+v^2)(1+v)^{t-2}.
\end{equation*}
In particular, the formul\ae\ show that the expression from the previous section can be simplified in general, which
could have been seen directly, of course.

\begin{theorem}
	The generating function of Deutsch path with lower boundary 0 and upper boundary $m-1$, starting at $(0,t)$ and ending at $(n,j)$ is given by
	\begin{gather*}
		\frac{(1+v)^{t-j-2}(1-v^{j+1})v(1-v^{m-t})(1+v+v^2)}{(1-v)(1-v^{m+2})}
		,\quad\text{for}
		\ j<t,\\
		\frac{v^{j-t}(1-v^{t+2})(1-v^{1-j+m})(1+v+v^2)}{(1-v)(1+v)^{j-t+2}(1-v^{m+2})}
		,\quad\text{for}
		\ j\ge t,		
	\end{gather*}
	with the substitution $z=\dfrac{v}{1+v+v^2}$.
\end{theorem}
By shifting everything down, we can interpret the results as Deutsch walks between boundaries $-t$ and $m-1-t$, starting at the origin $(0,0)$ and ending at $(n,j-t)$.

\begin{theorem}
	The generating function of Deutsch path with lower boundary $-t$ and upper boundary $h$, starting at $(0,0)$ and ending at $(n,i)$ with $-t\le i\le h$ is given by
	\begin{gather*}
		\frac{(1+v)^{i-2}(1-v^{i+t+1})v(1-v^{h+1})(1+v+v^2)}{(1-v)(1-v^{h+t+3})}
		,\quad\text{for}
		\ i<0,\\	
		\frac{v^{i}(1-v^{t+2})(1-v^{2-i+h})(1+v+v^2)}{(1-v)(1+v)^{i+2}(1-v^{h+t+3})}
		,\quad\text{for}
		\ i\ge 0.		
	\end{gather*}
\end{theorem}
It is possible to consider the limits $t\to\infty$ and/or $h\to\infty$ resulting in simplified formul\ae.
\begin{theorem}
	The generating function of Deutsch path with lower boundary $-t$ and upper boundary $\infty$, starting at $(0,0)$ and ending at $(n,i)$ with $-t\le i$ is given by
	\begin{gather*}
		\frac{(1+v)^{i-2}(1-v^{i+t+1})v(1+v+v^2)}{(1-v)}
		,\quad\text{for}
		\ i<0,\\	
		\frac{v^{i}(1-v^{t+2})(1+v+v^2)}{(1-v)(1+v)^{i+2}}
		,\quad\text{for}
		\ i\ge 0.		
	\end{gather*}
\end{theorem}
\begin{theorem}
	The generating function of Deutsch path with lower boundary $-\infty$ and upper boundary $h$, starting at $(0,0)$ and ending at $(n,i)$ with $\le i\le h$ is given by
	\begin{gather*}
		\frac{(1+v)^{i-2}v(1-v^{h+1})(1+v+v^2)}{(1-v)}
		,\quad\text{for}
		\ i<0,\\	
		\frac{v^{i}(1-v^{2-i+h})(1+v+v^2)}{(1-v)(1+v)^{i+2}}
		,\quad\text{for}
		\ i\ge 0.		
	\end{gather*}
\end{theorem}
\begin{theorem}
	The generating function of unbounded Deutsch path  starting at $(0,0)$ and ending at $(n,i)$ is given by
	\begin{gather*}
		\frac{(1+v)^{i-2}v(1+v+v^2)}{(1-v)}
		,\quad\text{for}
		\ i<0,\\	
		\frac{v^{i}(1+v+v^2)}{(1-v)(1+v)^{i+2}}
		,\quad\text{for}
		\ i\ge 0.		
	\end{gather*}
\end{theorem}

\section{Publication status of the problems treated in this walk in the garden}
\begin{itemize}
	\item The section on Hoppy walks is new; an earlier version appeared as\\
	arXiv:2009.13474
	
	\item The section on Combinatorics of sequence A002212 is new; an earlier version appeared as
	arXiv:2106.14782

\item	The section on Binary trees and Horton-Strahler numbers is classical; it was included since it can
be used almost verbatim for unary-binary trees when employing the proper substitution

\item The section on marked ordered trees is new and basically included in\\ 	arXiv:2106.14782

\item The bijection between multi-edge trees and 3-coloured Motzkin paths is new; early version in
arXiv:2105.03350

\item The section on combinatorics of skew Dyck paths was submitted to journal, and no answer was received; 
current version in arXiv:2108.09785

\item The short section on ternary trees was, after a waiting period of 15 months, accepted by a journal in
december 2021.

\item The short analysis of Retakh's Motzkin paths appeared in ECA.

\item The amplitude of Motzkin paths was submitted to a journal; no answer. 

\item The section on oscillations in honour of Rainer Kemp was written for this survey.

\item The enumeration of Deutsch-paths in a strip was exclusively written for this survey article;
 arXiv:2108.12797

\end{itemize}

\clearpage

\bibliographystyle{plain}

%\bibliography{hex}

\end{document}